\documentclass[11pt]{article}
\usepackage{color}
\usepackage{latexsym}
\usepackage{amssymb}
\usepackage{amsmath,amsfonts,theorem,euscript,array,enumerate,mathrsfs}
\usepackage{graphicx}
\newtheorem{Theorem}{Theorem}[section]
\newtheorem{Definition}[Theorem]{Definition}
\newtheorem{Proposition}[Theorem]{Proposition}

\newtheorem{Lemma}[Theorem]{Lemma}
\newtheorem{Corollary}[Theorem]{Corollary}
\newtheorem{Remark}[Theorem]{Remark}

\numberwithin{equation}{section}

\def\esssup_#1{\underset{#1}{\mathrm{ess\,sup\, }}}
\def\essinf_#1{\underset{#1}{\mathrm{ess\,inf\, }}}

\def \trans{^{\scriptscriptstyle{\intercal}}}

\def \trans{^{\scriptscriptstyle{\intercal }}}

\def \A{\mathbb{A}}

\def \N{\mathbb{N}}
\def \R{\mathbb{R}}

\def \E{\mathbb{E}}
\def \F{\mathbb{F}}

\def \H{\mathbb{H}}
\def \L{\mathbb{L}}
\def \P{\mathbb{P}}
\def \Q{\mathbb{Q}}

\def \S{\mathbb{S}}

\def \X{\mathbb{X}}

\def \Bc{{\cal B}}

\def \Fc{{\cal F}}

\def \Kc{{\cal K}}

\def \Uc{{\cal U}}

\def \eps{\varepsilon}

\def \ep{\hbox{ }\hfill$\Box$}

\allowdisplaybreaks

\begin{document}

\title{\small{\textbf{STRONG-VISCOSITY SOLUTIONS: CLASSICAL AND 
PATH-DEPENDENT PDEs}}}

\author{
\textsc{Andrea Cosso}\thanks{University of Bologna, Dipartimento di Matematica, Piazza di Porta S. Donato, 5, 40126 Bologna, Italy, \sf andrea.cosso at unibo.it}
\qquad\quad
\textsc{Francesco Russo}\thanks{ENSTA ParisTech, Universit\'e Paris-Saclay, Unit\'e de Math\'ematiques Appliqu\'ees, 828, boulevard des Mar\'echaux, 91120 Palaiseau, France, \sf francesco.russo at ensta-paristech.fr}
}

\date{December 7th 2017}

\maketitle

\begin{abstract}
\noindent The aim of the present work is the introduction of a
 viscosity type solution, called \emph{strong-viscosity solution}
emphasizing also a similarity with the existing notion of \emph{strong solution}
in the literature. It has  the following peculiarities: it is a purely analytic object; it can be easily adapted to more general
 equations than classical partial differential equations. 
First, we introduce the notion of strong-viscosity solution for semilinear parabolic partial differential equations, defining it, in a few words, as the pointwise limit of classical solutions to perturbed semilinear parabolic partial differential equations; we compare it with the standard definition of viscosity solution. Afterwards, we extend the concept of strong-viscosity solution to the case of semilinear parabolic path-dependent partial differential equations, providing an existence and uniqueness result.

\vspace{4mm}


\vspace{4mm}

\end{abstract}

\section{Introduction}

As it is well-known, viscosity solutions represent a cornerstone in the theory of Partial Differential Equations (PDEs) and their range of application is enormous, see the user's guide \cite{crandishiilions92}. Here, we just emphasize the important role they played in the study of semilinear parabolic
 partial differential equations. We also emphasize the role of  Backward Stochastic Differential Equations 
(BSDEs), which constitute a probabilistic counterpart of viscosity solutions of 
 semilinear parabolic partial differential equation, see the seminal paper \cite{pardoux_peng92}.

\vspace{1mm}

The aim of the present work is the definition of a variant of viscosity type solution,
called \emph{strong-viscosity solution} 
to distinguish it from the classical one. Compared to this latter, for several aspects it seems easier to handle and it can be easily adapted to a large class of equations.



\vspace{1mm}

In recent years, there has been a renewed interest in the study of generalized partial differential equations, motivated by the study of Markovian stochastic control problems with state variable living in an infinite dimensional space (see \cite{daprato_zabczyk14}) or path-dependent problems, for example, stochastic control problems with delay, see \cite{fabbrigozziswiech14}. The theory of backward stochastic differential equations is flexible enough to be extended to deal with both problems, see, e.g., \cite{fuhrman_pham13}, \cite{fuhrman_tessitore02}, \cite{peng_wang}. From an analytic point of view, regarding infinite dimensional Markovian problems, there exists in general a corresponding partial differential equation in infinite dimension, and also the notion of viscosity solution has been extended to deal with this case, see \cite{crandall_lions91}, \cite{swiech94}, and \cite{fabbrigozziswiech14}. However, uniqueness for viscosity solutions revealed to be arduous to extend to the infinite dimensional setting and requires, in general, strong assumptions on the coefficients of the partial differential equation. 

\vspace{1mm}

Concerning path-dependent problems, it is still not clear what should be the corresponding analytic characterization in terms of partial 
differential equations, whose probabilistic counterpart is represented by the backward stochastic differential equation. A possible solution to 
this problem is represented by the class of equations introduced in Chapter 9 of \cite{DGR} within the framework of Banach space valued calculus, 
for which we refer also to \cite{flandoli_zanco13}. Alternatively, \cite{dupire} introduced the concept of Path-dependent Partial Differential 
Equation (PPDE), which could do the job. Even if it is still not completely definite in the literature what a path-dependent partial differential 
equation is (indeed, it mainly depends on the definition of functional derivatives adopted), the issue of providing a suitable definition of 
viscosity solution for path-dependent partial differential equations has already attracted a great interest, see for example 
\cite{ektz,etzI,etzII,rtz1,tangzhang13}, motivated by the fact that regular solutions to path-dependent PDEs in general
 exist only under strong assumptions, see Remark \ref{R:MotivationSV}. We drive the attention in particular to the definition of viscosity
 solution to path-dependent PDEs provided by \cite{ektz,etzI,etzII,CFGRT,rtz1}, where the classical minimum/maximum property, appearing in the standard 
definition of viscosity solution, is replaced with an optimal stopping problem under nonlinear expectation~\cite{etzOptStop}. 
Notice that probability plays an essential role in this latter definition, which can, more properly, be interpreted as a probabilistic version of 
the standard definition of viscosity solution, rather than a purely analytic object; indeed, quite interestingly, the proof of the comparison 
principle turns out to be nearly a ``translation'' into probabilistic terms of the classical proof of the comparison principle, see \cite{rtz1}. 
We also emphasize that a similar notion of solution, called stochastic weak solution, has been introduced in the recent paper 
\cite{leao_ohashi_simas14} in the context of variational inequalities for the Snell envelope associated to a non-Markovian continuous process $X$.
 Those authors also revisit functional It\^o
calculus, making use of stopping times. This approach seems very promising.

\vspace{1mm}
The paper is organized as follows. First, in Section \ref{S:SV_Markov}, we develop the theory of strong-viscosity solutions in the finite dimensional
Markovian case, applying it to semilinear parabolic partial differential equations. A strong-viscosity supersolution (resp. subsolution) is defined,
in a few words, as the pointwise limit of classical supersolutions (resp. subsolutions) to perturbed semilinear parabolic PDEs. 
A generalized strong-viscosity solution is both a strong-viscosity supersolution and a strong-viscosity subsolution. This definition is more in the
 spirit of the standard definition of viscosity solution. We also introduce another definition, simply called strong-viscosity solution, which is 
defined as the pointwise limit of classical solutions to perturbed semilinear parabolic PDEs. We notice that the definition of strong-viscosity 
solution is similar in spirit to the vanishing 
viscosity method, which represents one of the primitive ideas leading to the conception of the modern definition of viscosity solution and
 justifies the term \emph{viscosity} in the name, which is also justified by the fact that a strong-viscosity solution is not assumed to be 
differentiable.
 Our definition is likewise inspired by the notion of strong solution (which explains the presence of the term \emph{strong} in the name), as defined for example in \cite{cerrai01}, \cite{gozzi_russo06a}, and \cite{gozzi_russo06b}, even though strong solutions are in general required to be more regular than viscosity type solutions. Finally, we observe that the notion of strong-viscosity solution has also some similarities with the concept of 
\emph{good solution}, which turned out to be equivalent to the definition of $L^p$-viscosity solution for certain fully nonlinear partial differential equations, see, e.g., \cite{cerutti_escauriaza_fabes93}, \cite{crandall_kocan_soravia_swiech94}, \cite{jensen96}, and \cite{jensen_kocan_swiech02}. 

\vspace{1mm}

We prove in Section \ref{S:SV_Markov}, Theorem \ref{T:RepresentationSuperSub}, that every strong-viscosity supersolution (resp. subsolution) can be represented in terms of a supersolution (resp. subsolution) to a backward stochastic differential equation. This in turn implies that a comparison principle (Corollary \ref{C:CompThm}) for strong-viscosity sub and supersolutions holds and follows from the comparison theorem for backward stochastic differential equations. In particular, the proof of the comparison principle is probabilistic and easier to extend to different contexts than the corresponding analytical proof
for classical viscosity solutions, which is based on real analysis' tools as Ishii's lemma and the doubling of variables technique. We conclude Section \ref{S:SV_Markov} providing two existence results (Theorem \ref{T:ExistSV_Markov} and Theorem \ref{T:ExistSV_Markov2}) for strong-viscosity solutions under quite general assumptions.

\vspace{1mm}

In Section \ref{S:SVPath} we extend the notion of strong-viscosity solution to the case of semilinear parabolic path-dependent partial differential equations, leaving to future research other possible extensions, e.g., the case of partial differential equations in infinite dimension. For PPDEs, as already said, a viscosity type solution, meant as a purely analytic object, is still missing, so we try to fill the gap. As previously noticed, the concept of path-dependent partial differential equation is still not definite in the literature and, in the present paper, we adopt the setting developed in the companion paper \cite{cosso_russo15a}. However, we notice that, if we had worked with the definition of functional derivatives and path-dependent partial differential equation used, e.g., in \cite{dupire,contfournie13}, the same results would hold 
in that context without any change, but for notational ones, 
see \cite{cosso_russo15a} for some insights on the link between these different settings. Let us explain the reasons why we adopt the definitions of \cite{cosso_russo15a}. First, in \cite{cosso_russo15a} the time and space variables $(t,\eta)\in[0,T]\times C([-T,0])$ play two distinct roles; moreover the space variable $\eta$ (i.e., the path) always represents the past trajectory of the process. This is in complete accordance with the literature on stochastic control problems with delay (see, e.g., \cite{cho} and \cite{fabbrigozziswiech14}), which is, for us, one of the main applications of path-dependent partial differential equations. On the contrary, in \cite{contfournie13} the time and space variables are strictly related to each other; moreover, the path represents the entire trajectory (past, present, and future) of the process, so that the notion of non-anticipative functional is required, see Definition 2.1 in \cite{contfournie13}.

\vspace{1mm}

We prove in Section \ref{S:SVPath}, Theorem \ref{T:UniqSV}, a uniqueness result for strong-viscosity solutions to path-dependent PDEs proceeding as in the finite dimensional Markovian case, i.e., by means of probabilistic methods based on the theory of backward stochastic differential equations. We also prove an existence result (Theorem \ref{T:ExistSV}) for strong-viscosity solutions in a more restrictive framework, which is based on the idea that a candidate solution to the path-dependent PDE is deduced from the corresponding backward stochastic differential equation. The existence proof consists in building a sequence of strict solutions (we prefer to use the term \emph{strict} in place of \emph{classical}, because even the notion of smooth solution can not be considered classical for path-dependent partial differential equations; indeed, all the theory is very recent) to perturbed path-dependent PDEs converging to our strong-viscosity solution. This regularization procedure is performed having in mind the following simple property: when the coefficients of the path-dependent partial differential equation are smooth enough the solution is smooth as well, i.e., the solution is strict. In the path-dependent case, smooth coefficients means \emph{cylindrical coefficients}, i.e., smooth maps of integrals of regular functions with respect to the path, as in the statement of Theorem \ref{T:ExistenceStrict}.

\vspace{1mm}

Finally, we defer some technical results to the Appendix. More precisely, we prove some basic estimates for path-dependent stochastic differential equations in Lemma \ref{L:AppendixX}. Then, we state a standard (but, to our knowledge, not at disposal in the literature) estimate for supersolutions to non-Markovian backward stochastic differential equations, see Proposition \ref{P:EstimateBSDEAppendix}. Afterwards, we prove the limit Theorem \ref{T:LimitThmBSDE} for supersolutions to backward stochastic differential equations. We conclude the Appendix with a technical result, Lemma \ref{L:StabilityApp}, of real analysis.

\section{Strong-viscosity solutions in the Markovian case}
\label{S:SV_Markov}

In the present section we introduce the notion of strong-viscosity solution in the non-path-depen\-dent case, for the  semilinear parabolic PDE
\begin{equation}
\label{KolmEq_Markov}
\begin{cases}
\partial_t u(t,x) + \langle b(t,x),D_x u(t,x)\rangle + \frac{1}{2}\textup{tr}(\sigma\sigma\trans(t,x)D_x^2 u(t,x)) & \\
\hspace{2.8cm}+\, f(t,x,u(t,x),\sigma\trans(t,x)D_x u(t,x)) \ = \ 0, &\forall\,(t,x)\in[0,T[\times\R^d, \\
u(T,x) \ = \ h(x), &\forall\,x\in\R^d,
\end{cases}
\end{equation}
where $b\colon[0,T]\times\R^d\rightarrow\R^d$, $\sigma\colon[0,T]\times\R^d\rightarrow\R^{d\times d}$, $f\colon[0,T]\times\R^d\times\R\times\R^d\rightarrow\R$, and $h\colon\R^d\rightarrow\R$ satisfy the following assumptions. 

\vspace{3mm}

\noindent\textbf{(A0)} \hspace{3mm} $b$, $\sigma$, $f$, $h$ are Borel measurable functions satisfying, for some positive constants $C$ and $m$,
\begin{align*}
|b(t,x)-b(t,x')| + |\sigma(t,x)-\sigma(t,x')| \ &\leq \ C|x-x'|, \\
|f(t,x,y,z)-f(t,x,y',z')| \ &\leq \ C\big(|y-y'| + |z-z'|\big), \\
|b(t,0)| + |\sigma(t,0)| \ &\leq \ C, \\
|f(t,x,0,0)| + |h(x)| \ &\leq \ C\big(1 + |x|^m\big),
\end{align*}
for all $t\in[0,T]$, $x,x'\in\R^d$, $y,y'\in\R$, and $z,z'\in\R^d$.

\subsection{Notations}
\label{SubS:Notation}

We denote by $\R^{d\times d}$ the linear space of real matrices of order $d$. In the all paper, $|\cdot|$ denotes the absolute value of a real number or the usual Euclidean norm in $\R^d$ or the Frobenius norm in $\R^{d\times d}$.

\vspace{1mm}

We fix a complete probability space $(\Omega,\Fc,\P)$ on which a $d$-dimensional Brownian motion $W=(W_t)_{t\geq0}$ is defined. Let $\F=(\Fc_t)_{t\geq0}$ denote the completion of the natural filtration generated by $W$. We introduce the following spaces of stochastic processes.

\begin{itemize}
\item $\S^p(t,T)$, $p\geq1$, $t \leq T$, the set  of real c\`adl\`ag adapted processes $Y=(Y_s)_{t\leq s\leq T}$ such that
\[
\|Y\|_{_{\S^p(t,T)}}^p := \ \E\Big[ \sup_{t\leq s\leq T} |Y_s|^p \Big] \ < \ \infty.
\]
\item $\H^p(t,T)^d$, $p$ $\geq$ $1$, $t \leq T$, the set of  $\R^d$-valued predictable processes $Z=(Z_s)_{t\leq s\leq T}$ such that
\[
\|Z\|_{_{\H^p(t,T)^d}}^p := \ \E\bigg[\bigg(\int_t^T |Z_s|^2 ds\bigg)^{\frac{p}{2}}\bigg] \ < \ \infty.
\]
We simply write $\H^p(t,T)$ when $d=1$.
\item $\A^{+,2}(t,T)$, $t \leq T$, the  set of real nondecreasing predictable processes $K$ $=$ $(K_s)_{t\leq s\leq T}\in\S^2(t,T)$ with $K_t$ $=$ $0$, so that
\[
\|K\|_{_{\S^2(t,T)}}^2 = \ \E\big[K_T^2\big].
\]
\item $\L^p(t,T;\R^{d'})$, $p\geq1$, $t \leq T$, the set of $\R^{d'}$-valued adapted processes $\varphi = (\varphi_s)_{t \leq s \leq T}$ such that
\[
\|\varphi\|_{_{\L^p(t,T;\R^{d'})}}^p := \ \E\bigg[\int_t^T |\varphi_s|^p ds\bigg] \ < \ \infty.
\]
\end{itemize}

We also consider, for every $(t,x)\in[0,T]\times\R^d$, the stochastic differential equation
\begin{equation}
\label{SDE_Markov}
\begin{cases}
dX_s = \ b(s,X_s)dt + \sigma(s,X_s)dW_s, \qquad\qquad s\in[t,T], \\
X_t \ = \ x.
\end{cases}
\end{equation}

It is well-known (see, e.g., Theorem 14.23 in \cite{jacod79}) that, under Assumption {\bf (A0)}, there exists a unique $($up to indistinguishability$)$ $\F$-adapted continuous process $X^{t,x}=(X_s^{t,x})_{s\in[t,T]}$ strong solution to equation \eqref{SDE_Markov}.

\subsection{First definition of strong-viscosity solution}

We begin recalling the standard definition of classical solution.

\begin{Definition}
A function $u\colon[0,T]\times\R^d\rightarrow\R$ is called  \textbf{classical solution} to equation \eqref{KolmEq_Markov} if $u\in C^{1,2}([0,T[\times\R^d)\cap C([0,T]\times\R^d)$ and solves \eqref{KolmEq_Markov}.
\end{Definition}

\noindent We state a uniqueness result for classical solutions.

\begin{Proposition}
\label{P:UniqStrictFiniteStandard}
Suppose that Assumption {\bf (A0)} holds. Let $u\colon[0,T]\times\R^d\rightarrow\R$ be a classical solution to equation \eqref{KolmEq_Markov}, satisfying the polynomial growth condition
\begin{equation}
\label{PolGrowthCondMarkov}
|u(t,x)| \ \leq \ C'\big(1 + |x|^{m'}\big), \qquad \forall\,(t,x)\in[0,T]\times\R^d,
\end{equation}
for some positive constants $C'$ and $m'$. Then, the following Feynman-Kac
 formula holds:
\begin{equation}
\label{Identification2}
u(t,x) \ = \ Y_t^{t,x}, \qquad \forall\,(t,x)\in[0,T]\times\R^d,
\end{equation}
where $(Y_s^{t,x},Z_s^{t,x})_{s\in[t,T]}=(u(s,X_s^{t,x}),\sigma\trans(s,X_s^{t,x})D_x u(s,X_s^{t,x})1_{[t,T[}(s))_{s\in[t,T]}\in\S^2(t,T)\times\H^2(t,T)^d$ is the unique solution to the backward stochastic differential equation: $\P$-a.s.,
\begin{equation}
\label{BSDE_Markov2}
Y_s^{t,x} \ = \ h(X_T^{t,x}) + \int_s^T f(r,X_r^{t,x},Y_r^{t,x},Z_r^{t,x}) dr - \int_s^T Z_r^{t,x} dW_r, \qquad t \leq s \leq T.
\end{equation}
In particular, there exists at most one classical solution to equation \eqref{KolmEq_Markov} satisfying a polynomial growth condition as in \eqref{PolGrowthCondMarkov}.
\end{Proposition}
\textbf{Proof.}
The proof is standard, even if we have not found an exact reference for it in the literature. We just give the main ideas. Fix $(t,x)\in[0,T[\times\R^d$ and set, for all $t \leq s \leq T$,
\[
Y_s^{t,x} \ = \ u(s,X_s^{t,x}), \qquad Z_s^{t,x} \ = \ D_x u(s,X_s^{t,x})1_{[t,T[}(s).
\]
Notice that identity \eqref{Identification2} holds taking $s=t$ in the first 
equality. Now, applying It\^o's formula to $u(s,X_s^{t,x})$ between $t$ and any $T_0\in[t,T[$, and using the fact that $u$ solves equation \eqref{KolmEq_Markov}, we see that \eqref{BSDE_Markov2} holds with $T_0$ in place of $T$ and $u(T_0,X_{T_0}^{t,x})$ in place of $h(X_T^{t,x})$. To conclude, it is enough to pass to the limit as $T_0\nearrow T$. This can be done using estimate \eqref{EstimateBSDE2} in Proposition \ref{P:EstimateBSDEAppendix} with $K\equiv0$. Finally, we notice that the present result is a slight generalization of Theorem 3.1 in \cite{pardoux_peng92}, 
since $u\in C^{1,2}([0,T[\times\R^d)\cap C([0,T]\times\R^d)$ instead of $u\in C^{1,2}([0,T]\times\R^d)$.
\ep

\vspace{3mm}

\noindent We can now present our first definition of strong-viscosity solution to equation \eqref{KolmEq_Markov}.

\begin{Definition}
\label{D:ViscosityFinite}
A function $u\colon[0,T]\times\R^d\rightarrow\R$ is called a \textbf{strong-viscosity solution} to equation \eqref{KolmEq_Markov} if there exists a sequence $(u_n,h_n,f_n,b_n,\sigma_n)_n$ of Borel measurable functions $u_n\colon[0,T]\times\R^d\rightarrow\R$, $h_n\colon\R^d\rightarrow\R$, $f_n\colon[0,T]\times\R^d\times\R\times\R^d\rightarrow\R$, $b_n\colon[0,T]\times\R^d\rightarrow\R^d$, and $\sigma_n\colon[0,T]\times\R^d\rightarrow\R^{d\times d}$, such that the following holds.
\begin{enumerate}
\item[\textup{(i)}]  For some positive constants $C$ and $m$,
\begin{align*}
|b_n(t,x)-b_n(t,x')| + |\sigma_n(t,x)-\sigma_n(t,x')| \ &\leq \ C|x-x'|, \notag \\
|f_n(t,x,y,z)-f_n(t,x,y',z')| \ &\leq \ C\big(|y-y'| + |z-z'|\big), \notag \\
|b_n(t,0)| + |\sigma_n(t,0)| \ &\leq \ C, \notag \\
|u_n(t,x)| + |h_n(x)| + |f_n(t,x,0,0)| \ &\leq \ C\big(1 + |x|^m\big),
\end{align*}
for all $t\in[0,T]$, $x,x'\in\R^d$, $y,y'\in\R$, and $z,z'\in\R^d$. Moreover, the functions $u_n(t,\cdot)$, $h_n(\cdot)$, $f_n(t,\cdot,\cdot,\cdot)$, $n\in\N$, are equicontinuous on compact sets, uniformly with respect to $t\in[0,T]$.
\item[\textup{(ii)}] $u_n$ is a classical solution to
\begin{equation}
\label{KolmEq_Markov_n}
\begin{cases}
\partial_t u_n(t,x) + \langle b_n(t,x),D_x u_n(t,x)\rangle + \frac{1}{2}\textup{tr}(\sigma_n\sigma_n\trans(t,x)D_x^2 u_n(t,x)) & \\
+\, f_n(t,x,u_n(t,x),\sigma_n\trans(t,x)D_x u_n(t,x)) \ = \ 0, &\hspace{-1.8cm}\forall\,(t,x)\in[0,T[\times\R^d, \\
u_n(T,x) \ = \ h_n(x), &\hspace{-1.8cm}\forall\,x\in\R^d.
\end{cases}
\end{equation}
\item[\textup{(iii)}] $(u_n,h_n,f_n,b_n,\sigma_n)$ converges pointwise to $(u,h,f,b,\sigma)$ as $n\rightarrow\infty$.
\end{enumerate}
\end{Definition}

\begin{Remark}
\label{R:UnifConvergence}
{\rm
(i) Notice that, for all $t\in[0,T]$, asking equicontinuity on compact sets of $(u_n(t,\cdot))_n$ together with its pointwise convergence to $u(t,\cdot)$ is equivalent to requiring the uniform convergence on compact sets of $(u_n(t,\cdot))_n$ to $u(t,\cdot)$. The same remark applies to $(h_n(\cdot))_n$ and $(f_n(t,\cdot,\cdot,\cdot))_n$. 

\noindent(ii) In Definition \ref{D:ViscosityFinite} we do not assume {\bf (A0)} for the functions $b,\sigma,f,h$. However, we can easily see that they satisfy automatically {\bf (A0)} as a consequence of point (i) of Definition \ref{D:ViscosityFinite}. See also Section \ref{R:Disc}.

\item(iii) 
We observe that a strong-viscosity solution to equation \eqref{KolmEq_Markov} in the sense of Definition \ref{D:ViscosityFinite} is a standard viscosity solution; for a definition we refer, e.g., to \cite{crandishiilions92}. Indeed, since a strong-viscosity solution $u$ to \eqref{KolmEq_Markov} is the limit of classical solutions (so, in particular, viscosity solutions) to perturbed equations, then from stability results for viscosity solutions (see, e.g., Lemma 6.1 and Remark 6.3 in \cite{crandishiilions92}), it follows that $u$ is a viscosity solution to equation \eqref{KolmEq_Markov}. On the other hand, if a strong-viscosity solution exists and a uniqueness result for viscosity solutions is in force, then a viscosity solution is a strong-viscosity solution, see also
 Remark \ref{R:Ishii}.
\ep
}
\end{Remark}

\begin{Theorem}
\label{T:UniqSV_Markov}
Let Assumption {\bf (A0)} hold and let $u\colon[0,T]\times\R^d\rightarrow\R$ be a strong-viscosity solution to equation \eqref{KolmEq_Markov}. Then, the following Feynman-Kac formula holds
\[
u(t,x) \ = \ Y_t^{t,x}, \qquad \forall\,(t,x)\in[0,T]\times\R^d,
\]
where $(Y_s^{t,x},Z_s^{t,x})_{s\in[t,T]}\in\S^2(t,T)\times\H^2(t,T)^d$, with $Y_s^{t,x}=u(s,X_s^{t,x})$, is the unique solution to the backward stochastic differential equation: $\P$-a.s.,
\begin{equation}
\label{BSDE_Markov}
Y_s^{t,x} \ = \ h(X_T^{t,x}) + \int_s^T f(r,X_r^{t,x},Y_r^{t,x},Z_r^{t,x}) dr - \int_s^T Z_r^{t,x} dW_r,
\end{equation}
for all $t \leq s \leq T$. In particular, there exists at most one strong-viscosity solution to equation \eqref{KolmEq_Markov}.
\end{Theorem}

The uniqueness Theorem \ref{T:UniqSV_Markov} will be proved in Section \ref{SubS:SecondDefnSV_Markov}, see Remark \ref{R:UniquenessSV_Markov}.

\subsection{Remarks in the case of discontinuous coefficients}
\label{R:Disc}

In the present section we will need the following additional assumption:

\vspace{3mm}

\noindent\textbf{(A0)'} \hspace{3mm} For every $(t,x)\in[0,T[\times\R^d$, $X_T^{t,x}$ has an absolutely continuous law with respect to the Lebesgue measure.

\vspace{3mm}

As noticed in Remark \ref{R:UnifConvergence}-(ii), in Definition \ref{D:ViscosityFinite} we easily see that the coefficients $b,\sigma,f,h$ satisfy automatically {\bf (A0)}. It also follows from point (i) of Definition \ref{D:ViscosityFinite} that $f(t,\cdot)$ and $h(\cdot)$ are continuous, uniformly with respect to $t\in[0,T]$. However, we can modify Definition \ref{D:ViscosityFinite} as follows in order to take into account the case where $f$ and $h$ are discontinuous.

\vspace{2mm}

\noindent\emph{Let Assumptions {\bf (A0)}-{\bf (A0)}' hold. A function $u\colon[0,T]\times\R^d\rightarrow\R$ is called a \textbf{strong-viscosity solution} to equation \eqref{KolmEq_Markov} if there exists a sequence $(u_n,b_n,\sigma_n)_n$ of Borel measurable functions $u_n\colon[0,T[\times\R^d\rightarrow\R$, $b_n\colon[0,T[\times\R^d\rightarrow\R^d$, and $\sigma_n\colon[0,T[\times\R^d\rightarrow\R^{d\times d}$, such that the following holds.
\begin{enumerate}
\item[\textup{(i)}]  For some positive constants $C$ and $m$
\begin{align*}
|b_n(t,x)-b_n(t,x')| + |\sigma_n(t,x)-\sigma_n(t,x')| \ &\leq \ C|x-x'|, \notag \\
|b_n(t,0)| + |\sigma_n(t,0)| \ &\leq \ C, \notag \\
|u_n(t,x)| \ &\leq \ C\big(1 + |x|^m\big),
\end{align*}
for all $t\in[0,T[$, $x,x'\in\R^d$. Moreover, the function $u_n(t,\cdot)$ is equicontinuous on compact sets, uniformly with respect to $t$ in any compact set of $[0,T[$.
\item[\textup{(ii)}] $u_n$ belongs to $C^{1,2}([0,T[\times\R^d)$ and satisfies
\begin{align*}
\partial_t u_n(t,x) + \langle b_n(t,x),D_x u_n(t,x)\rangle + \frac{1}{2}\textup{tr}(\sigma_n\sigma_n\trans(t,x)D_x^2 u_n(t,x)) & \\
+\,f(t,x,u_n(t,x),\sigma_n\trans(t,x)D_x u_n(t,x))& \ = \ 0, \quad \forall\,(t,x)\in[0,T[\times\R^d.
\end{align*}
\item[\textup{(iii)}] $(u_n,b_n,\sigma_n)$ converges pointwise to $(u,b,\sigma)$ on $[0,T[\times\R^d$ as $n\rightarrow\infty$.
\item[\textup{(iv)}] The set $D_T=\{x\in\R\colon\text{ $u(\cdot,\cdot)$ is discontinuous at }(T,x)\}$ has Lebesgue measure equal to zero.
\item[\textup{(v)}] For some positive constants $C$ and $m$
\[
|u(t,x)| \ \leq \ C\big(1 + |x|^m\big),
\]
for all $(t,x)\in[0,T]\times\R^d$.
\item[\textup{(vi)}] $u(T,x)=h(x)$ for almost every $x\in\R^d$.
\end{enumerate}
}

\vspace{2mm}


\noindent Notice that under Assumptions {\bf (A0)}-{\bf (A0)}' there exists at most one strong-viscosity solution (in the above sense) to equation \eqref{KolmEq_Markov}. Let us give an idea on how to prove this result: the proof is an adaptation of the proofs of Theorems \ref{T:UniqSV_Markov}. 
 Given $(t,x)\in[0,T[\times\R^d$, for every  $T_0\in[t,T[$ we apply Theorem \ref{T:UniqSV_Markov} on the time interval $[t,T_0]$ instead of $[t,T]$.
Indeed $u$ is a strong solution (in the sense of Definition 
\ref{D:ViscosityFinite})
of  \eqref{KolmEq_Markov},
replacing $T$ with $T_0$.
 Then, we get for $Y_s^{t,x}=u(s,X_s^{t,x})$, $s\in[t,T_0]$,
\begin{equation}\label{BSDE_Uniqueness}
Y_s^{t,x} \ = \ u(T_0,X_{T_0}^{t,x}) + \int_s^{T_0} f(r,X_r^{t,x},Y_r^{t,x},Z_r^{t,x}) dr - \int_s^{T_0} Z_r^{t,x} dW_r, \qquad t \leq s \leq T_0.
\end{equation}
Now, by item (iv) of and {\bf (A0)}'  that the event $N=\{\omega\in\Omega\colon X_T^{t,x}(\omega)\in D_T\}$ has probability zero. Then, up to a null subset of $\Omega$, we have
\[
\lim_{T_0\rightarrow T}u(T_0,X_{T_0}^{t,x}) \ = \ u(T,X_T^{t,x}) \ = \ h(X_T^{t,x}).
\]
We observe that the above limit also holds in $L^2$ since $u$ satisfies the polynomial growth condition (v). This allows to pass to the limit as $T_0\rightarrow T$ in \eqref{BSDE_Uniqueness} and to prove that $(Y_s^{t,x},Z_s^{t,x})_{s\in[t,T]}\in\S^2(t,T)\times\H^2(t,T)^d$, with $Y_s^{t,x}=u(s,X_s^{t,x})$, is the unique solution to the backward stochastic differential equation \eqref{BSDE_Markov}. From the uniqueness of $(Y^{t,x},Z^{t,x})$, and in particular of $Y_t^{t,x}$, we conclude that $u(t,x)=Y_t^{t,x}, t < T,$ is uniquely determined. Therefore, there exists at most one strong-viscosity solution (in the above sense) to equation \eqref{KolmEq_Markov}. This shows uniqueness for strong-viscosity solutions in the present sense.

\vspace{2mm}

Concerning the existence results, namely Theorems \ref{T:ExistSV_Markov} and \ref{T:ExistSV_Markov2}, they need $f$ and $h$ to be continuous. However, exploiting more refined results in the theory of regularity of parabolic equations (e.g., with $f$ and $h$ possibly discontinuous, but $\sigma$ such that uniform ellipticity holds), it is potentially possible to prove an existence result for a strong-viscosity solution in the sense of the above definition. We provide below a simple example where this can be done. This example also shows that there might be equations of the form \eqref{KolmEq_Markov} for which we have a unique strong-viscosity solution (in the above sense) but infinitely many viscosity solutions.

Take $d=1$, $b\equiv0$, $\sigma\equiv1$, $f\equiv0$, and $h(x)=1_{[1,\infty[}(x)$ for all $x\in\R$. Then equation \eqref{KolmEq_Markov} becomes 
\begin{equation}\label{PDE_Uniq}
\begin{cases}
\vspace{1mm}\partial_t u(t,x) + \frac{1}{2}D_x^2 u(t,x) \ = \ 0, \qquad &\forall\,(t,x)\in[0,T[\times\R, \\
u(T,x) \ = \ 1_{[1,\infty[}(x), &\forall\,x\in\R.
\end{cases}
\end{equation}
Notice that Assumptions {\bf (A0)}-{\bf (A0)'} hold, therefore there exists at most one strong-viscosity solution (in the above sense) to equation \eqref{PDE_Uniq}. Indeed, the unique strong-viscosity solution to equation \eqref{PDE_Uniq} is given by the following explicit formula:
\[
v(t,x) \ = \ 1 - \Phi\bigg(\frac{1 - x}{\sqrt{T - t}}\bigg), \quad \forall\,(t,x)\in[0,T[\times\R, \qquad\qquad v(T,x) \ = \ 1_{[1,\infty[}(x),
\quad \forall\,x\in\R,
\]
where $\Phi(z)=\int_{-\infty}^z \frac{1}{\sqrt{2\pi}}e^{-\frac{1}{2}z^2}dz$. As a matter of fact, $u=v$ fulfills previous items (iv), (v), (vi). Moreover
$v$ is $C^{1,2}([0,T[\times\R^d)$, so item (ii) holds with $b_n\equiv0$ and $\sigma_n\equiv 1$. Items (i) and (iii) are also trivially fulfilled:
we do not need to regularize the coefficients, so we take $u_n\equiv v$.
Finally $v$ a strong-viscosity solution in the sense of the present definition.

Let us now prove that there are infinitely many viscosity solutions to equation \eqref{PDE_Uniq}. Firstly, fix $a\in[0,1]$ and define
\[
v_a(t,x) \ = \ 1 - \Phi\bigg(\frac{1 - x}{\sqrt{T - t}}\bigg), \quad \forall\,(t,x)\in[0,T[\times\R, \qquad\qquad v_a(T,x) \ = \
\begin{cases}
0, \qquad & x \, < \, 1, \\
a, & x \, = \, 1, \\
1, & x \, > \, 1,
\end{cases}
\quad \forall\,x\in\R.
\]
Notice that $v\equiv v_1$. Let us prove that each $v_a$ is a (discontinuous) viscosity solution to equation \eqref{PDE_Uniq}: we refer for instance to Section 2 in \cite{ishii89} for the definition of discontinuous viscosity solution. As a matter of fact, consider the lower and upper semi-continuous envelopes of $v_a$:
\[
(v_a)_*(t,x) \ = \ \liminf_{\substack{(s,y)\rightarrow(t,x)\\ s<T}} v_a(s,y), \qquad\qquad (v_a)^*(t,x) \ = \ \limsup_{\substack{(s,y)\rightarrow(t,x)\\ s<T}} v_a(s,y),
\]
for every $(t,x)\in[0,T]\times\R$. Notice that $(v_a)_*\equiv v_0$ and $(v_a)^*\equiv v_1$, for every $a\in[0,1]$. It is easy to see that $v_0$ is a viscosity subsolution to equation \eqref{PDE_Uniq}, since $v_0(T,\cdot)\leq h$, moreover $v_0$ is $C^{1,2}([0,T[\times\R^d)$ and solves equation \eqref{PDE_Uniq} on $[0,T[\times\R$. Similarly, we see that $v_1$ is a viscosity supersolution to equation \eqref{PDE_Uniq}. This implies that each $v_a$ is a (discontinuous) viscosity solution to equation \eqref{PDE_Uniq}. We conclude that there is no uniqueness result for viscosity solutions to equation \eqref{PDE_Uniq}.

\subsection{Second definition of strong-viscosity solution}
\label{SubS:SecondDefnSV_Markov}

Our second definition of strong-viscosity solution to equation \eqref{KolmEq_Markov} is more in the spirit of the standard definition of viscosity solution, which is usually required to be both a viscosity subsolution and a viscosity supersolution. Indeed, we introduce the concept of \emph{generalized strong-viscosity solution}, which has to  be both a strong-viscosity subsolution and a strong-viscosity supersolution. As it will be clear from the definition, this new notion of solution is more general (in other words, weaker), than the concept of strong-viscosity solution given earlier in Definition \ref{D:ViscosityFinite}. For this reason, we added the adjective \emph{generalized} to its name.

\vspace{3mm}

\noindent First, we introduce the standard notions of classical sub and supersolution.
\begin{Definition}
A function $u\colon[0,T]\times\R^d\rightarrow\R$ is called a \textbf{classical subsolution} $($resp. \textbf{classical supersolution}$)$ to equation \eqref{KolmEq_Markov} if $u\in C^{1,2}([0,T[\times\R^d)\cap C([0,T]\times\R^d)$ and solves
\[
\begin{cases}
\partial_t u(t,x) + \langle b(t,x),D_x u(t,x)\rangle + \frac{1}{2}\textup{tr}(\sigma\sigma\trans(t,x)D_x^2 u(t,x)) & \\
\hspace{1.1cm}+\, f(t,x,u(t,x),\sigma\trans(t,x)D_x u(t,x)) \ \geq \ (\text{resp. $\leq$}) \ 0, &\forall\,(t,x)\in[0,T[\times\R^d, \\
u(T,x) \ \leq \ (\text{resp. $\geq$}) \ h(x), &\forall\,x\in\R^d.
\end{cases}
\]
\end{Definition}

\noindent We state the following probabilistic representation result for classical sub and supersolutions.

\begin{Proposition}
\label{P:UniqStrictFinite}
Suppose that Assumption {\bf (A0)} holds.\\
\textup{(i)} Let $u\colon[0,T]\times\R^d\rightarrow\R$ be a classical supersolution to equation \eqref{KolmEq_Markov}, satisfying the polynomial growth condition
\[
|u(t,x)| \ \leq \ C'\big(1 + |x|^{m'}\big), \qquad \forall\,(t,x)\in[0,T]\times\R^d,
\]
for some positive constants $C'$ and $m'$. Then, we have
\[
u(t,x) \ = \ Y_t^{t,x}, \qquad \forall\,(t,x)\in[0,T]\times\R^d,
\]
for some uniquely determined $(Y_s^{t,x},Z_s^{t,x},K_s^{t,x})_{s\in[t,T]}\in\S^2(t,T)\times\H^2(t,T)^d\times\A^{+,2}(t,T)$, with $(Y_s^{t,x},$ $Z_s^{t,x})=(u(s,X_s^{t,x}),\sigma\trans(s,X_s^{t,x})D_x u(s,X_s^{t,x})1_{[t,T[}(s))$, solving the backward stochastic differential equation, $\P$-a.s.,
\[
Y_s^{t,x} \ = \ h(X_T^{t,x}) + \int_s^T f(r,X_r^{t,x},Y_r^{t,x},Z_r^{t,x}) dr + K_T^{t,x} - K_s^{t,x} - \int_s^T \langle Z_r^{t,x},dW_r\rangle, \quad t \leq s \leq T.
\]
\textup{(ii)} Let $u\colon[0,T]\times\R^d\rightarrow\R$ be a classical subsolution to equation \eqref{KolmEq_Markov}, satisfying the polynomial growth condition
\[
|u(t,x)| \ \leq \ C'\big(1 + |x|^{m'}\big), \qquad \forall\,(t,x)\in[0,T]\times\R^d,
\]
for some positive constants $C'$ and $m'$. Then, we have
\[
u(t,x) \ = \ Y_t^{t,x}, \qquad \forall\,(t,x)\in[0,T]\times\R^d,
\]
for some uniquely determined $(Y_s^{t,x},Z_s^{t,x},K_s^{t,x})_{s\in[t,T]}\in\S^2(t,T)\times\H^2(t,T)^d\times\A^{+,2}(t,T)$, with $(Y_s^{t,x},$ $Z_s^{t,x})=(u(s,X_s^{t,x}),\sigma\trans(s,X_s^{t,x})D_x u(s,X_s^{t,x})1_{[t,T[}(s))$, solving the backward stochastic differential equation, $\P$-a.s.,
\[
Y_s^{t,x} \ = \ h(X_T^{t,x}) + \int_s^T f(r,X_r^{t,x},Y_r^{t,x},Z_r^{t,x}) dr - (K_T^{t,x} - K_s^{t,x}) - \int_s^T \langle Z_r^{t,x},dW_r\rangle, \quad t \leq s \leq T.
\]
\end{Proposition}
\textbf{Proof.}
The proof can be done proceeding as in the proof of Proposition \ref{P:UniqStrictFiniteStandard}, see Theorem 3.6 in \cite{cosso_russo14bis}.
\ep

\vspace{3mm}

\noindent We can now provide the definition of generalized strong-viscosity solution.

\begin{Definition}
\label{D:StrongSuperSub}
A function $u\colon[0,T]\times\R^d\rightarrow\R$ is called a \textbf{strong-viscosity supersolution} $($resp. \textbf{strong-viscosity subsolution}$)$ to equation \eqref{KolmEq_Markov} if there exists a sequence $(u_n,h_n,f_n,b_n,\sigma_n)_n$ of Borel measurable functions $u_n\colon[0,T]\times\R^d\rightarrow\R$, $h_n\colon\R^d\rightarrow\R$, $f_n\colon[0,T]\times\R^d\times\R\times\R^d\rightarrow\R$, $b_n\colon[0,T]\times\R^d\rightarrow\R^d$, and $\sigma_n\colon[0,T]\times\R^d\rightarrow\R^{d\times d}$, such that the following holds.
\begin{enumerate}
\item[\textup{(i)}] For some positive constants $C$ and $m$,
\begin{align*}
|b_n(t,x)-b_n(t,x')| + |\sigma_n(t,x)-\sigma_n(t,x')| \ &\leq \ C|x-x'|, \\
|f_n(t,x,y,z)-f_n(t,x,y',z')| \ &\leq \ C\big(|y-y'| + |z-z'|\big), \\
|b_n(t,0)| + |\sigma_n(t,0)| \ &\leq \ C, \\
|u_n(t,x)| + |h_n(x)| + |f_n(t,x,0,0)| \ &\leq \ C\big(1 + |x|^m\big),
\end{align*}
for all $t\in[0,T]$, $x,x'\in\R^d$, $y,y'\in\R$, and $z,z'\in\R^d$. Moreover, the functions $u_n(t,\cdot)$, $h_n(\cdot)$, $f_n(t,\cdot,\cdot,\cdot)$, $n\in\N$, are equicontinuous on compact sets, uniformly with respect to $t\in[0,T]$.
\item[\textup{(ii)}] $u_n$ is a classical supersolution $($resp. classical subsolution$)$ to
\[
\begin{cases}
\partial_t u_n(t,x) + \langle b_n(t,x),D_x u_n(t,x)\rangle + \frac{1}{2}\textup{tr}(\sigma_n\sigma_n\trans(t,x)D_x^2 u_n(t,x)) & \\
\hspace{2.8cm}+\, f_n(t,x,u_n(t,x),\sigma_n\trans(t,x)D_x u_n(t,x)) \ = \ 0, &\!\!\!\!\!\!\!\forall\,(t,x)\in[0,T[\times\R^d, \\
u_n(T,x) \ = \ h_n(x), &\!\!\!\!\!\!\!\forall\,x\in\R^d.
\end{cases}
\]
\item[\textup{(iii)}] $(u_n,h_n,f_n,b_n,\sigma_n)$ converges pointwise to $(u,h,f,b,\sigma)$ as $n\rightarrow\infty$.
\end{enumerate}
A function $u\colon[0,T]\times\R^d\rightarrow\R$ is called a \textbf{generalized strong-viscosity solution} to equation \eqref{KolmEq_Markov} if it is both a strong-viscosity supersolution and a strong-viscosity subsolution to \eqref{KolmEq_Markov}.
\end{Definition}

\begin{Remark}
\label{R:UnifConvergence2}
{\rm
Notice that a strong-viscosity subsolution (resp. supersolution) to equation \eqref{KolmEq_Markov} in the sense of Definition \ref{D:StrongSuperSub} is a standard viscosity subsolution (resp. supersolution), as it follows for instance from Lemma 6.1 and Remark 6.3 in \cite{crandishiilions92}. As a consequence, a generalized strong-viscosity solution is a standard viscosity solution.
}

\ep
\end{Remark}

We can now state the following probabilistic representation result for strong-viscosity sub and supersolutions, that is one of the main results of this paper, from which the comparison principle will follow in Corollary \ref{C:CompThm}.

\begin{Theorem}
\label{T:RepresentationSuperSub}
\textup{(1)} Let $u\colon[0,T]\times\R^d\rightarrow\R$ be a strong-viscosity supersolution to equation \eqref{KolmEq_Markov}. Then, we have
\[
u(t,x) \ = \ Y_t^{t,x}, \qquad \forall\,(t,x)\in[0,T]\times\R^d,
\]
for some uniquely determined $(Y_s^{t,x},Z_s^{t,x},K_s^{t,x})_{s\in[t,T]}\in\S^2(t,T)\times\H^2(t,T)^d\times\A^{+,2}(t,T)$, with $Y_s^{t,x}=u(s,X_s^{t,x})$, solving the backward stochastic differential equation, $\P$-a.s.,
\begin{equation}
\label{E:BSDE_StrongNonlinear_Super}
Y_s^{t,x} \ = \ Y_T^{t,x} + \int_s^T f(r,X_r^{t,x},Y_r^{t,x},Z_r^{t,x}) dr + K_T^{t,x} - K_s^{t,x} - \int_s^T \langle Z_r^{t,x},dW_r\rangle, \qquad t \leq s \leq T.
\end{equation}
\textup{(2)} Let $u\colon[0,T]\times\R^d\rightarrow\R$ be a strong-viscosity subsolution to equation \eqref{KolmEq_Markov}. Then, we have
\[
u(t,x) \ = \ Y_t^{t,x}, \qquad \forall\,(t,x)\in[0,T]\times\R^d,
\]
for some uniquely determined $(Y_s^{t,x},Z_s^{t,x},K_s^{t,x})_{s\in[t,T]}\in\S^2(t,T)\times\H^2(t,T)^d\times\A^{+,2}(t,T)$, with $Y_s^{t,x}=u(s,X_s^{t,x})$, solving the backward stochastic differential equation, $\P$-a.s.,
\begin{equation}
\label{E:BSDE_StrongNonlinear_Sub}
Y_s^{t,x} \ = \ Y_T^{t,x} + \int_s^T f(r,X_r^{t,x},Y_r^{t,x},Z_r^{t,x}) dr - \big(K_T^{t,x} - K_s^{t,x}\big) - \int_s^T \langle Z_r^{t,x},dW_r\rangle, \qquad t \leq s \leq T.
\end{equation}
\end{Theorem}
\textbf{Proof.}
We shall only prove statement (1), since (2) can be established similarly. To prove (1), consider a sequence $(u_n,h_n,f_n,b_n,\sigma_n)_n$ satisfying conditions (i)-(iii) of Definition \ref{D:StrongSuperSub}. For every $n\in\N$ and any $(t,x)\in[0,T]\times\R^d$, consider the stochastic equation, $\P$-a.s.,
\[
X_s \ = \ x + \int_t^s b_n(r,X_r) dr + \int_t^s \sigma_n(r,X_r) dW_r, \qquad t \leq s \leq T.
\]
It is well-known that there exists a unique solution $(X_s^{n,t,x})_{s\in[t,T]}$ to the above equation. Moreover, from Proposition \ref{P:UniqStrictFinite} we know that $u_n(t,x) = Y_t^{n,t,x}$, $(t,x)\in[0,T]\times\R^d$, for some $(Y_s^{n,t,x},Z_s^{n,t,x},$ $K_s^{n,t,x})_{s\in[t,T]}\in\S^2(t,T)\times\H^2(t,T)^d\times\A^{+,2}(t,T)$ solving the backward stochastic differential equation, $\P$-a.s.,
\[
Y_s^{n,t,x} \ = \ Y_T^{n,t,x} + \int_s^T f_n(r,X_r^{n,t,x},Y_r^{n,t,x},Z_r^{n,t,x}) dr + K_T^{n,t,x} - K_s^{n,t,x} - \int_s^T \langle Z_r^{n,t,x}, dW_r\rangle, \quad t \leq s \leq T.
\]
Notice that, from the uniform polynomial growth condition of $(u_n)_n$ and estimate \eqref{EstimateSupXn} in Lemma \ref{L:AppendixX} (for the particular case when $b_n$ and $\sigma_n$ only depend on the current value of path, rather than on all its past trajectory) we have, for any $p\geq1$,
\[
\sup_{n\in\N} \|Y^{n,t,x}\|_{\S^p(t,T)} \ < \ \infty.
\]
Then, it follows from Proposition \ref{P:EstimateBSDEAppendix}, the polynomial growth condition of $(f_n)_n$ in $x$, and the linear growth condition of $(f_n)_n$ in $(y,z)$, that
\[
\sup_n\big(\|Z^{n,t,x}\|_{\H^2(t,T)^d} + \|K^{n,t,x}\|_{\S^2(t,T)}\big) \ < \ \infty.
\]
Set $Y_s^{t,x}=u(s,X_s^{t,x})$, for any $s\in[t,T]$. Then, from the polynomial growth condition that $u$ inherits from the sequence $(u_n)_n$, and using estimate \eqref{EstimateSupXn} in Lemma \ref{L:AppendixX} (for the particular case of non-path-dependent $b_n$ and $\sigma_n$), we deduce that $\|Y^{t,x}\|_{\S^p(t,T)}<\infty$, for any $p\geq1$. In particular, $Y\in\S^2(t,T)$ and it is a continuous process. We also have, using
the convergence result \eqref{limXn-X} in Lemma \ref{L:AppendixX} (for the particular case of non-path-dependent $b_n$ and $\sigma_n$), that there exists a subsequence of $(X^{n,t,x})_n$, which we still denote $(X^{n,t,x})_n$, such that
\begin{equation}
\label{E:supX^n-X-->0}
\sup_{t\leq s\leq T}|X_s^{n,t,x}(\omega)-X_s^{t,x}(\omega)| \ \overset{n\rightarrow\infty}{\longrightarrow} \ 0, \qquad \forall\,\omega\in\Omega\backslash N,
\end{equation}
for some null measurable set $N\subset\Omega$. Moreover, from  estimate \eqref{EstimateSupXn} in Lemma \ref{L:AppendixX} (for the particular case of non-path-dependent $b_n$ and $\sigma_n$) it follows that, possibly enlarging $N$, $\sup_{t\leq s\leq T}(|X_s^{t,x}(\omega)|$ $+$ $|X_s^{n,t,x}(\omega)|)<\infty$, for any $n\in\N$ and any $\omega\in\Omega\backslash N$. Now, fix 
$\omega\in\Omega\backslash N$; then
\begin{align*}
|Y_s^{n,t,x}(\omega)-Y_s^{t,x}(\omega)| \ &= \ |u_n(s,X_s^{n,t,x}(\omega))-u(s,X_s^{t,x}(\omega))| \\
&= \ |u_n(s,X_s^{n,t,x}(\omega))-u_n(s,X_s^{t,x}(\omega))| + |u_n(s,X_s^{t,x}(\omega))-u(s,X_s^{t,x}(\omega))|.
\end{align*}
For any $\eps>0$, from point (iii) of Definition \ref{D:StrongSuperSub} it follows that there exists $n'\in\N$ such that
\[
|u_n(s,X_s^{t,x}(\omega))-u(s,X_s^{t,x}(\omega))| \ < \ \frac{\eps}{2}, \qquad \forall\,n\geq n'.
\]
On the other hand, from the equicontinuity on compact sets of $(u_n)_n$, we see that there exists $\delta>0$, independent of $n$, such that
\[
|u_n(s,X_s^{n,t,x}(\omega))-u_n(s,X_s^{t,x}(\omega))| \ < \ \frac{\eps}{2}, \qquad \text{if }|X_s^{n,t,x}(\omega)-X_s^{t,x}(\omega)|<\delta.
\]
Using \eqref{E:supX^n-X-->0}, we can find $n''\in\N$, $n''\geq n'$, such that
\[
\sup_{t\leq s\leq T}|X_s^{n,t,x}(\omega)-X_s^{t,x}(\omega)| \ < \ \delta, \qquad \forall\,n\geq n''.
\]
In conclusion, for any $\omega\in\Omega\backslash N$ and any $\eps>0$ there exists $n''\in\N$ such that
\[
|Y_s^{n,t,x}(\omega)-Y_s^{t,x}(\omega)| \ < \ \eps, \qquad \forall\,n\geq n''.
\]
Therefore, $Y_s^{n,t,x}(\omega)$ converges to $Y_s^{t,x}(\omega)$, as $n$ tends to infinity, for any $(s,\omega)\in[t,T]\times(\Omega\backslash N)$. In a similar way, we can prove that there exists a null measurable set $N'\subset\Omega$ such that $f_n(s,X_s^{n,t,x}(\omega),y,$ $z)\rightarrow f(s,X_s^{t,x}(\omega),y,z)$, for any $(s,\omega,y,z)\in[t,T]\times(\Omega\backslash N')\times\R\times\R^d$. As a consequence, the claim follows from Theorem \ref{T:LimitThmBSDE}.
\ep

\vspace{3mm}

We can finally state a comparison principle for strong-viscosity sub and supersolutions, which follows directly from the comparison theorem for BSDEs, for which we refer for instance to Theorem 1.3 in \cite{peng00}.

\begin{Corollary}[Comparison principle]
\label{C:CompThm}
Let $\check u\colon[0,T]\times\R^d\rightarrow\R$ $($resp. $\hat u\colon[0,T]\times\R^d\rightarrow\R$$)$ be a strong-viscosity subsolution $($resp. strong-viscosity supersolution$)$ to equation \eqref{KolmEq_Markov}. Then $\check u \leq \hat u$ on $[0,T]\times\R^d$. In particular, there exists at most one generalized strong-viscosity solution to equation \eqref{KolmEq_Markov}.
\end{Corollary}
\begin{Remark}
\label{R:UniquenessSV_Markov}
{\rm
\begin{description}
\item{(i)}
Notice that Theorem \ref{T:UniqSV_Markov} follows from Corollary \ref{C:CompThm}, since a strong-viscosity solution (Definition \ref{D:ViscosityFinite}) is in particular a generalized strong-viscosity solution.
\item{(ii)} There is no universal result concerning 
uniqueness for (classical) viscosity solutions. There are only partial
results, which impose several assumptions on the coefficients,
for instance  Theorem 7.4 in \cite{ishii89}.
\end{description}
}
\ep

\end{Remark}
\textbf{Proof.}
We know that $\check u(T,x) \leq g(x) \leq \hat u(T,x)$, for all $x\in\R^d$. Moreover, from Theorem \ref{T:RepresentationSuperSub} we have
\[
\check u(t,x) \ = \ \check Y_t^{t,x}, \qquad \hat u(t,x) \ = \ \hat Y_t^{t,x}, \qquad \text{for all }(t,x)\in[0,T]\times\R^d,
\]
for some $(\check Y_s^{t,x},\check Z_s^{t,x},\check K_s^{t,x})_{s\in[t,T]},(\hat Y_s^{t,x},\hat Z_s^{t,x},\hat K_s^{t,x})_{s\in[t,T]}\in\S^2(t,T)\times\H^2(t,T)^d\times\A^{+,2}(t,T)$ satisfying \eqref{E:BSDE_StrongNonlinear_Sub} and \eqref{E:BSDE_StrongNonlinear_Super}, respectively. Then, the result follows from a direct application of the comparison theorem for backward stochastic differential equations, see, e.g., Theorem 1.3 in \cite{peng00}.
\ep

\vspace{3mm}

Now, we present two existence (and uniqueness) results
 for strong-viscosity solutions to equation \eqref{KolmEq_Markov}.

\begin{Theorem}
\label{T:ExistSV_Markov}
Let Assumption {\bf (A0)} hold and suppose that $b=b(x)$ and $\sigma=\sigma(x)$ do not depend on $t$. Suppose also that the functions $f$ and $h$ are continuous. Then, the function $u$ given by
\begin{equation}
\label{Identification}
u(t,x) \ = \ Y_t^{t,x}, \qquad \forall\,(t,x)\in[0,T]\times\R^d,
\end{equation}
where $(Y_s^{t,x},Z_s^{t,x})_{s\in[t,T]}\in\S^2(t,T)\times\H^2(t,T)^d$ is the unique solution to \eqref{BSDE_Markov}, is a strong-viscosity solution to equation \eqref{KolmEq_Markov}.
\end{Theorem}
\begin{Remark}
\label{R:Ishii}
{\rm
 Since the seminal paper \cite{pardoux_peng92}, we know that the
function defined in \eqref{Identification}
is a viscosity solution.
\ep
}
\end{Remark}
\textbf{Proof (of Theorem \ref{T:ExistSV_Markov}).}
Let us fix some notations. Let $q\in\N\backslash\{0\}$ and consider the function $\phi_q\in C^\infty(\R^q)$ given by
\[
\phi_q(w) \ = \ c \exp\bigg(\frac{1}{|w|^2 - 1}\bigg) 1_{\{|w|<1\}}, \qquad \forall\,w\in\R^q,
\]
with $c>0$ such that $\int_{\R^q} \phi_q(w) dw = 1$. Then, we define $\phi_{q,n}(w) = n^q\phi_q(nw)$, $\forall\,w\in\R^q$, $n\in\N$. Let us now define, for any $n\in\N$,
\begin{align*}
b_n(x) \ &= \ \int_{\R^d} \phi_{d,n}(x') b(x - x') dx', \qquad \sigma_n(x) \ = \ \int_{\R^d} \phi_{d,n}(x') \sigma(x - x') dx', \\
f_n(t,x,y,z) \ &= \ \int_{\R^d\times\R\times\R^d} \phi_{2d+1,n}(x',y',z') f(t,x-x',y-y',z-z') dx'dy'dz', \\
h_n(x) \ &= \ \int_{\R^d} \phi_{d,n}(x') h(x - x') dx',
\end{align*}
for all $(t,x,y,z)\in[0,T]\times\R^d\times\R\times\R^d$. Then, we see that the sequence of continuous functions $(b_n,\sigma_n,f_n,h_n)_n$ satisfies assumptions (i) and (iii) of Definition \ref{D:ViscosityFinite}. Moreover, for any $n\in\N$
we have the following.
\begin{itemize}
\item $b_n$ and $\sigma_n$ are of class $C^3$ with partial derivatives from order $1$ up to order $3$ bounded.
\item For all $t\in[0,T]$, $f_n(t,\cdot,\cdot,\cdot)\in C^3(\R^d\times\R\times\R^d)$ and the two properties below.
\begin{itemize}
\item $f_n(t,\cdot,0,0)$ belongs to $C^3$ and its third order partial derivatives satisfy a polynomial growth condition uniformly in $t$.
\item $D_y f_n$, $D_z f_n$ are bounded on $[0,T]\times\R^d\times\R\times\R^d$, as well as their derivatives of order one and second with respect to $x,y,z$.
\end{itemize}
\item $h_n\in C^3(\R^d)$ and its third order partial derivatives satisfy a polynomial growth condition.
\end{itemize}
Therefore, it follows from Theorem 3.2 in \cite{pardoux_peng92} that a classical solution to equation \eqref{KolmEq_Markov_n} is given by
\begin{equation}
\label{Identification_un}
u_n(t,x) \ = \ Y_t^{n,t,x}, \qquad \forall\,(t,x)\in[0,T]\times\R^d,
\end{equation}
where $(Y_s^{n,t,x},Z_s^{n,t,x})_{s\in[t,T]}\in\S^2(t,T)\times\H^2(t,T)^d$ is the unique solution to the backward stochastic differential equation: $\P$-a.s.,
\[
Y_s^{n,t,x} \ = \ h_n(X_T^{n,t,x}) + \int_s^T f_n(r,X_r^{n,t,x},Y_r^{n,t,x},Z_r^{n,t,x}) dr - \int_s^T Z_r^{n,t,x} dW_r, \qquad t \leq s \leq T,
\]
with
\[
X_s^{n,t,x} \ = \ x + \int_t^s b_n(r,X_r^{n,t,x}) dr + \int_t^s \sigma_n(r,X_r^{n,t,x}) dW_r, \qquad t \leq s \leq T.
\]
From \eqref{Identification_un}, Proposition \ref{P:EstimateBSDEAppendix}, and estimate \eqref{EstimateSupXn}, we see that $u_n$ satisfies a polynomial growth condition uniform in $n$. It remains to prove that the sequence $(u_n)_n$ converges pointwise to $u$ as $n\rightarrow\infty$, and that the functions $u_n(t,\cdot)$, $n\in\N$, are equicontinuous on compact sets, uniformly with respect to $t\in[0,T]$. Concerning this latter property, fix $t\in[0,T]$, a compact subset $K\subset\R^d$, and $\eps>0$. We have to prove that there exists $\delta=\delta(\eps,K)$ such that
\begin{equation}
\label{UniformContinuity_un}
|u_n(t,x) - u_n(t,x')| \ \leq \ \eps, \qquad \text{if }|x-x'|\leq\delta,\,x,x'\in K.
\end{equation}
To this end, we begin noting that from estimate \eqref{EstimateBSDE2} we have that there exists a constant $C$, independent of $n$, such that
\begin{align*}
|u_n(t,x) - u_n(t,x')|^2 \ &\leq \ C\,\E\big[\big|h_n(X_T^{n,t,x}) - h_n(X_T^{n,t,x'})\big|^2\big] \\
&\quad \ + C\int_t^T \E\big[\big|f_n(s,X_s^{n,t,x},Y_s^{n,t,x},Z_s^{n,t,x}) - f_n(s,X_s^{n,t,x'},Y_s^{n,t,x},Z_s^{n,t,x})\big|^2\big] ds, \notag
\end{align*}
for all $t\in[0,T]$ and $x,x'\in\R^d$. In order to prove \eqref{UniformContinuity_un}, we also recall the following standard estimate: for any $p\geq2$ there exists a positive constant $C_p$, independent of $n$, such that
\[
\E\big[\big|X_s^{n,t,x} - X_s^{n,t,x'}\big|^p\big] \ \leq \ C_p |x - x'|^p,
\]
for all $t\in[0,T]$, $s\in[t,T]$, $x,x'\in\R^d$, $n\in\N$. Now, choose $p>d$, $R>0$, and $\alpha\in]0,p-d[$. Then, it follows from Garsia-Rodemich-Rumsey lemma (see, in particular, formula (3a.2) in \cite{barlow_yor82}) that, for all $t\in[0,T]$, $s\in[t,T]$, $x,x'\in\R^d$, $n\in\N$,
\begin{equation}
\label{UniformContinuityX}
|X_s^{n,t,x} - X_s^{n,t,x'}| \ \leq \ (\Gamma_s^{n,t})^{1/p} |x - x'|^{\alpha/p},
\end{equation}
for some process $\Gamma^{n,t}=(\Gamma_s^{n,t})_{s\in[t,T]}$ given by
\[
\Gamma_s^{n,t} \ = \ C_d\,8^p\,2^\alpha\bigg(1 + \frac{2d}{\alpha}\bigg)\int_{\{(y,y')\in\R^{2d}\colon|y|,|y'|\leq R\}} \frac{|X_s^{n,t,y} - X_s^{n,t,y'}|}{|y - y'|^{\alpha + 2d}} dydy'
\]
and
\begin{equation}
\label{ExpectedValueGamma}
\E\big[\Gamma_s^{n,t}\big] \ \leq \ C_d\,C_p \frac{1}{(p - d) - \alpha} R^{p - \alpha},
\end{equation}
where $C_d$ is a universal constant depending only on $d$.

Now, let us prove that
\begin{equation}
\label{UniformContinuity_hn}
\E\big[\big|h_n(X_T^{n,t,x}) - h_n(X_T^{n,t,x'})\big|^2\big] \ \leq \ \eps, \qquad \text{if }|x-x'|\leq\delta,\,x,x'\in K.
\end{equation}
Let $x,x'\in K$ and let $m$ be a strictly positive integer to be chosen later. Then, consider the event (we omit the dependence on $t,x$)
\[
\Omega_{n,m} \ = \ \big\{\omega\in\Omega\colon\Gamma_T^{n,t}(\omega)\leq m,\,|X_T^{n,t,x}(\omega)|\leq m\big\}.
\]
From \eqref{UniformContinuityX} we see that, on $\Omega_{n,m}$, $X_T^{n,t,x'}$ is also uniformly bounded by a constant independent of $n,t,x,x'$, since $x,x'\in K$. In particular, from the equicontinuity on compact sets of the sequence $(h_n)_n$, it follows that there exists a continuity modulus $\rho$ (depending on $K$, but independent of $n$) such that
\begin{align*}
\E\big[\big|h_n(X_T^{n,t,x}) - h_n(X_T^{n,t,x'})\big|^2\big] \ &\leq \ \E\big[\rho^2(|X_T^{n,t,x} - X_T^{n,t,x'}|) 1_{\Omega_{n,m}}\big] \\
&\quad \ + \E\big[\big|h_n(X_T^{n,t,x}) - h_n(X_T^{n,t,x'})\big|^2 1_{\Omega_{n,m}^c}\big].
\end{align*}
By \eqref{UniformContinuityX} and Cauchy-Schwarz inequality
\begin{align*}
\E\big[\big|h_n(X_T^{n,t,x}) - h_n(X_T^{n,t,x'})\big|^2\big] \ &\leq \ \rho^2\big(m^{1/p}|x - x'|^{\alpha/p}\big) \\
&+ \sqrt{\E\big[\big|h_n(X_T^{n,t,x}) - h_n(X_T^{n,t,x'})\big|^4\big]} \sqrt{\P(\Gamma_T^{n,t}>m) + \P(X_T^{n,t,x}>m)}.
\end{align*}
From the standard inequalities $|a-b|^4\leq8(a^4+b^4)$, $\forall\,a,b\in\R$, and $\sqrt{c+d}\leq\sqrt{c}+\sqrt{d}$, $\forall\,c,d\geq0$, we see that
\[
\sqrt{\E\big[\big|h_n(X_T^{n,t,x}) - h_n(X_T^{n,t,x'})\big|^4\big]} \ \leq \ \sqrt{8\E\big[\big|h_n(X_T^{n,t,x})\big|^4\big]} + \sqrt{8\E\big[\big|h_n(X_T^{n,t,x'})\big|^4\big]}.
\]
Now, using this estimate, the polynomial growth condition of $h_n$ (uniform in $n$), estimate \eqref{EstimateSupXn}, estimate \eqref{ExpectedValueGamma}, and Chebyshev's inequality, we obtain
\begin{align*}
\sqrt{\E\big[\big|h_n(X_T^{n,t,x}) - h_n(X_T^{n,t,x'})\big|^4\big]} \ &\leq \ C_K, \\
\P(\Gamma_T^{n,t}>m) \ &\leq \ \frac{\E\big[\Gamma_T^{n,t}\big]}{m} \ \leq \ \frac{C_K}{m}, \\
\P(X_T^{n,t,n}>m) \ &\leq \ \frac{\E\big[|X_T^{n,t,x}|\big]}{m} \ \leq \ \frac{C_K}{m},
\end{align*}
for some positive constant $C_K$, possibly depending on $K$ (in particular, on $x$ and $x'$), but independent of $n,t$. Therefore, we see that we can find $m=m(\eps,K)$ large enough such that
\[
\E\big[\big|h_n(X_T^{n,t,x}) - h_n(X_T^{n,t,x'})\big|^2\big] \ \leq \ \rho^2\big(m^{1/p}|x - x'|^{\alpha/p}\big) + \frac{\eps}{2}.
\]
Then, there exists $\delta=\delta(\eps,K)>0$ such that \eqref{UniformContinuity_hn} holds. In a similar way we can prove that, possibly taking a smaller $\delta=\delta(\eps,K)>0$, we have
\begin{equation}
\label{UniformContinuity_fn}
\E\big[\big|f_n(s,X_s^{n,t,x},Y_s^{n,t,x},Z_s^{n,t,x}) - f_n(s,X_s^{n,t,x'},Y_s^{n,t,x},Z_s^{n,t,x})\big|^2\big] \ \leq \ \eps,
\end{equation}
if $|x-x'|\leq\delta$, $x,x'\in K$, $\forall\,s\in[t,T]$. By \eqref{UniformContinuity_hn} and \eqref{UniformContinuity_fn} we deduce the validity of \eqref{UniformContinuity_un}.

Finally, let us prove the pointwise convergence of the sequence $(u_n)_n$ to $u$. Using again estimate \eqref{EstimateBSDE2}, we find
\begin{align}
|u_n(t,x) - u(t,x)|^2 \ &\leq \ C\,\E\big[\big|h_n(X_T^{n,t,x}) - h(X_T^{t,x})\big|^2\big] \label{unx-ux} \\
&\quad \ + C\int_t^T \E\big[\big|f_n(s,X_s^{n,t,x},Y_s^{t,x},Z_s^{t,x}) - f(s,X_s^{t,x},Y_s^{t,x},Z_s^{t,x})\big|^2\big] ds, \notag
\end{align}
$\forall\,(t,x)\in[0,T]\times\R^d$, $n\in\N$, for some constant $C$, independent of $n$  and depending only on the (uniform in $n$) Lipschitz constant of $f_n$ with respect to $(y,z)$. By the uniform convergence on compact sets of $(h_n(\cdot),f_n(t,\cdot,y,z))_n$ to $(h(\cdot),f(t,\cdot,y,z))$, we have, $\P$-a.s.,
\begin{align}
h_n(X_T^{n,t,x}) \ &\overset{n\rightarrow\infty}{\longrightarrow} \ h(X_T^{t,x}), \label{Conv_h_n} \\
f_n(s,X_s^{n,t,x},Y_s^{t,x},Z_s^{t,x}) \ &\overset{n\rightarrow\infty}{\longrightarrow} \ f(s,X_s^{t,x},Y_s^{t,x},Z_s^{t,x}), \label{Conv_f_n}
\end{align}
for all $s\in[t,T]$. By Assumption \textbf{(A0)} and the polynomial growth condition of $h_n$, $f_n$, $u_n$ (uniform in $n$), estimates \eqref{EstimateSupX} and \eqref{EstimateSupXn}, Proposition \ref{P:EstimateBSDEAppendix}, we can prove the uniform integrability of the sequences $(|h_n(X_T^{n,t,x}) - h(X_T^{t,x})|^2)_n$ and $(|f_n(s,X_s^{n,t,x},Y_s^{t,x},Z_s^{t,x}) - f(s,X_s^{t,x},Y_s^{t,x},Z_s^{t,x})|^2)_n$, $\forall\,s\in[t,T]$. This, together with \eqref{Conv_h_n}-\eqref{Conv_f_n}, implies that
\begin{align*}
\E\big[\big|h_n(X_T^{n,t,x}) - h(X_T^{t,x})\big|^2\big] \ &\overset{n\rightarrow\infty}{\longrightarrow} \ 0, \\
\E\big[\big|f_n(s,X_s^{n,t,x},Y_s^{t,x},Z_s^{t,x}) - f(s,X_s^{t,x},Y_s^{t,x},Z_s^{t,x})\big|^2\big] \ &\overset{n\rightarrow\infty}{\longrightarrow} \ 0,
\end{align*}
for all $s\in[t,T]$. From the second convergence, the polynomial growth condition of $f$ and $f_n$ (uniform in $n$), estimates \eqref{EstimateSupX} and \eqref{EstimateSupXn}, it follows that
\[
\lim_{n\rightarrow\infty} \int_t^T \E\big[\big|f_n(s,X_s^{n,t,x},Y_s^{t,x},Z_s^{t,x}) - f(s,X_s^{t,x},Y_s^{t,x},Z_s^{t,x})\big|^2\big] ds \ = \ 0.
\]
In conclusion, we can pass to the limit in \eqref{unx-ux} as $n\rightarrow\infty$, and we obtain the pointwise convergence of $(u_n)_n$ to $u$.
\ep

\begin{Remark}
\label{R:TwoDefinitions}
{\rm
Notice that Theorem \ref{T:ExistSV_Markov} gives an existence result for strong-viscosity solutions (see Definition \ref{D:ViscosityFinite}) to equation \eqref{KolmEq_Markov}, which implies an existence result for generalized strong-viscosity solutions (see Definition \ref{D:StrongSuperSub}). In Section \ref{S:SVPath} we will consider only Definition \ref{D:ViscosityFinite} and extend it to the path-dependent case.
\ep
}
\end{Remark}

We conclude this section providing another existence result for strong-viscosity solutions to equation \eqref{KolmEq_Markov} under a different set of assumptions with respect to Theorem \ref{T:ExistSV_Markov}. In particular, $f=f(t,x)$ does not depend on $(y,z)$, while $b$ and $\sigma$ can depend on $t$.

\begin{Theorem}
\label{T:ExistSV_Markov2}
Let Assumption {\bf (A0)} hold and suppose that $f=f(t,x)$ does not depend on $(y,z)$. Suppose also that the functions $f$ and $h$ are continuous. Then, the function $u$ given by
\[
u(t,x) \ = \ Y_t^{t,x}, \qquad \forall\,(t,x)\in[0,T]\times\R^d,
\]
where $(Y_s^{t,x},Z_s^{t,x})_{s\in[t,T]}\in\S^2(t,T)\times\H^2(t,T)^d$ is the unique solution to \eqref{BSDE_Markov}, is a strong-viscosity solution to equation \eqref{KolmEq_Markov}.
\end{Theorem}
\textbf{Proof.}
The proof can be done proceeding as in the proof of Theorem \ref{T:ExistSV_Markov}, by smoothing the coefficients, but using Theorem 6.1, Chapter 5, in \cite{friedman75vol1} instead of Theorem 3.2 in \cite{pardoux_peng92}.
\ep

\section{Strong-viscosity solutions in the path-dependent case}
\label{S:SVPath}

One of the goals of the present section is to show that the notion of strong-viscosity solution is very flexible and easy to extend, with respect to the standard notion of viscosity solution, to more general settings than the Markovian one. In particular, we focus on semilinear parabolic path-dependent PDEs.

\subsection{Semilinear parabolic path-dependent PDEs}
\label{Ss:PPDE}

Let us denote by $C([-T,0])$ the Banach space of all continuous paths $\eta\colon[-T,0]\rightarrow\R$ endowed with the supremum norm $\|\eta\|=\sup_{t\in[-T,0]}|\eta(t)|$. Let us consider the following semilinear parabolic path-dependent PDE (for simplicity of notation, we consider the unidimensional case, with $\eta$ taking values in $\R$):
\begin{equation}
\label{KolmEq}
\begin{cases}
\partial_t \Uc + D^H \Uc + b(t,\eta)D^V \Uc + \frac{1}{2}\sigma(t,\eta)^2D^{VV} \Uc \\
\hspace{2.3cm}+\, F(t,\eta,\Uc,\sigma(t,\eta)D^V \Uc) \ = \ 0, \;\;\; &\quad\forall\,(t,\eta)\in[0,T[\times C([-T,0]), \\
\Uc(T,\eta) \ = \ {H}(\eta), &\quad\forall\,\eta\in C([-T,0]),
\end{cases}
\end{equation}
where $D^H\Uc$, $D^V\Uc$, $D^{VV}\Uc$ are the functional derivatives introduced in \cite{cosso_russo15a}, whose definition is recalled below. Concerning the coefficients $b\colon[0,T]\times C([-T,0])\rightarrow\R$, $\sigma\colon[0,T]\times C([-T,0])\rightarrow\R$, $F\colon[0,T]\times C([-T,0])\times\R\times\R\rightarrow\R$, and ${H}\colon C([-T,0])\rightarrow\R$ of equation \eqref{KolmEq}, we shall impose the following assumptions.

\vspace{3mm}

\noindent\textbf{(A1)} \hspace{3mm} $b$, $\sigma$, $F$, ${H}$ are Borel measurable functions satisfying, for some positive constants $C$ and $m$,
\begin{align*}
|b(t,\eta)-b(t,\eta')| + |\sigma(t,\eta)-\sigma(t,\eta')| \ &\leq \ C\|\eta-\eta'\|, \\
|F(t,\eta,y,z)-F(t,\eta,y',z')| \ &\leq \ C\big(|y-y'| + |z-z'|\big), \\
|b(t,0)| + |\sigma(t,0)| \ &\leq \ C, \\
|F(t,\eta,0,0)| + |{H}(\eta)| \ &\leq \ C\big(1 + \|\eta\|^m\big),
\end{align*}
for all $t\in[0,T]$, $\eta,\eta'\in C([-T,0])$, $y,y',z,z'\in\R$.

\subsection{Recall on functional It\^o calculus}

In the present subsection we recall the results of functional It\^o calculus needed later, without pausing on the technicalities and focusing on the intuition. For all technical details and rigorous definitions, we refer to \cite{cosso_russo15a}.

\vspace{1mm}

We begin introducing the \emph{functional derivatives}. To this end, it is useful to think of $\Uc=\Uc(t,\eta)$ as $\Uc=\Uc(t,\eta(\cdot)1_{[-T,0[}+\eta(0)1_{\{0\}})$, in order to emphasize the past $\eta(\cdot)1_{[-T,0[}$ and present $\eta(0)$ of the path $\eta$. Then, we can give, at least formally, the following definitions, see Definition 2.23 in \cite{cosso_russo15a}.

\vspace{1mm}

\noindent$\bullet$ \emph{Horizontal derivative.} We look at the sensibility of $\Uc$ with respect to a constant extension of the past $\eta(\cdot)1_{[-T,0[}$, keeping fixed the present value at $\eta(0)$:
\[
D^H \Uc(t,\eta) \ := \ \lim_{\eps\rightarrow0^+} \frac{\Uc(t,\eta(\cdot)1_{[-T,0[}+\eta(0)1_{\{0\}}) - \Uc(t,\eta(\cdot-\eps)1_{[-T,0[}+\eta(0)1_{\{0\}})}{\eps}.
\]
$\bullet$ \emph{First vertical derivative.} We look at the first variation with respect to the present, with the past fixed:
\[
D^V \Uc(t,\eta) \ := \ \lim_{\eps\rightarrow0} \frac{\Uc(t,\eta(\cdot)1_{[-T,0[}+(\eta(0)+\eps)1_{\{0\}}) - \Uc(t,\eta(\cdot)1_{[-T,0[}+\eta(0)1_{\{0\}})}{\eps}.
\]
$\bullet$ \emph{Second vertical derivative.} We look at the second variation with respect to the present, with the past fixed:
\[
D^{VV} \Uc(t,\eta) \ := \ \lim_{\eps\rightarrow0} \frac{D^V \Uc(t,\eta(\cdot)1_{[-T,0[}+(\eta(0)+\eps)1_{\{0\}}) - D^V \Uc(t,\eta(\cdot)1_{[-T,0[}+\eta(0)1_{\{0\}})}{\eps}.
\]

\vspace{1mm}

Given $I=[0,T[$ or $I=[0,T]$, we say that $\Uc\colon I\times C([-T,0])\rightarrow\R$ is of class $C^{1,2}((I\times\textup{past})\times\textup{present}))$ if, roughly speaking, $\partial_t \Uc$, $D^H \Uc$, $D^V \Uc$, and $D^{VV}\Uc$ exist and are continuous together with $\Uc$, for a rigorous definition we refer to \cite{cosso_russo15a}, Definition 2.28.

\vspace{1mm}    
    
We can finally state the \emph{functional It\^o formula}. Firstly, we fix some notation. As in Section \ref{S:SV_Markov}, we consider a complete probability space $(\Omega,\Fc,\P)$. Given a real-valued continuous process $X = (X_t)_{t\in[0,T]}$ on $(\Omega,\Fc,\P)$, we extend it to all $t\in\R$ in a canonical way as follows: $X_t:=X_0$, $t<0$, and $X_t:=X_T$, $t>T$; then, we associate to $X$ the so-called \emph{window process} $\X=(\X_t)_{t\in\R}$, which is a $C([-T,0])$-valued process given by
\[
\X_t := \{X_{t+s},\,s\in[-T,0]\}, \qquad t\in\R.
\]

\begin{Theorem}
\label{T:ItoTime}
Let $\Uc\colon[0,T]\times C([-T,0])\rightarrow\R$ be of class $C^{1,2}(([0,T]\times\textup{past})\times\textup{present})$ and $X=(X_t)_{t\in[0,T]}$ be a real continuous finite quadratic variation process. Then, the following \textbf{functional It\^o formula} holds, $\P$-a.s.,
\begin{align*}
\Uc(t,\X_t) \ &= \ \Uc(0,\X_0) + \int_0^t \big(\partial_t \Uc(s,\X_s) + D^H \Uc(s,\X_s)\big)ds + \int_0^t D^V \Uc(s,\X_s) d^- X_s \notag \\
&\quad \ + \frac{1}{2}\int_0^t D^{VV}\Uc(s,\X_s)d[X]_s,
\end{align*}
for all $0 \leq t \leq T$.
\end{Theorem}

\begin{Remark}
{\rm
(i) The term $\int_0^t D^V \Uc(s,\X_s) d^- X_s$ denotes the \emph{forward integral} of $D^V \Uc(\cdot,\X_\cdot)$ with respect to $X$ defined by \emph{regularization} (see, e.g., \cite{russovallois91,russovallois93,russovallois07}), which coincides with the classical stochastic integral whenever $X$ is a semimartingale.

\vspace{1mm}

\noindent(ii) In the non-path-dependent case $\Uc(t,\eta)=F(t,\eta(0))$, for any $(t,\eta)\in[0,T]\times C([-T,0])$, with $F\in C^{1,2}([0,T]\times\R)$, we retrieve the finite-dimensional It\^o formula, see Theorem 2.1 of \cite{russovallois95}.
\ep
}
\end{Remark}

\subsection{Recall on strict solutions}

We recall the concept of strict solution to equation \eqref{KolmEq} from Section 3 in \cite{cosso_russo15a}.

\begin{Definition}
A map $\Uc\colon[0,T]\times C([-T,0])\rightarrow\R$ in $C^{1,2}(([0,T[\times\textup{past})\times\textup{present})\cap C([0,T]\times C([-T,0]))$, satisfying equation \eqref{KolmEq}, is called a \textbf{strict solution} to equation \eqref{KolmEq}.
\end{Definition}

We present now a probabilistic representation result, for which we adopt the same notations as
 in Section \ref{SubS:Notation}, with dimension $d=1$. First, we recall some preliminary results. More precisely, for any $(t,\eta)\in[0,T]\times C([-T,0])$, we consider the \emph{path-dependent SDE}
\begin{equation}
\label{SDE}
\begin{cases}
dX_s = \ b(s,\mathbb X_s)dt + \sigma(s,\mathbb X_s)dW_s, \qquad\qquad & s\in[t,T], \\
X_s \ = \ \eta(s-t), & s\in[-T+t,t].
\end{cases}
\end{equation}

\begin{Proposition}
\label{P:SDE}
Under Assumption {\bf (A1)}, for any $(t,\eta)\in[0,T]\times C([-T,0])$ there exists a unique $($up to indistinguishability$)$ $\F$-adapted continuous process $X^{t,\eta}=(X_s^{t,\eta})_{s\in[-T+t,T]}$ strong solution to equation \eqref{SDE}. Moreover, for any $p\geq1$ there exists a positive constant $C_p$ such that
\begin{equation}
\label{EstimateX}
\E\Big[\sup_{s\in[-T+t,T]}\big|X_s^{t,\eta}\big|^p\Big] \ \leq \ C_p \big(1 + \|\eta\|^p\big).
\end{equation}
\end{Proposition}
\textbf{Proof.}
See Lemma \ref{L:SDE}.
\ep

\begin{Theorem}
\label{T:UniqStrict}
Suppose that Assumption {\bf (A1)} holds. Let $\Uc\colon[0,T]\times C([-T,0])\rightarrow\R$ be a strict solution to equation \eqref{KolmEq} satisfying the polynomial growth condition
\begin{equation}
\label{PolGrowthCondPath}
|\Uc(t,\eta)| \ \leq \ C'\big(1 + \|\eta\|^{m'}\big), \qquad \forall\,(t,\eta)\in[0,T]\times C([-T,0]),
\end{equation}
for some positive constants $C'$ and $m'$. Then, the following Feynman-Kac formula holds
\[
\Uc(t,\eta) \ = \ Y_t^{t,\eta}, \qquad \forall\,(t,\eta)\in[0,T]\times C([-T,0]),
\]
where $(Y_s^{t,\eta},Z_s^{t,\eta})_{s\in[t,T]} = (\Uc(s,\mathbb X_s^{t,\eta}),\sigma(s,\mathbb X_s^{t,\eta})D^V\Uc(s,\mathbb X_s^{t,\eta})1_{[t,T[}(s))_{s\in[t,T]}\in\S^2(t,T)\times\H^2(t,T)$ is the unique solution to the backward stochastic differential equation: $\P$-a.s.,
\[
Y_s^{t,\eta} \ = \ {H}(\mathbb X_T^{t,\eta}) + \int_s^T F(r,\mathbb X_r^{t,\eta},Y_r^{t,\eta},Z_r^{t,\eta}) dr - \int_s^T Z_r^{t,\eta} dW_r, \qquad t \leq s \leq T.
\]
In particular, there exists at most one strict solution to equation \eqref{KolmEq} satisfying a polynomial growth condition as in \eqref{PolGrowthCondPath}.
\end{Theorem}
\textbf{Proof.}
See Theorem 3.4 in \cite{cosso_russo15a}.
\ep

\vspace{3mm}

We state the following existence result.

\begin{Theorem}
\label{T:ExistenceStrict}
Suppose that there exists $N\in\N\backslash\{0\}$ such that, for all $(t,\eta,y,z)\in[0,T]\times C([-T,0])$ $\times\R\times\R$ 
\begin{align*}
b(t,\eta) \ &= \ \bar b\bigg(\int_{[-t,0]}\varphi_1(x+t)d^-\eta(x),\ldots,\int_{[-t,0]}\varphi_N(x+t)d^-\eta(x)\bigg), \\
\sigma(t,\eta) \ &= \ \bar\sigma\bigg(\int_{[-t,0]}\varphi_1(x+t)d^-\eta(x),\ldots,\int_{[-t,0]}\varphi_N(x+t)d^-\eta(x)\bigg), \\
F(t,\eta,y,z) \ &= \ \bar F\bigg(t,\int_{[-t,0]}\varphi_1(x+t)d^-\eta(x),\ldots,\int_{[-t,0]}\varphi_N(x+t)d^-\eta(x),y,z\bigg), \\
{H}(\eta) \ &= \ \bar{H}\bigg(\int_{[-T,0]}\varphi_1(x+T)d^-\eta(x),\ldots,\int_{[-T,0]}\varphi_N(x+T)d^-\eta(x)\bigg),
\end{align*}
where (we refer to Definition 2.4(i) in the companion paper \cite{cosso_russo15a} for a definition of the forward integral with respect to $\eta$) the following assumptions are made.
\begin{itemize}
\item[\textup{(i)}] $\bar b$, $\bar\sigma$, $\bar F$, $\bar{H}$ are continuous and satisfy Assumption {\bf (A0)}.
\item[\textup{(ii)}] $\bar b$ and $\bar\sigma$ are of class $C^3$ with partial derivatives from order $1$ up to order $3$ bounded.
\item[\textup{(iii)}] For all $t\in[0,T]$, $\bar F(t,\cdot,\cdot,\cdot)\in C^3(\R^{N+2})$ and moreover 
we assume the validity of the properties below.
\begin{enumerate}
\item[\textup{(a)}] $\bar F(t,\cdot,0,0)$ belongs to $C^3$ and its third order partial derivatives satisfy a polynomial growth condition uniformly in $t$.
\item[\textup{(b)}] $D_y\bar F$, $D_z\bar F$ are bounded on $[0,T]\times\R^N\times\R\times\R$, as well as their derivatives of order one and second with respect to $x_1,\ldots,x_N,y,z$.
\end{enumerate}
\item[\textup{(iv)}] $\bar H\in C^3(\R^N)$ and its third order partial derivatives satisfy a polynomial growth condition.
\item[\textup{(v)}] $\varphi_1,\ldots,\varphi_N\in C^2([0,T])$.
\end{itemize}
Then, the map $\Uc$ given by
\[
\Uc(t,\eta) \ = \ Y_t^{t,\eta}, \qquad \forall\,(t,\eta)\in[0,T]\times C([-T,0]),
\]
where $(Y_s^{t,\eta},Z_s^{t,\eta})_{s\in[t,T]}\in\S^2(t,T)\times\H^2(t,T)$ is the unique solution to \eqref{BSDE_SV}, is a strict solution to equation \eqref{KolmEq}.
\end{Theorem}
\textbf{Proof.}
See Theorem 3.6 in \cite{cosso_russo15a}.
\ep

\begin{Remark}
\label{R:b_sigma_time-dep}
{\rm
Notice that in Theorem \ref{T:ExistenceStrict} the functions $\bar b$ and $\bar\sigma$ do not depend on time. For the case where $\bar b$ and $\bar\sigma$ are time-dependent, we refer to Theorem 3.5 in \cite{cosso_russo15a}. Notice that, in this case, $F=F(t,\eta)$ does not depend on $(y,z)$.
\ep
}
\end{Remark}

\subsection{Strong-viscosity solutions}

In the present section, we introduce the notion of (path-dependent) strong-viscosity solution to equation \eqref{KolmEq}. To do it, we extend in a natural way Definition \ref{D:ViscosityFinite} to the present path-dependent case, see also Remark \ref{R:TwoDefinitions}.

\begin{Remark}
\label{R:MotivationSV}
{\rm
As a motivation for the introduction of a viscosity type solution for path-dependent PDEs, let us consider the following hedging example in mathematical finance, taken from Section 3.2 in the survey paper \cite{cosso_russo14}. Let $b\equiv 0$, $\sigma\equiv 1$, $F\equiv0$ and consider the lookback-type payoff
\[
{H}(\eta) \ = \ \sup_{x\in[-T,0]} \eta(x), \qquad \forall\,\eta \in C([-T,0]).
\]
Then, we look for a solution to the following linear parabolic path-dependent PDE:
\begin{equation}
\label{KolmEqLinear}
\begin{cases}
\partial_t \Uc + D^H \Uc + \frac{1}{2}D^{VV} \Uc \ = \ 0, \qquad &\quad\forall\,(t,\eta)\in[0,T[\times C([-T,0]), \\
\Uc(T,\eta) \ = \ {H}(\eta), &\quad\forall\,\eta\in C([-T,0]).
\end{cases}
\end{equation}
We refer to \eqref{KolmEqLinear} as path-dependent heat equation. Notice that, however, \eqref{KolmEqLinear} does not have the smoothing effect characterizing the classical heat equation, in spite of some regularity properties illustrated in Section 3.2 of \cite{cosso_russo14}. Indeed, let us consider the functional
\[
\Uc(t,\eta) \ = \ \E\big[{H}(\mathbb W_T^{t,\eta})\big] \ = \E\Big[\sup_{-T\leq x\leq0}\mathbb W_T^{t,\eta}(x)\Big], \qquad \forall\,(t,\eta)\in[0,T]\times C([-T,0]),
\]
where, for any $t \leq s \leq T$,
\[
\mathbb W_s^{t,\eta}(x) \ = \
\begin{cases}
\eta(x+s-t), &\quad-T \leq x \leq t-s, \\
\eta(0) + W_{x+s} - W_t, \qquad &\quad t-s < x \leq 0.
\end{cases}
\]
If $\Uc\in C^{1,2}(([0,T[\times\textup{past})\times\textup{present})\cap C([0,T]\times C([-T,0]))$, then $\Uc$ could be proved to solve equation \eqref{KolmEqLinear}. However, as claimed in \cite{cosso_russo14}, $\Uc$ is not a strict solution to \eqref{KolmEqLinear}. On the other hand, since ${H}$ is continuous and has linear growth, it follows from Theorems \ref{T:UniqSV} and \ref{T:ExistSV} that $\Uc$ is the unique strong-viscosity solution to equation \eqref{KolmEqLinear}.
\ep
}
\end{Remark}

\begin{Definition}
\label{D:Strong-Visc}
A function $\Uc\colon[0,T]\times C([-T,0])\rightarrow\R$ is called a \textbf{(path-dependent) strong-viscosity solution} to equation \eqref{KolmEq} if there exists a sequence $(\Uc_n,{H}_n,F_n,b_n,\sigma_n)_n$ of Borel measurable functions $\Uc_n\colon[0,T]\times C([-T,0])\rightarrow\R$, ${H}_n\colon C([-T,0])\rightarrow\R$, $F_n\colon[0,T]\times C([-T,0])\times\R\times\R\rightarrow\R$, $b_n\colon[0,T]\times C([-T,0])\rightarrow\R$, $\sigma_n\colon[0,T]\times C([-T,0])\rightarrow\R$, such that the following holds.
\begin{itemize}
\item[\textup{(i)}] For some positive constants $C$ and $m$,
\begin{align*}
|b_n(t,\eta)| + |\sigma_n(t,\eta)| \ &\leq \ C(1 + \|\eta\|), \\
|b_n(t,\eta) - b_n(t,\eta')| + |\sigma_n(t,\eta) - \sigma_n(t,\eta')| \ &\leq \ C\|\eta - \eta'\|, \\
|F_n(t,\eta,y,z) - F_n(t,\eta,y',z')| \ &\leq \ C\big(|y-y'| + |z-z'|\big), \\
|{H}_n(\eta)| + |F_n(t,\eta,0,0)| + |\Uc_n(t,\eta)| \ &\leq \ C\big(1 + \|\eta\|^m\big),
\end{align*}
for all $t\in[0,T]$, $\eta,\eta'\in C([-T,0])$, $y,y',z,z'\in\R$. Moreover, the functions $\Uc_n(t,\cdot)$, ${H}_n(\cdot)$, $F_n(t,\cdot,\cdot,\cdot)$, $n\in\N$, are equicontinuous on compact sets, uniformly with respect to $t\in[0,T]$.
\item[\textup{(ii)}] $\Uc_n$ is a strict solution to
\[
\begin{cases}
\partial_t\Uc_n + D^H\Uc_n + b_n(t,\eta)D^V\Uc_n + \frac{1}{2}\sigma_n(t,\eta)^2D^{VV}\Uc_n \\
\hspace{2.5cm} +\, F_n(t,\eta,\Uc_n,\sigma_n(t,\eta)D^V\Uc_n) \ = \ 0, \;\;\; &\;\;\forall\,(t,\eta)\in[0,T[\times C([-T,0]), \\
\Uc_n(T,\eta) \ = \ {H}_n(\eta), &\;\;\forall\,\eta\in C([-T,0]).
\end{cases}
\]
\item[\textup{(iii)}] $(\Uc_n,{H}_n,F_n,b_n,\sigma_n)_n$ converges pointwise to $(\Uc,{H},F,b,\sigma)$ as $n\rightarrow\infty$.
\end{itemize}
\end{Definition}

We present a Feynman-Kac type representation for a generic strong-viscosity solution to equation \eqref{KolmEq}, which, as a consequence, yields a uniqueness
 result.

\begin{Theorem}
\label{T:UniqSV}
Let Assumption {\bf (A1)} hold and let $\Uc\colon[0,T]\times C([-T,0])\rightarrow\R$ be a strong-viscosity solution to equation \eqref{KolmEq}. Then, the following Feynman-Kac formula holds
\begin{equation}
\label{Feynman-Kac}
\Uc(t,\eta) \ = \ Y_t^{t,\eta}, \qquad \forall\,(t,\eta)\in[0,T]\times C([-T,0]),
\end{equation}
where $(Y_s^{t,\eta},Z_s^{t,\eta})_{s\in[t,T]}\in\S^2(t,T)\times\H^2(t,T)$, with $Y_s^{t,\eta}=\Uc(s,\X_s^{t,\eta})$, is the unique solution in $\S^2(t,T)\times\H^2(t,T)$ to the backward stochastic differential equation: $\P$-a.s.
\begin{equation}
\label{BSDE_SV}
Y_s^{t,\eta} \ = \ {H}(\X_T^{t,\eta}) + \int_s^T F(r,\X_r^{t,\eta},Y_r^{t,\eta},Z_r^{t,\eta}) dr - \int_s^T Z_r^{t,\eta} dW_r, \qquad t\leq s\leq T.
\end{equation}
In particular, there exists at most one strong-viscosity solution to equation \eqref{KolmEq}. 
\end{Theorem}
\textbf{Proof.}
Let $(\Uc_n,{H}_n,F_n,b_n,\sigma_n)_n$ be as in Definition \ref{D:Strong-Visc} and, for any $(t,\eta)\in[0,T]\times C([-T,0])$, denote by $X^{n,t,\eta}=(X_s^{n,t,\eta})_{s\in[t,T]}$ the unique solution to equation \eqref{Xn}. Then, from Theorem \ref{T:UniqStrict}, $(Y_s^{n,t,\eta},Z_s^{n,t,\eta})_{s\in[t,T]}=(\Uc_n(s,\X_s^{n,t,\eta}),\sigma_n(s,\X_s^{n,t,\eta})D^V\Uc_n(s,\X_s^{n,t,\eta})1_{[t,T[}(s))_{s\in[t,T]}$ is the unique solution to the backward stochastic differential equation: $\P$-a.s.,
\[
Y_s^{n,t,\eta} \ = \ {H}_n(\X_T^{n,t,\eta}) + \int_s^T F_n(r,\X_r^{n,t,\eta},Y_r^{n,t,\eta},Z_r^{n,t,\eta}) dr - \int_s^T Z_r^{n,t,\eta} dW_r, \qquad t\leq s\leq T.
\]
We wish now to take the limit when $n$ goes to infinity in the above equation. We make use of Theorem \ref{T:LimitThmBSDE}, for which we check the assumptions. From the polynomial growth condition of $\Uc_n$ together with estimate \eqref{EstimateSupXn}, there exists, for every $p\geq1$, a constant $\tilde C_p\geq0$ such that
\begin{equation}
\label{EstimateSupYn}
\big\|Y^{n,t,\eta}\big\|_{_{\S^p(t,T)}}^p \ \leq \ \tilde C_p \big(1 + \|\eta\|^p\big), \qquad \forall\,n\in\N.
\end{equation}
Now, from Proposition \ref{P:EstimateBSDEAppendix}, it follows that there exists a constant $\tilde c\geq0$ (depending only on $T$ and on the Lipschitz constant $C$ of $F_n$ with respect to $(y,z)$ appearing in Definition \ref {D:Strong-Visc}(i)) such that
\[
\big\|Z^{n,t,\eta}\big\|_{_{\H^2(t,T)}}^2 \ \leq \ \tilde c \bigg(\big\|Y^{n,t,\eta}\big\|_{_{\S^2(t,T)}}^2 + \E \int_t^T |F_n(s,\X_s^{n,t,\eta},0,0)|^2 ds\bigg).
\]
Therefore, from \eqref{EstimateSupYn}, the polynomial growth condition of $F_n$, and estimate \eqref{EstimateSupXn}, we find that $\sup_n\|Z^{n,t,\eta}\|_{_{\H^2(t,T)}}^2<\infty$. Moreover, from \eqref{limXn-X} we see that, for any $s\in[t,T]$, $\|\X_s^{n,t,\eta}(\omega)-\X_s^{t,\eta}(\omega)\|\rightarrow0$, as $n\rightarrow\infty$, for $\P$-a.e. $\omega\in\Omega$. Fix such an $\omega$ and consider the set $K_\omega\subset C([-T,0])$ given by
\[
K_\omega \ := \ \big(\cup_{n\in\N} \big\{\X_s^{n,t,\eta}(\omega)\big\}\big) \cup \big\{\X_s^{t,\eta}(\omega)\big\}.
\]
Then, $K_\omega$ is a compact subset of $C([-T,0])$. Since the sequence $(F_n(s,\cdot,\cdot,\cdot))_n$ is equicontinuous on compact sets and converges pointwise to $F(s,\cdot,\cdot,\cdot)$, it follows that $(F_n(s,\cdot,\cdot,\cdot))_n$ converges to $F(s,\cdot,\cdot,\cdot)$ uniformly on compact sets. In particular, we have
\begin{align*}
&\big|F_n(s,\X_s^{n,t,\eta}(\omega),0,0) - F(s,\X_s^{t,\eta}(\omega),0,0)\big| \\
&\leq \ \sup_{\eta\in K_\omega} \big|F_n(s,\eta,0,0) - F(s,\eta,0,0)\big| + \big|F(s,\X_s^{n,t,\eta}(\omega),0,0) - F(s,\X_s^{t,\eta}(\omega),0,0)\big| \ \overset{n\rightarrow\infty}{\longrightarrow} \ 0.
\end{align*}
Similarly, we have
\begin{align*}
&\big|\Uc_n(s,\X_s^{n,t,\eta}(\omega)) - \Uc(s,\X_s^{t,\eta}(\omega))\big| \\
&\leq \ \sup_{\eta\in K_\omega} \big|\Uc_n(s,\eta) - \Uc(s,\eta)\big| + \big|\Uc(s,\X_s^{n,t,\eta}(\omega)) - \Uc(s,\X_s^{t,\eta}(\omega))\big| \ \overset{n\rightarrow\infty}{\longrightarrow} \ 0.
\end{align*}
Let us now define $Y_s^{t,\eta}:=\Uc(s,\X_s^{t,\eta})$, for all $s\in[t,T]$. We can then apply Theorem \ref{T:LimitThmBSDE} (notice that, in this case, for every $n\in\N$, the process $K^n$ appearing in Theorem \ref{T:LimitThmBSDE} is identically zero, so that $K$ is also identically zero), from which it follows that there exists $Z^{t,\eta}\in\H^2(t,T)$ such that the pair $(Y^{t,\eta},Z^{t,\eta})$ solves equation \eqref{BSDE_SV}. From Theorem 3.1 in \cite{parpen90} we have that $(Y^{t,\eta},Z^{t,\eta})$ is the unique pair in $\S^2(t,T)\times\H^2(t,T)$ satisfying equation \eqref{BSDE_SV}. This concludes the proof.
\ep

\vspace{3mm}

By Theorem \ref{T:UniqSV} we deduce Lemma \ref{L:Stability} below, which says that in Definition \ref{D:Strong-Visc} the convergence of $(\Uc_n)_n$ is indeed a consequence of the convergence of the coefficients $({H}_n,F_n,b_n,\sigma_n)_n$. This result is particularly useful to establish the existence of strong-viscosity solutions, as in the proof of Theorem \ref{T:ExistSV}.

\begin{Lemma}
\label{L:Stability}
Suppose that Assumption {\bf (A1)} holds and let $(\Uc_n,{H}_n,F_n,b_n,\sigma_n)_n$ be as in Definition \ref{D:Strong-Visc}, except that
 we do not assume the convergence of $(\Uc_n)_n$. Then, there exists $\Uc\colon[0,T]\times C([-T,0])\rightarrow\R$ such that $(\Uc_n)_n$ converges pointwise to $\Uc$. In particular, $\Uc$ is a strong-viscosity solution to equation \eqref{KolmEq} and is given by formula \eqref{Feynman-Kac}.
\end{Lemma}
\textbf{Proof.}
Let us prove the pointwise convergence of the sequence $(\Uc_n)_{n\in\N}$ to the function $\Uc$ given by formula \eqref{Feynman-Kac}. To this end, we notice that, from Theorem \ref{T:UniqStrict}, for every $n\in\N$, $\Uc_n$ is given by
\[
\Uc_n(t,\eta) \ = \ Y_t^{n,t,\eta}, \qquad \forall\,(t,\eta)\in[0,T]\times C([-T,0]),
\]
where $(Y^{n,t,\eta},Z^{n,t,\eta}) = (\Uc_n(\cdot,\mathbb X^{n,t,\eta}),\sigma_n(\cdot,\mathbb X^{n,t,\eta})D^V\Uc_n(\cdot,\mathbb X^{n,t,\eta})1_{[t,T[})\in\S^2(t,T)\times\H^2(t,T)$ is the unique solution to the backward stochastic differential equation: $\P$-a.s.,
\[
Y_s^{n,t,\eta} \ = \ {H}_n(\mathbb X_T^{n,t,\eta}) + \int_s^T F_n(r,\mathbb X_r^{n,t,\eta},Y_r^{n,t,\eta},Z_r^{n,t,\eta}) dr - \int_s^T Z_r^{n,t,\eta} dW_r, \quad t \leq s \leq T,
\]
with
\[
X_s^{n,t,\eta} = \ \eta(0\wedge(s-t)) + \int_t^{t\vee s} b_n(r,\mathbb X_r^{n,t,\eta})dr + \int_t^{t\vee s} \sigma_n(r,\mathbb X_r^{n,t,\eta})dW_r, \;\; - T + t \leq s \leq T.
\]
Consider the function $\Uc$ given by formula \eqref{Feynman-Kac}. From estimate \eqref{EstimateBSDE2}, there exists a constant $C$, independent of $n\in\N$, such that
\begin{align*}
|\Uc_n(t,\eta) - \Uc(t,\eta)|^2 \ &\leq \ C\,\E\big[\big|H_n(\X_T^{n,t,\eta}) - H(\X_T^{t,\eta})\big|^2\big] \\
&\quad \ + C\int_t^T \E\big[\big|F_n(s,\X_s^{n,t,\eta},Y_s^{t,\eta},Z_s^{t,\eta}) - F(s,\X_s^{t,\eta},Y_s^{t,\eta},Z_s^{t,\eta})\big|^2\big] ds, \notag
\end{align*}
for all $t\in[0,T]$ and $\eta\in C([-T,0])$. Now we recall the following.
\begin{itemize}
\item[(i)] $(H_n,F_n,b_n,\sigma_n)_{n\in\N}$ converges pointwise to $({H},F,b,\sigma)$ as $n\rightarrow\infty$.
\item[(ii)] The functions $H_n(\cdot),F_n(t,\cdot,\cdot,\cdot),b_n(t,\cdot),\sigma_n(t,\cdot)$, $n\in\N$, are equicontinuous on compact sets, uniformly with respect to $t\in[0,T]$.
\end{itemize}
We notice that (i) and (ii) imply the following property:
\begin{itemize}
\item[(iii)] $(H_n(\eta_n),F_n(t,\eta_n,y,z),b_n(t,\eta_n),\sigma_n(t,\eta_n))$ converges to $(H(\eta),F(t,\eta,y,z),b(t,\eta),\sigma(t,\eta))$ as $n\rightarrow\infty$, $\forall\,(t,y,z)\in[0,T]\times\R\times\R$, $\forall\,(\eta_n)_{n\in\N}\subset C([-T,0])$ with $\eta_n\rightarrow\eta\in C([-T,0])$.
\end{itemize}
Let us now remind that, for any $r\in[t,T]$, we have
\[
\X_s^{n,t,\eta}(x) \ =
\begin{cases}
\eta(s-t+x), \qquad & x\in[-T,t-s], \\
X_{s+x}^{n,t,\eta}, & x\in\;]t-s,0],
\end{cases}
\quad
\X_s^{t,\eta}(x) \ =
\begin{cases}
\eta(s-t+x), \qquad & x\in[-T,t-s], \\
X_{s+x}^{t,\eta}, & x\in\;]t-s,0].
\end{cases}
\]
Therefore, for every $p\geq1$,
\begin{equation}
\label{ConvX}
\E\Big[\sup_{t\leq s\leq T} \|\X_s^{n,t,\eta} - \X_s^{t,\eta}\|_\infty^p\Big] \ = \ \E\Big[\sup_{t\leq s\leq T} |X_s^{n,t,\eta} - X_s^{t,\eta}|^p\Big] \ \overset{n\rightarrow\infty}{\longrightarrow} \ 0,
\end{equation}
where the convergence follows from \eqref{limXn-X}. Then, we claim that the following convergences in probability hold:
\begin{align}
&\big|H_n(\X_T^{n,t,\eta}) - H(\X_T^{t,\eta})\big|^2 \ \overset{\P}{\underset{n\rightarrow\infty}{\longrightarrow}} \ 0, \label{ConvH} \\
&\big|F_n(s,\X_s^{n,t,\eta},Y_s^{t,\eta},Z_s^{t,\eta}) - F(s,\X_s^{t,\eta},Y_s^{t,\eta},Z_s^{t,\eta})\big|^2 \ \overset{\P}{\underset{n\rightarrow\infty}{\longrightarrow}} \ 0, \label{ConvF}
\end{align}
for all $s\in[t,T]$. Concerning \eqref{ConvH}, we begin noting that it is enough to prove that, for every subsequence $(|H_{n_m}(\X_T^{n_m,t,\eta}) - H(\X_T^{t,\eta})|^2)_{m\in\N}$ there exists a subsubsequence which converges to zero. From \eqref{ConvX} and property (iii) above, it follows that there exists a subsubsequence $(|H_{n_{m_\ell}}(\X_T^{n_{m_\ell},t,\eta}) - H(\X_T^{t,\eta})|^2)_{\ell\in\N}$ which converges $\P$-a.s., and therefore in probability, to zero. This concludes the proof of \eqref{ConvH}. In a similar way we can prove \eqref{ConvF}.

From \eqref{ConvH} and \eqref{ConvF}, together with the uniform integrability of the sequences $(|H_n(\X_T^{n,t,\eta}) - H(\X_T^{t,\eta})|^2)_{n\in\N}$ and $(|F_n(s,\X_s^{n,t,\eta},Y_s^{t,\eta},Z_s^{t,\eta}) - F(s,\X_s^{t,\eta},Y_s^{t,\eta},Z_s^{t,\eta})|^2)_{n\in\N}$, for every $s\in[t,T]$, we deduce that
\begin{align*}
&\lim_{n\rightarrow\infty} \E\big[\big|H_n(\X_T^{n,t,\eta}) - H(\X_T^{t,\eta})\big|^2\big] \ = \ 0, \\
&\lim_{n\rightarrow\infty} \E\big[\big|F_n(s,\X_s^{n,t,\eta},Y_s^{t,\eta},Z_s^{t,\eta}) - F(s,\X_s^{t,\eta},Y_s^{t,\eta},Z_s^{t,\eta})\big|^2\big] \ = \ 0.
\end{align*}
From the second convergence, the polynomial growth condition of $F$ and $F_n$ (uniform in $n$), and standard moment estimates for $\|\X^{n,t,\eta}\|_\infty\leq\sup_{t\leq s\leq T}|X_s^{n,t,\eta}|$ (see estimate \eqref{EstimateSupXn}), it follows that
\[
\lim_{n\rightarrow\infty} \int_t^T \E\big[\big|F_n(s,\X_s^{n,t,\eta},Y_s^{t,\eta},Z_s^{t,\eta}) - F(s,\X_s^{t,\eta},Y_s^{t,\eta},Z_s^{t,\eta})\big|^2\big] ds \ = \ 0.
\]
As a consequence, we have $|\Uc_n(t,\eta)-\Uc(t,\eta)|^2\rightarrow0$ as $n\rightarrow\infty$, which concludes the proof.
\ep

\vspace{3mm}

We can now state an existence result. Notice that it holds under quite general conditions on the terminal condition $H$ of equation \eqref{KolmEq}.

\begin{Theorem}
\label{T:ExistSV}
Let Assumption {\bf (A1)} hold and suppose that ${H}$ is continuous. Suppose also that there exists a nondecreasing sequence $(N_n)_{n\in\N}\subset\N\backslash\{0\}$ such that, for all $n\in\N$ and $(t,\eta,y,z)\in[0,T]\times C([-T,0])\times\R\times\R$,
\begin{align*}
b_n(t,\eta) \ &= \ \bar b_n\bigg(\int_{[-t,0]}\varphi_1(x+t)d^-\eta(x),\ldots,\int_{[-t,0]}\varphi_{N_n}(x+t)d^-\eta(x)\bigg), \\
\sigma_n(t,\eta) \ &= \ \bar\sigma_n\bigg(\int_{[-t,0]}\varphi_1(x+t)d^-\eta(x),\ldots,\int_{[-t,0]}\varphi_{N_n}(x+t)d^-\eta(x)\bigg), \\
F_n(t,\eta,y,z) \ &= \ \bar F_n\bigg(t,\int_{[-t,0]}\varphi_1(x+t)d^-\eta(x),\ldots,\int_{[-t,0]}\varphi_{N_n}(x+t)d^-\eta(x),y,z\bigg),
\end{align*}
where the following holds.
\begin{itemize}
\item[\textup{(i)}] $\bar b_n$, $\bar\sigma_n$, $\bar F_n$ are continuous and satisfy Assumption {\bf (A0)} with constants $C$ and $m$ independent of $n$.
\item[\textup{(ii)}] For every $n\in\N$, $\bar b_n,\bar\sigma_n,\bar F_n$ satisfy items (ii) and (iii) of Theorem \ref{T:ExistenceStrict}.
\item[\textup{(iii)}] The functions $b_n(t,\cdot)$, $\sigma_n(t,\cdot)$, $F_n(t,\cdot,\cdot,\cdot)$, $n\in\N$, are equicontinuous on compact sets, uniformly with respect to $t\in[0,T]$.
\item[\textup{(iv)}] $\varphi_1,\ldots,\varphi_{N_n}\in C^2([0,T])$ are uniformly bounded with respect to $n\in\N$, and their first derivative are bounded in $L^1([0,T])$ uniformly with respect to $n\in\N$.
\item[\textup{(v)}] $(b_n,\sigma_n,F_n)_n$ converges pointwise to $(b,\sigma,F)$ as $n\rightarrow\infty$.
\end{itemize}
Then, the map $\Uc$ given by
\begin{equation}
\label{Feynman-Kac2}
\Uc(t,\eta) \ = \ Y_t^{t,\eta}, \qquad \forall\,(t,\eta)\in[0,T]\times C([-T,0]),
\end{equation}
where $(Y_s^{t,\eta},Z_s^{t,\eta})_{s\in[t,T]}\in\S^2(t,T)\times\H^2(t,T)$ is the unique solution to \eqref{BSDE_SV}, is a (path-dependent) strong-viscosity solution to equation \eqref{KolmEq}.
\end{Theorem}
\textbf{Proof.} We divide the proof into four steps. In the first three steps we construct an approximating sequence of smooth functions for $H$. We conclude the proof in the fourth step.

\vspace{1mm}

\noindent\textbf{Step I.} \emph{Approximation of $\eta\in C([-t,0])$, $t\in\,]0,T]$, with Fourier partial sums.} Consider the sequence $(e_i)_{i\in\N}$ of $C^{\infty}([-T,0])$ functions:
\[
e_0 = \frac{1}{\sqrt{T}}, \qquad e_{2i-1}(x) = \sqrt{\frac{2}{T}}\sin\bigg(\frac{2\pi}{T}(x+T)i\bigg), \qquad e_{2i}(x) = \sqrt{\frac{2}{T}}\cos\bigg(\frac{2\pi}{T}(x+T)i\bigg),
\]
for all $i\in\N\backslash\{0\}$. Then $(e_i)_{i\in\N}$ is an orthonormal basis of $L^2([-T,0])$. Let us define the linear operator $\Lambda\colon C([-T,0])\rightarrow C([-T,0])$ by
\[
(\Lambda\eta)(x) \ = \ \frac{\eta(0)-\eta(-T)}{T}x, \qquad x\in[-T,0],\,\eta\in C([-T,0]).
\]
Notice that $(\eta-\Lambda\eta)(-T) = (\eta-\Lambda\eta)(0)$, therefore $\eta-\Lambda\eta$ can be extended to the entire real line in a periodic way with period $T$, so that we can expand it in Fourier series. In particular, for each $n\in\N$ and $\eta\in C([-T,0])$, consider the Fourier partial sum
\begin{equation}
\label{E:s_n}
s_n(\eta-\Lambda\eta) \ = \ \sum_{i=0}^n (\eta_i-(\Lambda\eta)_i) e_i, \qquad \forall\,\eta\in C([-T,0]),
\end{equation}
where (denoting $\tilde e_i(x) = \int_{-T}^x e_i(y) dy$, for any $x\in[-T,0]$), by the integration by parts formula (2.4) of \cite{cosso_russo15a},
\begin{align}
\label{E:eta_i}
\eta_i \ = \ \int_{-T}^0 \eta(x)e_i(x) dx \ &= \ \eta(0)\tilde e _i(0) - \int_{[-T,0]} \tilde e_i(x)d^-\eta(x) = \ \int_{[-T,0]} (\tilde e_i(0) - \tilde e_i(x)) d^- \eta(x),
\end{align}
since $\eta(0) = \int_{[-T,0]}d^-\eta(x)$. Moreover we have
\begin{align}
\label{E:Lambda_eta_i}
(\Lambda\eta)_i \ &= \ \int_{-T}^0 (\Lambda\eta)(x)e_i(x) dx \ = \ \frac{1}{T} \int_{-T}^0 xe_i(x) dx \bigg(\int_{[-T,0]} d^-\eta(x) - \eta(-T)\bigg).
\end{align}
Define $\sigma_n=\frac{s_0 + s_1 + \cdots + s_n}{n+1}$. Then, by \eqref{E:s_n},
\[
\sigma_n(\eta-\Lambda\eta) \ = \ \sum_{i=0}^n \frac{n+1-i}{n+1} (\eta_i-(\Lambda\eta)_i) e_i, \qquad \forall\,\eta\in C([-T,0]).
\]
We know from Fej\'er's theorem on Fourier series (see, e.g., Theorem 3.4, Chapter III, in \cite{zygmund02}) that, for any $\eta\in C([-T,0])$, $\sigma_n(\eta-\Lambda\eta)\rightarrow\eta-\Lambda\eta$ uniformly on $[-T,0]$, as $n$ tends to infinity, and $\|\sigma_n(\eta-\Lambda\eta)\|_\infty \leq \|\eta-\Lambda\eta\|_\infty$. Let us define the linear operator $T_n\colon C([-T,0])\rightarrow C([-T,0])$ by (denoting $e_{-1}(x) = x$, for any $x\in[-T,0]$)
\begin{align*}
T_n\eta \ = \ \sigma_n(\eta-\Lambda\eta) + \Lambda\eta \ &= \ \sum_{i=0}^n \frac{n+1-i}{n+1} (\eta_i-(\Lambda\eta)_i) e_i + \frac{\eta(0) - \eta(-T)}{T}e_{-1} \notag \\
&= \ \sum_{i=0}^n \frac{n+1-i}{n+1} y_ie_i + y_{-1}e_{-1},
\end{align*}
where, using \eqref{E:eta_i} and \eqref{E:Lambda_eta_i},
\begin{align*}
y_{-1} \ &= \int_{[-T,0]} \frac{1}{T} d^-\eta(x) - \frac{1}{T}\eta(-T), \\
y_i \ &= \ \int_{[-T,0]}\bigg (\tilde e_i(0)-\tilde e_i(x)- \frac{1}{T}\int_{-T}^0 xe_i(x)dx\bigg )d^-\eta(x) + \frac{1}{T}\int_{-T}^0 xe_i(x)dx\,\eta(-T),
\end{align*}
for $i=0,\ldots,n$. Then, for any $\eta\in C([-T,0])$, $T_n\eta\rightarrow\eta$ uniformly on $[-T,0]$, as $n$ tends to infinity. Furthermore, there exists a positive constant $M$ such that
\begin{equation}
\label{E:UniformBoundT_n}
\|T_n\eta\|_\infty \ \leq \ M\|\eta\|_\infty, \qquad \forall\,n\in\N,\,\forall\,\eta\in C([-T,0]).
\end{equation}
Then, we define
\[
\tilde{H}_n(\eta) \ := \ {H}(T_n\eta), \qquad \forall\,(t,\eta)\in[0,T]\times C([-T,0]).
\]
Notice that $\tilde{H}_n$ satisfies a polynomial growth condition as in Assumption {\bf (A1)} with constants $C$ and $m$ independent of $n$. Moreover, since ${H}$ is uniformly continuous on compact sets, from \eqref{E:UniformBoundT_n} we see that $(\tilde{H}_n)_n$ is equicontinuous on compact sets. Now, we define the function $\bar{H}_n\colon\R^{n+2}\rightarrow\R$ as follows
\[
\bar{H}_n(y_{-1},\ldots,y_n) \ := \ {H}\bigg(\sum_{i=0}^n \frac{n+1-i}{n+1} y_ie_i + y_{-1}e_{-1}\bigg), \qquad \forall\,(y_{-1},\ldots,y_n)\in\R^{n+2}.
\]
Then, we have
\[
\tilde{H}_n(\eta) = \bar{H}_n\bigg(\int_{[-T,0]} \!\! \psi_{-1}(x + T) d^-\eta(x) + a_{-1}\eta(-T),\ldots,\int_{[-T,0]} \!\! \psi_n(x + T) d^-\eta(x) + a_n\eta(-T)\bigg),
\]
for all $\eta\in C([-T,0])$, $n\in\N$, where
\begin{align*}
\psi_{-1}(x) \ &= \ \frac{1}{T}, \qquad \psi_i(x) \ = \ \tilde e_i(0) - \tilde e_i(x - T) - \frac{1}{T} \int_{-T}^0 xe_i(x)dx, \qquad x\in[-T,0], \\
a_{-1} \ &= \ - \frac{1}{T}, \qquad\;\;\, a_i \ = \ \frac{1}{T} \int_{-T}^0 xe_i(x)dx.
\end{align*}

\vspace{1mm}

\noindent\textbf{Step II.} \emph{Smoothing of $\eta(-T)$ through mollifiers.} Consider the function $\phi\in C^\infty([0,\infty[)$ given by
\[
\phi(x) \ = \ c \exp\bigg(\frac{1}{x^2 - T^2}\bigg) 1_{[0,T[}(x), \qquad \forall\,x\geq0,
\]
with $c>0$ such that $\int_0^\infty \phi(x) dx=1$. Then, we define $\phi_m(x)=m\phi(mx)$, $\forall\,x\geq0$, $m\in\N$. Notice that
\begin{align*}
\int_{-T}^0 \eta(x)\phi_m(x + T) dx \ &= \ \eta(0)\tilde\phi_m(T) - \int_{[-T,0]} \tilde\phi_m(x + T) d^-\eta(x) \\
&= \ \int_{[-T,0]} \big(\tilde\phi_m(T) - \tilde\phi_m(x + T)\big) d^-\eta(x),
\end{align*}
where $\tilde\phi_m(x)=\int_0^x \phi_m(z) dz$, $x\in[0,T]$. In particular, we have
\[
\lim_{m\rightarrow\infty} \int_{[-T,0]} \big(\tilde\phi_m(T) - \tilde\phi_m(x + T)\big) d^-\eta(x) \ = \ \lim_{m\rightarrow\infty} \int_{-T}^0 \eta(x)\phi_m(x + T) dx \ = \ \eta(-T).
\]
Then, we define
\begin{align}
\label{barHn=H}
H_n(\eta) \ &:= \ \bar{H}_n\bigg(\ldots,\int_{[-T,0]} \psi_i(x + T) d^-\eta(x) + a_i\int_{[-T,0]} \big(\tilde\phi_n(T) - \tilde\phi_n(x + T)\big) d^-\eta(x),\ldots\bigg) \notag \\
&= \ {H}\bigg(T_n\eta + \bigg(\sum_{i=0}^n \frac{n+1-i}{n+1} a_i e_i + a_{-1} e_{-1}\bigg)\int_{-T}^0 \big(\eta(x) - \eta(-T)\big) \phi_n(x + T) dx\bigg) \notag \\
&= \ {H}\bigg(T_n\eta + \bigg(T_n\gamma + \frac{1}{T(T-1)} e_{-1}\bigg)\int_{-T}^0 \big(\eta(x) - \eta(-T)\big) \phi_n(x + T) dx\bigg),
\end{align}
for all $\eta\in C([-T,0])$ and $n\in\N$, where $\gamma(x):=-x/(T-1)$, $\forall\,x\in[-T,0]$. Then, the sequence $(H_n)_n$ is equicontinuous on compact sets and converges pointwise to ${H}$ as $n\rightarrow\infty$.

\vspace{1mm}

\noindent\textbf{Step III.} \emph{Smoothing of $\bar{H}_n(\cdot)$.} From \eqref{barHn=H} it follows that for any compact subset $K\subset C([-T,0])$ there exists a continuity modulus $m_K$, independent of $n\in\N$, such that
\begin{align}
\label{UniformContinuitybarH_n}
&\bigg|\bar{H}_n\bigg(\ldots,\int_{[-T,0]} \psi_i(x + T) d^-\eta_1(x) + a_i\int_{[-T,0]} \big(\tilde\phi_n(T) - \tilde\phi_n(x + T)\big) d^-\eta_1(x) + \xi_i,\ldots\bigg) \notag \\
&- \bar{H}_n\bigg(\ldots,\int_{[-T,0]} \psi_i(x + T) d^-\eta_2(x) + a_i\int_{[-T,0]} \big(\tilde\phi_n(T) - \tilde\phi_n(x + T)\big) d^-\eta_2(x) + \xi_i,\ldots\bigg)\bigg| \notag \\
&\leq \ m_K(\|\eta_1-\eta_2\|_\infty),
\end{align}
for all $\eta_1,\eta_2\in K$, $n\in\N$, $\xi=(\xi_{-1},\ldots,\xi_n)\in E_{n+2}$, where $E_{n+2}:=\{\xi=(\xi_{-1},\ldots,\xi_n)\in\R^{n+2}\colon|\xi_i|\leq 2^{-(i+1)},\,i=-1,\ldots,n\}$. Indeed, set
\[
\mathcal K \ := \ K\cup\tilde K,
\]
where
\begin{align*}
\tilde K \ := \ \bigg\{\eta\in C([-T,0])\colon\eta \ &= \ T_n\eta_1 + \Big(T_n\gamma + \frac{1}{T(T-1)} e_{-1}\Big)\int_{-T}^0 \big(\eta_1(x) - \eta_1(-T)\big) \phi_n(x + T) dx \\
&\quad \ + \sum_{i=0}^n \frac{n+1-i}{n+1} \xi_ie_i + \xi_{-1}e_{-1},\,\text{for some }\eta_1\in K,\,n\in\N,\,\xi\in E_{n+2}\bigg\}.
\end{align*}

\vspace{1mm}

\noindent\textbf{Digression.} \emph{$\mathcal K$ is a relatively compact subset of $C([-T,0])$.} Since $K$ is compact, it is enough to prove that $\tilde K$ is relatively compact. To this end, define
\begin{align*}
K_1 \ := \ \bigg\{\eta\in C([-T,0])\colon\eta \ &= \ T_n\eta_1 + \Big(T_n\gamma + \frac{1}{T(T-1)} e_{-1}\Big)\int_{-T}^0 \big(\eta_1(x) - \eta_1(-T)\big) \phi_n(x + T) dx \\
&\quad \ \text{for some }\eta_1\in K,\,n\in\N\bigg\}, \\
K_2 \ := \ \bigg\{\eta\in C([-T,0])\colon\eta \ &= \ \sum_{i=-1}^n \xi_ie_i,\,\text{for some }\,n\in\N,\,\xi\in E_{n+2}\bigg\}.
\end{align*}
Then $\tilde K\subset K_1+K_2$, where $K_1+K_2$ denotes the sum of the sets $K_1$ and $K_2$, i.e., $K_1+K_2=\{\eta\in C([-T,0])\colon\eta=\eta_1+\eta_2,\,\text{for some }\eta_1\in K_1,\,\eta_2\in K_2\}$. In order to prove that $\tilde K$ is relatively compact, it is enough to show that both $K_1$ and $K_2$ are relatively compact sets.

Firstly, let us prove that $K_1$ is relatively compact. Take a sequence $(\eta_\ell)_{\ell\in\N}$ in $K_1$. Our aim is to prove that $(\eta_\ell)_{\ell\in\N}$ admits a convergent subsequence. We begin noting that, for every $\ell\in\N$, there exist $\eta_{1,\ell}\in C([-T,0])$ and $n_\ell\in\N$ such that
\[
\eta_\ell \ = \ T_{n_\ell}\eta_{1,\ell} + \Big(T_{n_\ell}\gamma + \frac{1}{T(T-1)} e_{-1}\Big)\int_{-T}^0 \big(\eta_{1,\ell}(x) - \eta_{1,\ell}(-T)\big) \phi_{n_\ell}(x + T) dx.
\]
Let us suppose that $(n_\ell)_{\ell\in\N}$ admits a subsequence diverging 
 to infinity (the other cases can be treated even simpler),
still denoted by $(n_\ell)_{\ell\in\N}$. Then $T_{n_\ell}\gamma\rightarrow\gamma$ in $C([-T,0])$. Since $(\eta_{1,\ell})_{\ell\in\N}\subset K$ and $K$ is compact, there exists a subsequence, still denoted by $(\eta_{1,\ell})_{\ell\in\N}$, which converges to some $\eta_{1,\infty}\in K$. Then, $T_{n_\ell}\eta_{1,\ell}\rightarrow\eta_{1,\infty}$ as $\ell\rightarrow\infty$. Indeed
\[
\|T_{n_\ell}\eta_{1,\ell}-\eta_{1,\infty}\|_\infty \ \leq \ \|T_{n_\ell}\eta_{1,\ell}-T_{n_\ell}\eta_{1,\infty}\|_\infty + \|T_{n_\ell}\eta_{1,\infty}-\eta_{1,\infty}\|_\infty.
\]
Then, the claim follows since $T_{n_\ell}\eta_{1,\infty}\rightarrow\eta_{1,\infty}$ in $C([-T,0])$ and
\[
\|T_{n_\ell}\eta_{1,\ell}-T_{n_\ell}\eta_{1,\infty}\|_\infty \ \overset{\text{by \eqref{E:UniformBoundT_n}}}{\leq} \ M \|\eta_{1,\ell}-\eta_{1,\infty}\|_\infty \ \overset{\ell\rightarrow\infty}{\longrightarrow} \ 0.
\]
Proceeding in a similar way, we see that
\begin{align*}
\int_{-T}^0 \big(\eta_{1,\ell}(x) - \eta_{1,\ell}(-T)\big) \phi_{n_\ell}(x + T) dx \ &= \ \int_{-T}^0 \eta_{1,\ell}(x)\phi_{n_\ell}(x + T) dx - \eta_{1,\ell}(-T) \\
&\overset{\ell\rightarrow\infty}{\longrightarrow} \ \eta_{1,\infty}(-T) - \eta_{1,\infty}(-T) \ = \ 0.
\end{align*}
In conclusion, we get $\eta_\ell\rightarrow\eta_{1,\infty}$, from which the claim follows.

Let us now prove that $K_2$ is relatively compact. Let $(\eta_\ell)_{\ell\in\N}$ be a sequence in $K_2$ and let us prove that $(\eta_\ell)_{\ell\in\N}$ admits a convergent subsequence in $C([-T,0])$. We first notice that, for every $\ell\in\N$, there exists $n_\ell\in\N$ and $\xi_\ell=(\xi_{-1,\ell},\ldots,\xi_{n_\ell,\ell})\in E_{n_\ell+2}$ such that
\[
\eta_\ell \ = \ \sum_{i=-1}^{n_\ell} \xi_{i,\ell}e_i.
\]
As we already did in the proof for $K_1$, we suppose that the sequence $(n_\ell)_{\ell\in\N}$ diverges to $\infty$. Notice that, for every $i\in\{-1,0,1,2,\ldots\}$, there exists a subsequence of $(\xi_{i,\ell})_\ell$ which converges to some $\xi_{i,\infty}$ satisfying $|\xi_{i,\infty}|\leq2^{-(i+1)}$. By a diagonalisation argument we construct a subsequence of $(\eta_\ell)_{\ell\in\N}$, still denoted by $(\eta_\ell)_{\ell\in\N}$, such that for every $i$ the sequence $(\xi_{i,\ell})_{\ell\in\N}$ converges to $\xi_{i,\infty}$. As a consequence, $\eta_\ell$ converges to $\eta_\infty=\sum_{i=-1}^\infty \xi_{i,\infty}e_i$ as $\ell\rightarrow\infty$. This proves the claim.

\vspace{1mm}

\noindent\textbf{Step III (Continued).} Since $\Kc$ is a relatively compact subset of $C([-T,0])$, property \eqref{UniformContinuitybarH_n} follows from the fact that $H$ is continuous on $C([-T,0])$, and consequently uniformly continuous on $\mathcal K$.

To alleviate the presentation, we suppose, without loss of generality, that $H_n$ has the following form (with the same functions $\varphi_i$ as in the expression of $b_n,\sigma_n,F_n$)
\[
H_n(\eta) \ = \ \bar H_n\bigg(\int_{[-T,0]}\varphi_1(x+T)d^-\eta(x),\ldots,\int_{[-T,0]}\varphi_{N_n}(x+T)d^-\eta(x)\bigg).
\]
So that $\bar H_n\colon\R^{N_n}\rightarrow\R$. Then, property \eqref{UniformContinuitybarH_n} can be written as follows: for any compact subset $K\subset C([-T,0])$ there exists a continuity modulus $\rho_K$, independent of $n\in\N$, such that
\begin{align}
\label{UniformContinuitybarH_n2}
\bigg|\bar{H}_n\bigg(\int_{[-T,0]}\varphi_1(x+T)d^-\eta_1(x) + \xi_1,\ldots\bigg) - \bar{H}_n\bigg(\int_{[-T,0]}\varphi_1(x+T)d^-\eta_2(x) + \xi_1,\ldots\bigg)\bigg| \notag \\
\leq \ m_K(\|\eta_1-\eta_2\|_\infty),
\end{align}
for all $\eta_1,\eta_2\in K$, $n\in\N$, $\xi=(\xi_1,\ldots,\xi_{N_n})\in E_{N_n}$, where we recall that $E_{N_n}=\{\xi=(\xi_1,\ldots,\xi_{N_n})\in\R^{N_n}\colon|\xi_i|\leq 2^{1-i},\,i=1,\ldots,N_n\}$. 

Now, for any $n$ consider the function $\rho_n\in C^\infty(\R^{N_n})$ given by
\begin{equation}\label{rho}
\rho_n(\xi) \ = \ c \prod_{i=1}^{N_n} \exp\bigg(\frac{1}{\xi_i^2 - 2^{2(i-1)}}\bigg) 1_{\{|\xi_i|<2^{i-1}\}}, \qquad \forall\,\xi=(\xi_1,\ldots,\xi_{N_n})\in\R^{N_n},
\end{equation}
with $c>0$ such that $\int_{\R^{N_n}} \rho_n(\xi) d\xi = 1$. Set $\rho_{n,k}(\xi) := k^{N_n}\rho_n(k\,\xi)$, $\forall\,\xi\in\R^{N_n}$, $k\in\N$. Let us now define, for any $n,k\in\N$,
\[
\bar H_{n,k}(x) \ = \ \int_{\R^{N_n}} \rho_{n,k}(\xi) \bar H_n(x - \xi) d\xi \ = \ \int_{E_{N_n}} \rho_{n,k}(\xi) \bar H_n(x - \xi) d\xi,
\]
for all $(t,x,y,z)\in[0,T]\times\R^d\times\R\times\R^d$. Notice that, for any $n\in\N$, the sequence $(\bar{H}_{n,k}(\cdot))_{k\in\N}$ is equicontinuous on compact subsets of $\R^{N_n}$, satisfies a polynomial growth condition (uniform in both $n$ and $k$), converges pointwise to $\bar{H}_n(\cdot)$, and satisfies item (iv) of Theorem \ref{T:ExistenceStrict}. Then, we define
\[
{H}_{n,k}(\eta) \ = \ \bar{H}_{n,k}\bigg(\int_{[-T,0]}\varphi_1(x+T)d^-\eta(x),\ldots,\int_{[-T,0]}\varphi_{N_n}(x+T)d^-\eta(x)\bigg),
\]
for all $\eta\in C([-T,0])$ and $n,k\in\N$. Notice that the functions $H_{n,k}$, $n,k\in\N$, are equicontinuous on compact subsets of $C([-T,0])$. Indeed, let $K$ be a compact subset of $C([-T,0])$ and $\eta_1,\eta_2\in K$, then (using property \eqref{UniformContinuitybarH_n2} and the fact that $\int_{E_{N_n}} \rho_{n,k}(\xi) d\xi = 1$)
\begin{align*}
&|H_{n,k}(\eta_1) - H_{n,k}(\eta_2)| \\
&= \ \bigg|\bar H_{n,k}\bigg(\int_{[-T,0]}\varphi_1(x+T)d^-\eta_1(x),\ldots\bigg) - \bar H_{n,k}\bigg(\int_{[-T,0]}\varphi_1(x+T)d^-\eta_2(x),\ldots\bigg)\bigg| \\
&\leq \ \int_{K_{N_n}} \rho_{n,k}(\xi) \bigg|\bar H_n\bigg(\int_{[-T,0]}\varphi_1(x+T)d^-\eta_1(x) + \xi_1,\ldots\bigg) \\
&\quad \ - \bar H_n\bigg(\int_{[-T,0]}\varphi_1(x+T)d^-\eta_2(x) + \xi_1,\ldots\bigg)\bigg| d\xi \ \leq \ m_K(\|\eta_1-\eta_2\|_\infty).
\end{align*}
This proves the equicontinuity on compact sets of $H_{n,k}$, $n,k\in\N$. Set $G:=H$, $G_n:=H_n$, and $G_{n,k}:=H_{n,k}$, for all $n,k\in\N$. Then, a direct application of Lemma \ref{L:StabilityApp} yields the existence of a subsequence $(H_{n,k_n})_{n\in\N}$ which converges pointwise to $H$. For simplicity of notation, we denote $(H_{n,k_n})_{n\in\N}$ simply by $(H_n)_{n\in\N}$.

\vspace{1mm}

\noindent\textbf{Step IV.} \emph{Conclusion.} Let us consider, for any $n\in\N$ and $(t,\eta)\in[0,T]\times C([-T,0])$, the following forward-backward system of stochastic differential equations:
\begin{equation}
\label{FBSDE}
\begin{cases}
X_s^{n,t,\eta} \ = \ \eta(0\wedge(s-t)) + \int_t^{t\vee s} b_n(r,\X_r^{n,t,\eta}) dr + \int_t^{t\vee s} \sigma_n(r,\X_r^{n,t,\eta}) dW_r, &s\in[t-T,T], \\
Y_s^{n,t,\eta} \ = \ H_n(\X_T^{n,t,\eta}) + \int_s^T F_n(r,\X_r^{n,t,\eta},Y_r^{n,t,\eta},Z_r^{n,t,\eta}) dr - \int_s^T Z_r^{n,t,\eta} dW_r, &s\in[t,T].
\end{cases}
\end{equation}
Under the assumptions on $b_n$ and $\sigma_n$, it follows from Proposition \ref{P:SDE} that there exists a unique continuous process $X^{n,t,\eta}$ strong solution to the forward equation in \eqref{FBSDE}. Moreover, from Theorem 4.1 in \cite{parpen90} it follows that, under the assumptions on $F_n$ and $H_n$, there exists a unique solution $(Y^{n,t,\eta},Z^{n,t,\eta})\in\S^2(t,T)\times\H^2(t,T)$ to the backward equation in \eqref{FBSDE}.

Then, it follows from Theorem \ref{T:ExistenceStrict} that, for any $n\in\N$, the function
\[
\Uc_n(t,\eta) \ = \ Y_t^{n,t,\eta} \ = \ \E\bigg[\int_t^T F_n(s,\mathbb X_s^{n,t,\eta},Y_s^{n,t,\eta},Z_s^{n,t,\eta})ds + {H}_n(\mathbb X_T^{n,t,\eta})\bigg],
\]
$\forall\,(t,\eta)\in[0,T]\times C([-T,0])$, is a strict solution to equation \eqref{KolmEq} with coefficients ${H}_n$, $F_n$, $b_n$, and $\sigma_n$. From estimates \eqref{EstimateSupXn} and \eqref{EstimateBSDE2} together with the polynomial growth condition of $F_n,H_n$ (uniform in $n$), we see that $\Uc_n$ satisfies a polynomial growth condition uniform in $n$. 

We can now apply Lemma \ref{L:Stability} to the sequence $(\Uc_n,H_n,F_n,b_n,\sigma_n)_{n\in\N}$, from which we deduce: first, the convergence of the sequence $(\Uc_n)_{n\in\N}$ to the map $\Uc$ given by \eqref{Feynman-Kac2}; secondly, that $\Uc$ is a strong-viscosity solution to equation \eqref{KolmEq}. This concludes the proof.
\ep

\begin{Remark}
{\rm
(i)
Here we notice that Theorem \ref{T:ExistSV} applies when $b$, $\sigma$, $F$ have a Markovian structure. More precisely, suppose that there exist $\bar b$, $\bar\sigma$, $\bar f$ satisfying Assumption {\bf (A0)}, with $\bar f$ continuous, such that
\[
b(t,\eta) \ = \ \bar b(\eta(0)), \qquad \sigma(t,\eta) \ = \ \bar\sigma(\eta(0)), \qquad F(t,\eta,y,z) \ = \ \bar f(t,\eta(0),y,z),
\]
for all $(t,\eta,y,z)\in[0,T]\times C([-T,0])\times\R\times\R$. Recalling from the integration by parts formula (2.4) in \cite{cosso_russo15a} that $\eta(0)=\int_{[-t,0]}1\,d^-\eta(x)$, we see that $b$, $\sigma$, $F$ have automatically a cylindrical form, as a matter of fact
\[
b(t,\eta) \ = \ \bar b\bigg(\int_{[-t,0]}1\,d^-\eta(x)\bigg)
\]
and similarly for $\sigma$ and $F$. Therefore, taking $b$, $\sigma$, $F$ of this form, and $H$ continuous (as in the statement of Theorem \ref{T:ExistSV}), we deduce from Theorem \ref{T:ExistSV} that the map $\Uc$ given by \eqref{Feynman-Kac2} is a 
strong-viscosity solution to equation \eqref{KolmEq}.

\vspace{1mm}
\noindent(ii) The result of Theorem \ref{T:ExistSV} can be improved as follows. Items (ii) and (iii) in Theorem \ref{T:ExistSV} can be replaced by the following weaker assumption: \emph{for every compact subset $K\subset C([-T,0])$, there exists a continuity modulus $m_K$, independent of $n\in\N$, such that
\begin{align*}
\bigg|\bar{F}_n\!\bigg(\!t,\!\int_{[-t,0]}\!\!\!\varphi_1(x+t)d^-\eta_1(x) + \xi_1,\ldots,y,z\!\bigg)\! - \bar{F}_n\!\bigg(\!t,\!\int_{[-t,0]}\!\!\!\varphi_1(x+t)d^-\eta_2(x) + \xi_1,\ldots,y,z\!\bigg)\bigg| \\
+ \bigg|\bar{b}_n\bigg(\int_{[-t,0]}\varphi_1(x+t)d^-\eta_1(x) + \xi_1,\ldots\bigg) - \bar{b}_n\bigg(\int_{[-t,0]}\varphi_1(x+t)d^-\eta_2(x) + \xi_1,\ldots\bigg)\bigg| \\
+ \bigg|\bar{\sigma}_n\bigg(\int_{[-t,0]}\varphi_1(x+t)d^-\eta_1(x) + \xi_1,\ldots\bigg) - \bar{\sigma}_n\bigg(\int_{[-t,0]}\varphi_1(x+t)d^-\eta_2(x) + \xi_1,\ldots\bigg)\bigg| \\
\leq \ m_K(\|\eta_1 - \eta_2\|_\infty),
\end{align*}
for all $n\in\N$, $\eta_1,\eta_2\in K$, $y,z\in\R$, $t\in[0,T]$, $\xi\in E_{N_n}$, where $E_{N_n}=\{\xi=(\xi_1,\ldots,\xi_{N_n})\in\R^{N_n}\colon|\xi_i|\leq 2^{1-i},\,i=1,\ldots,N_n\}$.}

In this case, we perform a smoothing of $(\bar b_n,\bar\sigma_n,\bar F_n)$ by means of convolutions as we did for $\bar H_n$ in Step III of the proof of Theorem \ref{T:ExistSV}, in order to end up with a sequence of regular coefficients satisfying items (ii) and (iii) in Theorem \ref{T:ExistSV}. Then, we conclude the proof proceeding as in Step IV of the proof of Theorem \ref{T:ExistSV}.

\noindent{(iii)}  The particular case $b\equiv0$, $\sigma\equiv1$, and $F\equiv0$ was addressed in Theorem 3.4 of \cite{cosso_russo14}. 
Concerning the case with general coefficients $b$, $\sigma$, $F$, we refer to Theorem \ref{T:ExistSV_3} below.
\ep
}
\end{Remark}

We also state the following existence result, which holds under slightly different assumptions than Theorem \ref{T:ExistSV}.

\begin{Theorem}
\label{T:ExistSV_2}
Let Assumption {\bf (A1)} hold. Suppose also that ${H}$ is continuous and $F$ does not depend on $(y,z)$. In addition, suppose that there exists a nondecreasing sequence $(N_n)_{n\in\N}\subset\N\backslash\{0\}$ such that, for all $n\in\N$ and $(t,\eta)\in[0,T]\times C([-T,0])$,
\begin{align*}
b_n(t,\eta) \ &= \ \bar b_n\bigg(t,\int_{[-t,0]}\varphi_1(x+t)d^-\eta(x),\ldots,\int_{[-t,0]}\varphi_{N_n}(x+t)d^-\eta(x)\bigg), \\
\sigma_n(t,\eta) \ &= \ \bar\sigma_n\bigg(t,\int_{[-t,0]}\varphi_1(x+t)d^-\eta(x),\ldots,\int_{[-t,0]}\varphi_{N_n}(x+t)d^-\eta(x)\bigg), \\
F_n(t,\eta) \ &= \ \bar F_n\bigg(t,\int_{[-t,0]}\varphi_1(x+t)d^-\eta(x),\ldots,\int_{[-t,0]}\varphi_{N_n}(x+t)d^-\eta(x)\bigg).
\end{align*}
We suppose that items \textup{(i)}, \textup{(iii)}, \textup{(iv)}, \textup{(v)} of Theorem \ref{T:ExistSV} hold, while item \textup{(ii)} is replaced by the following:
\begin{itemize}
\item[\textup{(ii)'}] For every $n\in\N$, $\bar b_n,\bar\sigma_n,\bar F_n$ satisfy:
\begin{itemize}
\item[\textup{(ii)'-(a)}] $\bar b_n$ and $\bar\sigma_n$ are continuous functions, with first and second spatial derivatives continuous and satisfying a polynomial growth condition.
\item[\textup{(ii)'-(b)}] $\bar F_n$ is continuous and, for all $t\in[0,T]$, the function $\bar F_n(t,\cdot)$ belongs to $C^2(\R^{N_n})$ and its second order spatial derivatives satisfy a polynomial growth condition uniformly in $t$.
\end{itemize}

\end{itemize}
Then, the map $\Uc$ given by \eqref{Feynman-Kac2} is a (path-dependent) strong-viscosity solution to equation \eqref{KolmEq}.
\end{Theorem}
\begin{Remark}\label{R:F(t,x)}
{\rm
Notice that the requirement that $F$ does not depend on $(y,z)$ imposed in Theorem \ref{T:ExistSV_2} is used only at a specific point in the proof of Theorem \ref{T:ExistSV_2}, namely when it is used Theorem 3.5 in \cite{cosso_russo15a}, which relies on Theorem 6.1, Chapter 5, in \cite{friedman75vol1}, that is a regularity result for linear (rather than semilinear) parabolic partial differential equations. However, if one would have at disposal a more general regularity result than Theorem 6.1, Chapter 5, in \cite{friedman75vol1}, which applies to semilinear parabolic partial differential equations, then we would be able to extend Theorem \ref{T:ExistSV_2} (as well as Theorem 3.5 in \cite{cosso_russo15a}) to the case where $F$ also depends on $(y,z)$.
\ep}
\end{Remark}
\textbf{Proof.}
The proof can be done proceeding along the same lines as in the proof of Theorem \ref{T:ExistSV}, the only difference being that in Step IV we rely on Theorem 3.5 in \cite{cosso_russo15a} rather than on Theorem \ref{T:ExistenceStrict}
of this paper.
\ep

\vspace{3mm}

We finally state the following existence result, which relies on the previous Theorem \ref{T:ExistSV_2}.
\begin{Theorem}
\label{T:ExistSV_3}
Let Assumption {\bf (A1)} hold. We also suppose the following.
\begin{enumerate}[\upshape (a)]
\item $b$, $\sigma$, $F$, $H$ are continuous;
\item $F$ does not depend on $(y,z)$;
\item for every $t\in[0,T]$, the map $F(t,\cdot)$ is continuous uniformly with respect to $t\in[0,T]$;
\item $b$, $\sigma$, $F$ satisfy the following property:
\begin{equation}\label{b_sigma_F_non_anticipative}
b(t,\eta) \ = \ b(t,\gamma), \qquad\quad \sigma(t,\eta) \ = \ \sigma(t,\gamma), \qquad\quad F(t,\eta) \ = \ F(t,\gamma),
\end{equation}
for every $t\in[0,T]$, $\eta,\gamma\in C([-T,0])$, with $\eta(x)=\gamma(x)$ for any $x\in[-t,0]$.
\end{enumerate}
Then, the map $\Uc$ given by \eqref{Feynman-Kac2} is a (path-dependent) strong-viscosity solution to equation \eqref{KolmEq}.
\end{Theorem}
\begin{Remark}
{\rm
Notice that the requirement that $F$ does not depend on $(y,z)$ is needed only because in the proof of Theorem \ref{T:ExistSV_3} we use Theorem \ref{T:ExistSV_2} above, for which in turn we refer to Remark \ref{R:F(t,x)}.
\ep}
\end{Remark}
\textbf{Proof.} By Theorem \ref{T:ExistSV_2}, it is enough to prove that there exists a nondecreasing sequence $(N_n)_{n\in\N}\subset\N\backslash\{0\}$ such that, for all $n\in\N$ and $(t,\eta)\in[0,T]\times C([-T,0])$,
\begin{align*}
b_n(t,\eta) \ &= \ \bar b_n\bigg(t,\int_{[-t,0]}\varphi_1(x+t)d^-\eta(x),\ldots,\int_{[-t,0]}\varphi_{N_n}(x+t)d^-\eta(x)\bigg), \\
\sigma_n(t,\eta) \ &= \ \bar\sigma_n\bigg(t,\int_{[-t,0]}\varphi_1(x+t)d^-\eta(x),\ldots,\int_{[-t,0]}\varphi_{N_n}(x+t)d^-\eta(x)\bigg), \\
F_n(t,\eta) \ &= \ \bar F_n\bigg(t,\int_{[-t,0]}\varphi_1(x+t)d^-\eta(x),\ldots,\int_{[-t,0]}\varphi_{N_n}(x+t)d^-\eta(x)\bigg),
\end{align*}
where $(b_n,\sigma_n,F_n)_n$ and $(\bar b_n,\bar\sigma_n,\bar F_n)_n$ satisfy items \textup{(i)}, \textup{(ii)'}, \textup{(iii)}, \textup{(iv), \textup{(v)}} of Theorem \ref{T:ExistSV_2}. Then, the fact that the map $\Uc$ given by \eqref{Feynman-Kac2} is a (path-dependent) strong-viscosity solution to equation \eqref{KolmEq} follows directly from Theorem \ref{T:ExistSV_2}. We divide the proof of the construction of the sequences $(b_n,\sigma_n,F_n)_n$ and $(\bar b_n,\bar\sigma_n,\bar F_n)_n$ into five steps.

\vspace{1mm}

\noindent\textbf{Step I.} \emph{Polygonal approximation of $\eta\in C([-T,0])$.} For every $n\in\N$, we consider the $n$-th dyadic subdivision of $[0,T]$, namely
\[
0 \ = \ t_0^n \ < \ t_1^n \ < \ \ldots \ < \ t_{N_n}^n \ = \ T, \qquad\quad \text{where }N_n \ := \ 2^n \text{ and }t_j^n \ := \ \frac{j}{2^n}T,\;\forall\,j=0,\ldots,2^n. 
\]
For every fixed $n\in\N$, $t\in[0,T]$, $\eta\in C([-T,0])$, we consider the $n$-th polygonal approximation $\tilde\eta_n^t\in C([-T,0])$ of path $\eta$ at time $t$, defined as 
\begin{itemize}
\item $\tilde\eta_n^t$ is constant on $[-T,-t]$ and equal to $\eta(-t)$;
\item $\tilde\eta_n^t$ is linear on $[(t_{j-1}^n\wedge t)-t,(t_j^n\wedge t)-t]$, for every $j=1,\ldots,N_n$;
\item $\tilde\eta_n^t((t_j^n\wedge t)-t)=\eta((t_j^n\wedge t)-t)$, for every $j=0,\ldots,N_n$.
\end{itemize}
Notice that $-t=(t_0^n\wedge t)-t<(t_1^n\wedge t)-t<\cdots<(t_{N_n}^n\wedge t)-t=0$, so that the finite sequence $((t_j^n\wedge t)-t)_j$ is a subdivision of $[-t,0]$. Then, on the interval $[-t,0]$ we see that the continuous function $\tilde\eta_n^t$ is given by
\begin{align*}
\tilde\eta_n^t(x) \ &= \ \frac{\eta((t_j^n\wedge t)-t) - \eta((t_{j-1}^n\wedge t)-t)}{(t_j^n\wedge t)-(t_{j-1}^n\wedge t)} x \\
&\quad \ + \frac{((t_j^n\wedge t)-t)\,\eta((t_{j-1}^n\wedge t)-t) - ((t_{j-1}^n\wedge t)-t)\,\eta((t_j^n\wedge t)-t)}{(t_j^n\wedge t)-(t_{j-1}^n\wedge t)},
\end{align*}
for every $x\in[(t_{j-1}^n\wedge t)-t,(t_j^n\wedge t)-t]$, whenever $t_{j-1}^n<t$, so that $(t_{j-1}^n\wedge t)-t<(t_j^n\wedge t)-t$. Notice that, for any $\eta,\gamma\in C([-T,0])$ with $\eta(x)=\gamma(x)$ for any $x\in[-t,0]$, we have $\tilde\eta_n^t\equiv\tilde\gamma_n^t$. Moreover, for any $\eta\in C([-T,0])$, let $\eta^t\in C([-T,0])$ denote the continuous path satisfying: $\eta^t\equiv\eta(-t)$ on $[-T,-t]$ and $\eta^t\equiv\eta$ on $[-t,0]$. Then, from the uniform continuity of $\eta$, and hence of $\eta^t$, we deduce that $\|\tilde\eta_n^t-\eta^t\|_\infty\rightarrow0$ as $n\rightarrow+\infty$.

Define the maps $\tilde b_n\colon[0,T]\times C([-T,0])\rightarrow\R$, $\tilde\sigma_n\colon[0,T]\times C([-T,0])\rightarrow\R$, $\tilde F_n\colon[0,T]\times C([-T,0])\rightarrow\R$ as
\begin{equation}\label{b_sigma_F_tilde}
\tilde b_n(t,\eta) \ := \ b(t,\tilde\eta_n^t), \qquad\quad \tilde\sigma_n(t,\eta) \ := \ \sigma(t,\tilde\eta_n^t), \qquad\quad \tilde F_n(t,\eta) \ := \ F(t,\tilde\eta_n^t),
\end{equation}
for every $(t,\eta)\in[0,T]\times C([-T,0])$. Since $\|\tilde\eta_n^t-\eta^t\|_\infty\rightarrow0$ and $b$, $\sigma$, $F$ satisfy property \eqref{b_sigma_F_non_anticipative}, we deduce that the sequence $(\tilde b_n,\tilde\sigma_n,\tilde F_n)_n$ converges pointwise to $(b,\sigma,F)$ as $n\rightarrow+\infty$.

Now, given $\mathbf y=(y_0,\ldots,y_{N_n})\in\R^{N_n+1}$ and $t\in[0,T]$, we define the polygonal $\tilde\eta_{n,\mathbf y}^t\in C([-T,0])$ associated with $\mathbf y$, defined as follows:
\begin{itemize}
\item $\tilde\eta_{n,\mathbf y}^t$ is constant on $[-T,-t]$ and equal to $y_0$;
\item $\tilde\eta_{n,\mathbf y}^t$ is linear on $[(t_{j-1}^n\wedge t)-t,(t_j^n\wedge t)-t]$, for every $j=1,\ldots,N_n$;
\item $\tilde\eta_n^t((t_j^n\wedge t)-t)=y_j$, for every $j=0,\ldots,N_n$.
\end{itemize}
More precisely, on the interval $[-t,0]$ the continuous function $\tilde\eta_{n,\mathbf y}^t$ is given by
\[
\tilde\eta_{n,\mathbf y}^t(x) \ := \ \frac{y_j - y_{j-1}}{(t_j^n\wedge t)-(t_{j-1}^n\wedge t)} x + \frac{((t_j^n\wedge t)-t)\,y_{j-1} - ((t_{j-1}^n\wedge t)-t)\,y_j}{(t_j^n\wedge t)-(t_{j-1}^n\wedge t)},
\]
for every $x\in[(t_{j-1}^n\wedge t)-t,(t_j^n\wedge t)-t]$, whenever $t_{j-1}^n<t$, so that $(t_{j-1}^n\wedge t)-t<(t_j^n\wedge t)-t$.

Define the maps $\hat b_n\colon[0,T]\times\R^{N_n+1}\rightarrow\R$, $\hat\sigma_n\colon[0,T]\times\R^{N_n+1}\rightarrow\R$, $\hat F_n\colon[0,T]\times\R^{N_n+1}\rightarrow\R$ as
\[
\hat b_n(t,y_0,\ldots,y_{N_n}) := b(t,\tilde\eta_{n,\mathbf y}^t), \quad \hat\sigma_n(t,y_0,\ldots,y_{N_n}) := \sigma(t,\tilde\eta_{n,\mathbf y}^t), \quad \hat F_n(t,y_0,\ldots,y_{N_n}) := F(t,\tilde\eta_{n,\mathbf y}^t),
\]
for every $t\in[0,T]$, $\mathbf y=(y_0,\ldots,y_{N_n})\in\R^{N_n+1}$. Notice that there exists the following relation between the maps $\tilde b_n,\tilde\sigma_n,\tilde F_n$ and $\hat b_n,\hat\sigma_n,\hat F_n$: 
\begin{align*}
\tilde b_n(t,\eta) \ &= \ \hat b_n(t,\eta((t_0^n\wedge t) - t),\ldots,\eta((t_{N_n}^n\wedge t) - t)), \\
\tilde\sigma_n(t,\eta) \ &= \ \hat\sigma_n(t,\eta((t_0^n\wedge t) - t),\ldots,\eta((t_{N_n}^n\wedge t) - t)), \\
\tilde F_n(t,\eta) \ &= \ \hat F_n(t,\eta((t_0^n\wedge t) - t),\ldots,\eta((t_{N_n}^n\wedge t) - t)),
\end{align*}
for every $(t,\eta)\in[0,T]\times C([-T,0])$. Recalling from the integration by parts formula (2.4) in \cite{cosso_russo15a} that $\eta((t_j^n\wedge t) - t)=\int_{[-t,0]}1_{[0,t_j^n]}(x+t)d^-\eta(x)$, we can rewrite the above equalities as follows:
\begin{align*}
\tilde b_n(t,\eta) \ &= \ \hat b_n\bigg(t,\int_{[-t,0]}1_{[0,t_0^n]}(x+t)d^-\eta(x),\ldots,\int_{[-t,0]}1_{[0,t_{N_n}^n]}(x+t)d^-\eta(x)\bigg), \\
\tilde\sigma_n(t,\eta) \ &= \ \hat\sigma_n\bigg(t,\int_{[-t,0]}1_{[0,t_0^n]}(x+t)d^-\eta(x),\ldots,\int_{[-t,0]}1_{[0,t_{N_n}^n]}(x+t)d^-\eta(x)\bigg), \\
\tilde F_n(t,\eta) \ &= \ \hat F_n\bigg(t,\int_{[-t,0]}1_{[0,t_0^n]}(x+t)d^-\eta(x),\ldots,\int_{[-t,0]}1_{[0,t_{N_n}^n]}(x+t)d^-\eta(x)\bigg),
\end{align*}
for every $(t,\eta)\in[0,T]\times C([-T,0])$.

\vspace{1mm}

\noindent\textbf{Step II.} \emph{The maps $\hat b_n$, $\hat\sigma_n$, $\hat F_n$ satisfy item (i) of Theorem \ref{T:ExistSV_2}.} We begin noting that, given $\mathbf y=(y_0,\ldots,y_{N_n})$, $\mathbf y'=(y_0',\ldots,y_{N_n}')\in\R^{N_n+1}$, we have
\begin{equation}\label{eta_y_Lipschitz}
\|\tilde\eta_{n,\mathbf y}^t\|_\infty \ \leq \ \max_j |y_j| \ \leq \ |\mathbf y|, \qquad\qquad \|\tilde\eta_{n,\mathbf y}^t - \tilde\eta_{n,\mathbf y'}^t\|_\infty \ \leq \ \max_j |y_j - y_j'| \ \leq \ |\mathbf y - \mathbf y'|,
\end{equation}
where $|\mathbf y|=(y_0^2+\cdots+y_{N_n}^2)^{1/2}$ denotes the Euclidean norm of $\mathbf y$. Then, denoting by $C$ and $m$ the constants appearing in assumption {\bf (A1)}, it follows that $\hat b_n$, $\hat\sigma_n$, $\hat F_n$ satisfy the following conditions (with the same constants $C$ and $m$):
\begin{align*}
|\hat b_n(t,\mathbf y)-\hat b_n(t,\mathbf y')| + |\hat\sigma_n(t,\mathbf y)-\hat\sigma_n(t,\mathbf y')| \ &\leq \ C|\mathbf y - \mathbf y'|, \\
|\hat b_n(t,0)| + |\hat\sigma_n(t,0)| \ &\leq \ C, \\
|\hat F_n(t,\mathbf y)| \ &\leq \ C\big(1 + |\mathbf y|^m\big),
\end{align*}
for all $t\in[0,T]$, $\mathbf y,\mathbf y'\in\R^{N_n+1}$. Now, fix $n\in\N$, $\mathbf y=(y_0,\ldots,y_{N_n})\in\R^{N_n+1}$, $s,t\in[0,T]$, with $s\leq t$. Notice that
\[
\tilde\eta_{n,\mathbf y}^s(x) \ = \
\begin{cases}
y_0, & x\in[-T,-T+(t-s)], \\
\tilde\eta_{n,\mathbf y}^t(x - (t - s)), \qquad & x\in[-T+(t-s),0].
\end{cases}
\]
This proves the continuity of the map $t\mapsto \tilde\eta_{n,\mathbf y}^t$, from $[0,T]$ to $C([-T,0])$. From this latter property, together with \eqref{eta_y_Lipschitz} and the continuity of $b$, $\sigma$, $F$, we deduce that the maps $\hat b_n$, $\hat\sigma_n$, $\hat F_n$ are continuous in both arguments. In conclusion, $\hat b_n$, $\hat\sigma_n$, $\hat F_n$ satisfy item (i) of Theorem \ref{T:ExistSV_2}.


\vspace{1mm}

\noindent\textbf{Step III.} \emph{The maps $\tilde b_n$, $\tilde\sigma_n$, $\tilde F_n$ satisfy item (iii) of Theorem \ref{T:ExistSV_2}.} Fix a compact set $K\subset C([-T,0])$. Our aim is to prove that the following subset of $C([-T,0])$ is relatively compact:
\begin{equation}\label{compact_set_K}
\mathcal K \ := \ \big\{\gamma\in C([-T,0])\colon\gamma\equiv\tilde\eta_n^t,\text{ for some }n\in\N,\,t\in[0,T],\,\eta\in K\big\}.
\end{equation}
To this end, take a sequence $(\gamma_k)_{k\in\N}\subset\mathcal K$. Then, for every $k\in\N$, there exist $n_k\in\N$, $t_k\in[0,T]$, $\eta_k\in C([-T,0])$ such that $\gamma_k\equiv\tilde\eta_{k,n_k}^{t_k}$. Since $K$ is compact, there exists $\eta\in K$ such that, up to a subsequence, $(\eta_k)_{k\in\N}$ converges to $\eta$ in $C([-T,0])$. Similarly, there exists $t\in[0,T]$ such that, up to a subsequence, $(t_k)_{k\in\N}$ converges to $t$. Finally, concerning the sequence $(n_k)_{k\in\N}$ we distinguish two cases: $(n_k)_{k\in\N}$ goes, up to a subsequence, to $+\infty$; $(n_k)_{k\in\N}$ is identically equal to some $n_0\in\N$, up to a subsequence. This latter case is easier to be treated, therefore we do not report the proof for this case, and in the sequel we suppose that $n_k\rightarrow+\infty$. Let us then prove that the sequence $(\tilde\eta_{k,n_k}^{t_k})_{k\in\N}$ converges to $\eta^t$ in $C([-T,0])$, where we recall from Step I that the continuous path $\eta^t\in C([-T,0])$ is equal to the constant $\eta(-t)$ on $[-T,-t]$ and coincides with $\eta$ on $[-t,0]$. From the triangular inequality, we have
\[
\|\tilde\eta_{k,n_k}^{t_k} - \eta^t\|_\infty \ \leq \ \|\tilde\eta_{k,n_k}^{t_k} - \tilde\eta_{n_k}^{t_k}\|_\infty + \|\tilde\eta_{n_k}^{t_k} - \eta^t\|_\infty,
\]
where we recall that $\tilde\eta_{n_k}^{t_k}$ denotes the $n_k$-th polygonal approximation of $\eta$ relative to time $t_k$, while $\tilde\eta_{k,n_k}^{t_k}$ denotes the $n_k$-th polygonal approximation of $\eta_k$ relative to time $t_k$. Then, from \eqref{eta_y_Lipschitz}, we obtain
\[
\|\tilde\eta_{k,n_k}^{t_k} - \tilde\eta_{n_k}^{t_k}\|_\infty \ \leq \ \max_j \big|\eta((t_j^{n_k}\wedge t_k) - t_k) - \eta((t_j^{n_k}\wedge t_k) - t_k)\big| \ \leq \ \|\eta - \eta_k\|_\infty.
\]
On the other hand, we notice that the term $\|\tilde\eta_{n_k}^{t_k} - \eta^t\|_\infty$ goes to zero as $k\rightarrow+\infty$, as a consequence of the uniform continuity of $\eta$ and of the fact that $\tilde\eta_{n_k}^{t_k}$ is a polygonal approximation of $\eta$ on $[-t_k,0]$. In conclusion, we deduce that $\|\tilde\eta_{k,n_k}^{t_k} - \eta^t\|_\infty$ goes to zero as $k\rightarrow+\infty$. This proves that the set $\mathcal K$ is relatively compact.

From the definition \eqref{b_sigma_F_tilde} of $\tilde b_n$, $\tilde \sigma_n$, $\tilde F_n$, it follows that they satisfy item (iii) of Theorem \ref{T:ExistSV_2}. As a matter of fact, fix a compact set $K\subset C([-T,0])$ and define the set $\mathcal K$ as in \eqref{compact_set_K}. Since $b$, $\sigma$, $F$ are continuous, they are uniformly continuous on the relatively compact set $[0,T]\times\mathcal K$. Therefore, by \eqref{b_sigma_F_tilde}, we see that there exists a continuity modulus $m_K$ (depending only on $b$, $\sigma$, $F$, and the compact set $K$) such that
\begin{equation}\label{equicontinuity_K}
|\tilde b_n(t,\eta) - \tilde b_n(t,\eta')| + |\tilde\sigma_n(t,\eta) - \tilde\sigma_n(t,\eta')| + |\tilde F_n(t,\eta) - \tilde F_n(t,\eta')| \ \leq \ m_K(\|\eta - \eta'\|_\infty),
\end{equation}
for all $t\in[0,T]$, $\eta,\eta'\in K$.

\vspace{1mm}

\noindent\textbf{Step IV.} \emph{Smooth approximation of $x\mapsto 1_{[0,t_j^n]}(x)$.} Our aim is to find, for every $n\in\N$ and $j=0,\ldots,N_n$, a sequence of functions $(\varphi_{n,j,k})_{k\in\N}\subset C^2([0,T])$, bounded uniformly with respect to $n$, $j$, $k$, with first derivatives bounded in $L^1([0,T])$ uniformly with respect to $n$, $j$, $k$, such that
\[
\eta((t_j^n\wedge t) - t) \ = \ \int_{[-t,0]}1_{[0,t_j^n]}(x + t)d^-\eta(x) \ = \ \lim_{k\rightarrow+\infty} \int_{[-t,0]}\varphi_{n,j,k}(x +  t)d^-\eta(x).
\]
We begin approximating the term $\int_{[-t,0]}1_{[0,t_j^n]}(x+t)d^-\eta(x)=\eta((t_j^n\wedge t)-t)$, for any $j=1,\ldots,N_n$. Let $\Phi(x)=\frac{1}{\sqrt{2\pi}}\int_{-\infty}^x \exp(-\frac{1}{2}z^2)dz$, for any $x\in\R$, be the cumulative distribution function of the standard Gaussian distribution. Then, we notice that, for any $t_j^n$, with $j\geq1$, and for any $t\in[0,T]$,
\[
\eta((t_j^n\wedge t) - t) \ = \ \int_{[-t,0]}1_{[0,t_j^n]}(x+t)d^-\eta(x) \ = \ \lim_{k\rightarrow+\infty} \int_{[-t,0]}\big[1 - \Phi(k(x + t - t_j^n))\big]d^-\eta(x).
\]
More precisely, we have (using the integration by parts formula (2.4) in \cite{cosso_russo15a})
\begin{align*}
&\sup_{t\in[0,T]}\bigg(\int_{[-t,0]}\big[1 - \Phi(k(x + t - t_j^n))\big]d^-\eta(x) - \eta((t_j^n\wedge t) - t)\bigg) \\
&= \ \sup_{t\in[0,T]}\bigg(\big[1 - \Phi(k(t - t_j^n))\big]\eta(0) + \int_{-t}^0 \eta(x) \frac{k}{\sqrt{2\pi}} e^{-\frac{1}{2}k^2(x + t - t_j^n)^2}dx - \eta((t_j^n\wedge t) - t)\bigg) \\
&= \ \sup_{t\in[0,T]}\bigg(\big[1 - \Phi(k(t - t_j^n))\big]\eta(0) + \int_{-kt_j^n}^{k(t - t_j^n)} \eta\Big(\frac{z}{k} + t_j^n - t\Big) \frac{1}{\sqrt{2\pi}} e^{-\frac{1}{2}z^2}dz - \eta((t_j^n\wedge t) - t)\bigg) \\
&= \ \sup_{t\in[0,T]}\bigg(\big[1 - \Phi(k(t - t_j^n))\big]\big[\eta(0) - \eta((t_j^n\wedge t) - t)\big] - \Phi(-k t_j^n) \eta((t_j^n\wedge t) - t) \\
&\quad \ + \int_{-kt_j^n}^{k(t - t_j^n)} \Big[\eta\Big(\frac{z}{k} + t_j^n - t\Big) - \eta((t_j^n\wedge t) - t)\Big] \frac{1}{\sqrt{2\pi}} e^{-\frac{1}{2}z^2}dz\bigg) \ \overset{k\rightarrow+\infty}{\longrightarrow} \ 0.
\end{align*}
It remains to consider the term corresponding to $j=0$, namely $\int_{[-t,0]}1_{[0,t_0^n]}(x+t)d^-\eta(x)=\eta(-t)$. We notice that we have (using the integration by parts formula (2.4) in \cite{cosso_russo15a})
\begin{align*}
&\sup_{t\in[0,T]}\bigg(\int_{[-t,0]}\big[1 - \Phi(k(x + t - 1/\sqrt{k}))\big]d^-\eta(x) - \eta(-t)\bigg) \\
&= \ \sup_{t\in[0,T]}\bigg(\big[1 - \Phi(k(t - 1/\sqrt{k}))\big]\eta(0) + \int_{-\sqrt{k}}^{k(t - 1/\sqrt{k})} \eta\Big(\frac{z}{k} + \frac{1}{\sqrt{k}} - t\Big) \frac{1}{\sqrt{2\pi}} e^{-\frac{1}{2}z^2}dz - \eta(-t)\bigg) \\
&= \ \sup_{t\in[0,T]}\bigg(\big[1 - \Phi(k(t - 1/\sqrt{k}))\big]\big[\eta(0) - \eta(-t)\big] - \Phi(-\sqrt{k}) \eta(-t) \\
&\quad \ + \int_{-\sqrt{k}}^{k(t - 1/\sqrt{k})} \Big[\eta\Big(\frac{z}{k} + \frac{1}{\sqrt{k}} - t\Big) - \eta(-t)\Big] \frac{1}{\sqrt{2\pi}} e^{-\frac{1}{2}z^2}dz\bigg) \ \overset{k\rightarrow+\infty}{\longrightarrow} \ 0.
\end{align*}
In conclusion, we obtain (the same property holds for the coefficients $\tilde\sigma_n$ and $\tilde F_n$)
\begin{align*}
\sup_{t\in[0,T]}\bigg[\hat b_n\bigg(t,\int_{[-t,0]}\big[1 - \Phi(k(x + t - 1/\sqrt{k}))\big]d^-\eta(x),\int_{[-t,0]}\big[1 - \Phi(k(x + t - t_1^n))\big]d^-\eta(x),\ldots\bigg)& \\
- \hat b_n\bigg(t,\int_{[-t,0]}1_{[0,t_0^n]}(x+t)d^-\eta(x),\ldots,\int_{[-t,0]}1_{[0,t_{N_n}^n]}(x+t)d^-\eta(x)\bigg)\bigg] \ \overset{k\rightarrow+\infty}{\longrightarrow} \ 0&,
\end{align*}
for every $\eta\in C([-T,0])$. Now, define
\[
\tilde b_{n,k}(t,\eta) \ := \ \hat b_n\bigg(t,\int_{[-t,0]}\big[1 - \Phi(k(x + t - 1/\sqrt{k}))\big]d^-\eta(x),\int_{[-t,0]}\big[1 - \Phi(k(x + t - t_1^n))\big]d^-\eta(x),\ldots\bigg),
\]
for all $n,k\in\N$, $t\in[0,T]$, $\eta\in C([-T,0])$. In a similar way, we define $\tilde\sigma_{n,k}$ and $\tilde F_{n,k}$. Proceeding along the same lines as in the previous Step III, we see that, for every compact set $K\in C([-T,0])$, we have (with the same continuity modulus as in \eqref{equicontinuity_K})
\[
|\tilde b_{n,k}(t,\eta) - \tilde b_{n,k}(t,\eta')| + |\tilde\sigma_{n,k}(t,\eta) - \tilde\sigma_{n,k}(t,\eta')| + |\tilde F_{n,k}(t,\eta) - \tilde F_{n,k}(t,\eta')| \ \leq \ m_K(\|\eta - \eta'\|_\infty),
\]
for all $t\in[0,T]$, $\eta,\eta'\in K$. In other words, the functions $\tilde b_{n,k}(t,\cdot)$, $\tilde\sigma_{n,k}(t,\cdot)$, $\tilde F_{n,k}(t,\cdot)$, $n,k\in\N$, are equicontinuous on compact sets, uniformly with respect to $t\in[0,T]$. Then, it follows from Lemma \ref{L:StabilityApp} that there exists a subsequence $(\tilde b_{n,k_n},\tilde\sigma_{n,k_n},\tilde F_{n,k_n})_{n\in\N}$ which converges pointwise to $(b,\sigma,F)$ on $[0,T]\times C([-T,0])$.

From now on, to alleviate the notation, we denote the subsequence $(\tilde b_{n,k_n},\tilde\sigma_{n,k_n},\tilde F_{n,k_n})_{n\in\N}$ simply by $(\tilde b_n,\tilde\sigma_n,\tilde F_n)_{n\in\N}$, with
\begin{align*}
\tilde b_n(t,\eta) \ &= \ \hat b_n\bigg(t,\int_{[-t,0]}\varphi_1(x+t)d^-\eta(x),\ldots,\int_{[-t,0]}\varphi_{N_n}(x+t)d^-\eta(x)\bigg), \\
\tilde\sigma_n(t,\eta) \ &= \ \hat\sigma_n\bigg(t,\int_{[-t,0]}\varphi_1(x+t)d^-\eta(x),\ldots,\int_{[-t,0]}\varphi_{N_n}(x+t)d^-\eta(x)\bigg), \\
\tilde F_n(t,\eta) \ &= \ \hat F_n\bigg(t,\int_{[-t,0]}\varphi_1(x+t)d^-\eta(x),\ldots,\int_{[-t,0]}\varphi_{N_n}(x+t)d^-\eta(x)\bigg),
\end{align*}
for every $(t,\eta)\in[0,T]\times C([-T,0])$, for some sequence of functions $(\varphi_j)_{j\in\N}\subset C^2([0,T])$, bounded uniformly with respect to $j$, with first derivatives bounded in $L^1([0,T])$ uniformly with respect to $j$.

\vspace{1mm}

\noindent\textbf{Step V.} \emph{Conclusion.} In order to conclude the proof, it remains to perform a smooth approximation of $\hat b_n$, $\hat\sigma_n$, $\hat F_n$. To this end, we proceed along the same lines as in Step III of the proof of Theorem \ref{T:ExistSV}. More precisely, for every $n\in\N$, we consider the function $\rho_n\in C^\infty(\R^{N_n})$ given by \eqref{rho}. As in Step III of the proof of Theorem \ref{T:ExistSV}, we set $\rho_{n,k}(\mathbf z) := k^{N_n}\rho_n(k\,\mathbf z)$, $\forall\,\mathbf z\in\R^{N_n}$, $k\in\N$, and we define
\[
\hat b_{n,k}(t,\mathbf y) \ = \ \int_{\R^{N_n}} \rho_{n,k}(\mathbf z) \hat b_n(t,\mathbf y - \mathbf z) d\mathbf z,
\]
for all $(t,\mathbf y)\in[0,T]\times\R^{N_n}$. In a similar way, we define $\hat\sigma_{n,k}$ and $\hat F_{n,k}$. We denote
\[
\check b_{n,k}(t,\eta) \ := \ \hat b_{n,k}\bigg(t,\int_{[-t,0]}\varphi_1(x+t)d^-\eta(x),\ldots,\int_{[-t,0]}\varphi_{N_n}(x+t)d^-\eta(x)\bigg),
\]
for all $n,k\in\N$, $t\in[0,T]$, $\eta\in C([-T,0])$. We define similarly $\check\sigma_{n,k}$ and $\check F_{n,k}$. By Lemma \ref{L:StabilityApp}, we deduce that there exists a subsequence $(\check b_{n,k_n},\check\sigma_{n,k_n},\check F_{n,k_n})_{n\in\N}$ which converges pointwise to $(b,\sigma,F)$. Then, we define $\bar b_n:=\hat b_{n,k_n}$, $\bar\sigma_n:=\hat\sigma_{n,k_n}$, $\bar F_n:=\hat F_{n,k_n}$, and also $b_n:=\check b_{n,k_n}$, $\sigma_n:=\check\sigma_{n,k_n}$, $F_n:=\check F_{n,k_n}$. We see that the sequences $(b_n,\sigma_n,F_n)_{n\in\N}$ and $(\bar b_n,\bar\sigma_n,\bar F_n)_{n\in\N}$ satisfy items \textup{(i)}, \textup{(ii)'}, \textup{(iii)}, \textup{(iv), \textup{(v)}} of Theorem \ref{T:ExistSV_2}. This concludes the proof.
\ep

\subsection{The Markovian case revisited}

In the present section we show that, in the Markovian case, Definition \ref{D:Strong-Visc} of (path-dependent) strong-viscosity solution is coherent with Definition \ref{D:ViscosityFinite}. In particular, consider the semilinear parabolic PDE
\begin{equation}
\label{KolmEq_Markov_2}
\begin{cases}
\partial_t u(t,x) + \langle b(t,x),D_x u(t,x)\rangle + \frac{1}{2}\textup{tr}(\sigma\sigma\trans(t,x)D_x^2 u(t,x)) & \\
\hspace{2.8cm}+\, f(t,x,u(t,x),\sigma\trans(t,x)D_x u(t,x)) \ = \ 0, &\forall\,(t,x)\in[0,T)\times\R, \\
u(T,x) \ = \ h(x), &\forall\,x\in\R,
\end{cases}
\end{equation}
with $b$, $\sigma$, $f$, $h$ satisfying Assumption {\bf (A0)}. Then, equation \eqref{KolmEq_Markov_2} can be rewritten as the following semilinear parabolic path-dependent PDE
\begin{equation}
\label{KolmEq_2}
\begin{cases}
\partial_t \Uc + D^H \Uc + \tilde b(t,\eta)D^V \Uc + \frac{1}{2}\tilde\sigma(t,\eta)^2D^{VV} \Uc \\
\hspace{2.3cm}+\, \tilde F(t,\eta,\Uc,\sigma(t,\eta)D^V \Uc) \ = \ 0, \;\;\; &\quad\forall\,(t,\eta)\in[0,T[\times C([-T,0]), \\
\Uc(T,\eta) \ = \ \tilde H(\eta), &\quad\forall\,\eta\in C([-T,0]),
\end{cases}
\end{equation}
where $\tilde b\colon[0,T]\times C([-T,0])\rightarrow\R$, $\tilde\sigma\colon[0,T]\times C([-T,0])\rightarrow\R$, $\tilde F\colon[0,T]\times C([-T,0])\times\R\times\R\rightarrow\R$, $\tilde H\colon C([-T,0])\rightarrow\R$ are defined as 
\[
\tilde b(t,\eta) \, := \, b(t,\eta(0)), \quad\; \tilde\sigma(t,\eta) \, := \, \sigma(t,\eta(0)), \quad\; \tilde F(t,\eta,y,z) \, := \, f(t,\eta(0),y,z), \quad\; \tilde H(\eta) \, := \, h(\eta(0)),
\]
for every $(t,\eta,y,z)\in[0,T]\times C([-T,0])\times\R\times\R$. Since $b$, $\sigma$, $f$, $h$ satisfy Assumption {\bf (A0)}, it follows that $\tilde b$, $\tilde\sigma$, $\tilde F$, $\tilde H$ satisfy Assumption {\bf (A1)}. Moreover, we have the following result.

\begin{Proposition}
Suppose that Assumption {\bf (A0)} holds.
\begin{enumerate}
\item[(1)] Every strong-viscosity solution $u\colon[0,T]\times\R\rightarrow\R$ to equation \eqref{KolmEq_Markov_2}, in the sense of Definition \ref{D:ViscosityFinite}, is such that the map $\Uc\colon[0,T]\times C([-T,0])\rightarrow\R$, defined by
\[
\Uc(t,\eta) \ := \ u(t,\eta(0)), \qquad \forall\,(t,\eta)\in[0,T]\times C([-T,0]),
\]
is a (path-dependent) strong-viscosity solution to equation \eqref{KolmEq_2} in the sense of Definition \ref{D:Strong-Visc}.
\item[(2)] Viceversa, every (path-dependent) strong-viscosity solution $\Uc\colon[0,T]\times C([-T,0])\rightarrow\R$ to equation \eqref{KolmEq_2}, in the sense of Definition \ref{D:Strong-Visc}, can be represented as 
\[
\Uc(t,\eta) \ = \ u(t,\eta(0)), \qquad \forall\,(t,\eta)\in[0,T]\times C([-T,0]),
\]
for some function $u\colon[0,T]\times\R\rightarrow\R$. Moreover, under the assumptions of either Theorem \ref{T:ExistSV_Markov} or Theorem \ref{T:ExistSV_Markov2}, $u$ is a strong-viscosity solution to equation \eqref{KolmEq_Markov_2} in the sense of Definition \ref{D:ViscosityFinite}.
\end{enumerate}
\end{Proposition}
\textbf{Proof.}
\emph{Proof of point (1).} Let $u\colon[0,T]\times\R\rightarrow\R$ be a strong-viscosity solution to equation \eqref{KolmEq_Markov_2} in the sense of Definition \ref{D:ViscosityFinite}, so that there exists a sequence $(u_n,h_n,f_n,b_n,\sigma_n)_n$ of Borel measurable functions $u_n\colon[0,T]\times\R\rightarrow\R$, $h_n\colon\R\rightarrow\R$, $f_n\colon[0,T]\times\R\times\R\times\R\rightarrow\R$, $b_n\colon[0,T]\times\R\rightarrow\R$, $\sigma_n\colon[0,T]\times\R\rightarrow\R$, satisfying points (i)-(ii)-(iii) of Definition \ref{D:ViscosityFinite}. Define $\Uc\colon[0,T]\times C([-T,0])\rightarrow\R$ as 
$\Uc(t,\eta):=u(t,\eta(0))$. Let us prove that $\Uc$ is a (path-dependent) strong-viscosity solution to equation \eqref{KolmEq_2} in the sense of Definition \ref{D:Strong-Visc}. For every $n$, define $\tilde b_n\colon[0,T]\times C([-T,0])\rightarrow\R$, $\tilde\sigma_n\colon[0,T]\times C([-T,0])\rightarrow\R$, $\tilde F_n\colon[0,T]\times C([-T,0])\times\R\times\R\rightarrow\R$, $\tilde H_n\colon C([-T,0])\rightarrow\R$ as
\[
\tilde b_n(t,\eta) := b_n(t,\eta(0)), \; \tilde\sigma_n(t,\eta) := \sigma_n(t,\eta(0)), \; \tilde F_n(t,\eta,y,z) := f_n(t,\eta(0),y,z), \; \tilde H_n(\eta) := h_n(\eta(0)),
\]
for every $(t,\eta,y,z)\in[0,T]\times C([-T,0])\times\R\times\R$. Moreover, let $\Uc_n\colon[0,T]\times C([-T,0])\rightarrow\R$ be given by $\Uc_n(t,\eta)=u_n(t,\eta(0))$. Then, the sequence $(\Uc_n,\tilde H_n,\tilde F_n,\tilde b_n,\tilde\sigma_n)_n$ satisfies points (i)-(ii)-(iii) of Definition \ref{D:Strong-Visc}, from which it follows that $\Uc$ is a (path-dependent) strong-viscosity solution to equation \eqref{KolmEq_2} in the sense of Definition \ref{D:Strong-Visc}.

\vspace{1mm}

\noindent\emph{Proof of point (2).} We begin recalling that, since Assumption {\bf (A1)} holds, by Theorem \ref{T:UniqSV} we have that $\Uc$ is given by
\begin{equation}\label{FK_Markov_Uc}
\Uc(t,\eta) \ = \ Y_t^{t,\eta}, \qquad \forall\,(t,\eta)\in[0,T]\times C([-T,0]),
\end{equation}
where $(Y_s^{t,\eta},Z_s^{t,\eta})_{s\in[t,T]}\in\S^2(t,T)\times\H^2(t,T)$ is the unique solution to the equation
\[
Y_s^{t,\eta} \ = \ \tilde H(\X_T^{t,\eta}) + \int_s^T \tilde F(r,\X_r^{t,\eta},Y_r^{t,\eta},Z_r^{t,\eta}) dr - \int_s^T Z_r^{t,\eta} dW_r, \qquad s\in[t,T],
\]
with $\X^{t,\eta}$ window process of $X^{t,\eta}=(X_s^{t,\eta})_{s\in[-T+t,T]}$, solution of the equation
\[
\begin{cases}
dX_s^{t,\eta} = \ \tilde b(s,\mathbb X_s^{t,\eta})dt + \tilde\sigma(s,\mathbb X_s^{t,\eta})dW_s, \qquad\qquad & s\in[t,T], \\
X_s^{t,\eta} \ = \ \eta(s-t), & s\in[-T+t,t].
\end{cases}
\]
From the definition of $\X^{t,\eta}$, we see that $\X_r^{t,\eta}(0)=X_r^{t,\eta}$, for every $r\in[t,T]$. Therefore, the equation solved by $X^{t,\eta}$ can be written as follows:
\[
\begin{cases}
dX_s^{t,\eta} = \ b(s,X_s^{t,\eta})dt + \sigma(s,X_s^{t,\eta})dW_s, \qquad\qquad & s\in[t,T], \\
X_s^{t,\eta} \ = \ \eta(s-t), & s\in[-T+t,t].
\end{cases}
\]
We see that $(X_s^{t,\eta})_{s\in[t,T]}$
 solves the same equation of the process $X^{t,\eta(0)}=(X_s^{t,\eta(0)})_{s\in[t,T]}$, namely equation \eqref{SDE_Markov} with $x=\eta(0)$. As a consequence, the processes $(X_s^{t,\eta})_{s\in[t,T]}$ and $X^{t,\eta(0)}=(X_s^{t,\eta(0)})_{s\in[t,T]}$ are $\P$-indistinguishable. From this result we deduce that the equation solved by $(Y_s^{t,\eta},Z_s^{t,\eta})_{s\in[t,T]}$ can be written as 
\[
Y_s^{t,\eta} \ = \ h(X_T^{t,\eta(0)}) + \int_s^T f(r,X_r^{t,\eta(0)},Y_r^{t,\eta},Z_r^{t,\eta}) dr - \int_s^T Z_r^{t,\eta} dW_r, \qquad t\leq s\leq T.
\]
We notice that $(Y_s^{t,\eta},Z_s^{t,\eta})_{s\in[t,T]}$ solves the same equation of $(Y_s^{t,\eta(0)},Z_s^{t,\eta(0)})_{s\in[t,T]}$, namely equation \eqref{BSDE_Markov2} with $x=\eta(0)$. By uniqueness, it follows that $\|Y^{t,\eta}-Y^{t,\eta(0)}\|_{_{\S^2(t,T)}}=0$ and $\|Z^{t,\eta}-Z^{t,\eta(0)}\|_{_{\H^2(t,T)}}=0$. Therefore, we deduce $Y_t^{t,\eta}=Y_t^{t,\eta(0)}$. Now, define the map $u\colon[0,T]\times\R\rightarrow\R$ as 
\begin{equation}
\label{FK_Markov_3}
u(t,x) \ := \ Y_t^{t,x}, \qquad \forall\,(t,x)\in[0,T]\times\R.
\end{equation}
Then, the following relation between the map $\Uc$ in \eqref{FK_Markov_Uc} and the function $u$ in \eqref{FK_Markov_3} holds:
\[
\Uc(t,\eta) \ = \ u(t,\eta(0)), \qquad \forall\,(t,\eta)\in[0,T]\times C([-T,0]).
\]
Finally, under the assumptions of either Theorem \ref{T:ExistSV_Markov} or Theorem \ref{T:ExistSV_Markov2}, we deduce that $u$ is a strong-viscosity solution to equation \eqref{KolmEq_Markov_2} in the sense of Definition \ref{D:ViscosityFinite}.
\ep

\appendix

\renewcommand\thesection{Appendix}

\section{}

In the present appendix we fix a complete probability space $(\Omega,\Fc,\P)$ on which a $d$-dimensional Brownian motion $W=(W_t)_{t\geq0}$ is defined. We denote $\F=(\Fc_t)_{t\geq0}$ the completion of the natural filtration generated by $W$.

\renewcommand\thesection{\Alph{subsection}}

\renewcommand\thesubsection{\Alph{subsection}.}

\subsection{Estimates for path-dependent stochastic differential equations}

Let $C([-T,0];\R^d)$ denote the Banach space of all continuous paths $\eta\colon[-T,0]\rightarrow\R^d$ endowed with the supremum norm $\|\eta\|=\sup_{t\in[0,T]}|\eta(t)|$. Notice that, when $d=1$, we simply write $C([-T,0])$ instead of $C([-T,0];\R)$. In the present section we consider the $d$-dimensional path-dependent SDE:
\begin{equation}
\label{SDE_d-dim}
\begin{cases}
dX_s = \ b(s,\mathbb X_s)dt + \sigma(s,\mathbb X_s)dW_s, \qquad\qquad & s\in[t,T], \\
X_s \ = \ \eta(s-t), & s\in[-T+t,t],
\end{cases}
\end{equation}
where on $b\colon[0,T]\times C([-T,0];\R^d)\rightarrow\R^d$ and $\sigma\colon[0,T]\times C([-T,0];\R^d)\rightarrow\R^d$ we shall impose the following assumptions.

\vspace{3mm}

\noindent\textbf{(A$_{b,\sigma}$)} \hspace{3mm} $b$ and $\sigma$ are Borel measurable functions satisfying, for some positive constant $C$,
\begin{align*}
|b(t,\eta)-b(t,\eta')| + |\sigma(t,\eta)-\sigma(t,\eta')| \ &\leq \ C\|\eta-\eta'\|, \\
|b(t,0)| + |\sigma(t,0)| \ &\leq \ C,
\end{align*}
for all $t\in[0,T]$ and $\eta,\eta'\in C([-T,0];\R^d)$.

\vspace{3mm}

Notice that equation \eqref{SDE_d-dim} on $[t,T]$ becomes equation \eqref{SDE_Markov} when $b=b(t,\eta)$ and $\sigma=\sigma(t,\eta)$ are non-path-dependent, so that they depend only on $\eta(t)$ at time $t$. On the other hand, when $d=1$ equation \eqref{SDE_d-dim} reduces to equation \eqref{SDE}.

\begin{Lemma}
\label{L:SDE}
Under Assumption {\bf (A$_{b,\sigma}$)}, for any $(t,\eta)\in[0,T]\times C([-T,0];\R^d)$ there exists a unique $($up to indistinguishability$)$ $\F$-adapted continuous process $X^{t,\eta}=(X_s^{t,\eta})_{s\in[-T+t,T]}$ strong solution to equation \eqref{SDE_Markov}. Moreover, for any $p\geq1$ there exists a positive constant $C_p$ such that
\begin{equation}
\label{EstimateSupX}
\E\Big[\sup_{s\in[-T+t,T]}\big|X_s^{t,\eta}\big|^p\Big] \ \leq \ C_p \big(1 + \|\eta\|^p\big).
\end{equation}
\end{Lemma}
\textbf{Proof.}
Existence and uniqueness follow from Theorem 14.23 in \cite{jacod79}. Concerning estimate \eqref{EstimateX} we refer to Proposition 3.1 in \cite{cosso_russo15a} (notice that in \cite{cosso_russo15a}, estimate \eqref{EstimateX} is proved for the case $d=1$; however, proceeding along the same lines, we can prove \eqref{EstimateX} for a generic $d\in\N\backslash\{0\}$).
\ep

\begin{Lemma}
\label{L:AppendixX}
Suppose that Assumption {\bf (A$_{b,\sigma}$)} holds and let $(b_n,\sigma_n)_n$ be a sequence satisfying Assumption {\bf (A$_{b,\sigma}$)} with a positive constant $C$ independent of $n$. Moreover, $(b_n,\sigma_n)$ converges pointwise to $(b,\sigma)$ as $n\rightarrow\infty$. For any $n\in\N$ and $(t,\eta)\in[0,T]\times C([-T,0];\R^d)$, denote by $X^{n,t,\eta}=(X_s^{n,t,\eta})_{s\in[-T+t,T]}$ the unique solution to the  path-dependent SDE
\begin{equation}
\label{Xn}
\begin{cases}
dX_s^{n,t,\eta} = \ b_n(s,\mathbb X_s^{n,t,\eta})dt + \sigma_n(s,\mathbb X_s^{n,t,\eta})dW_s, \qquad\qquad & t\leq s\leq T, \\
X_s^{n,t,\eta} \ = \ \eta(s-t), & -T+t\leq s\leq t.
\end{cases}
\end{equation}
Then, for every $p\geq1$, we have
\begin{equation}
\label{EstimateSupXn}
\E\Big[\sup_{t\leq s\leq T} |X_s^{n,t,\eta}|^p\Big] \ \leq \ C_p \big(1 + \|\eta\|^p\big), \qquad \forall\,(t,\eta)\in [0,T]\times C([-T,0];\R^d),\,\forall\,n\in\N,
\end{equation}
for some positive constant $C_p$, and
\begin{equation}
\label{limXn-X}
\lim_{n\rightarrow\infty} \E\Big[\sup_{t\leq s\leq T} |X_s^{n,t,\eta} - X_s^{t,\eta}|^p\Big] \ = \ 0, \qquad \forall\,(t,\eta)\in [0,T]\times C([-T,0];\R^d).
\end{equation}
\end{Lemma}
\textbf{Proof.}
For any $n\in\N$ and $(t,\eta)\in[0,T]\times C([-T,0];\R^d)$, the existence and uniqueness of $(X_s^{n,t,\eta})_{s\in[-T+t,T]}$, as well as estimate \eqref{EstimateSupXn}, can be proved proceeding as in Lemma \ref{L:SDE}. It remains to prove \eqref{limXn-X}. Observe that
\[
X_s^{n,t,\eta} - X_s^{t,\eta} \ = \ \int_t^s \big(b_n(r,\X_r^{n,t,\eta}) - b(r,\X_r^{t,\eta})\big) dr + \int_t^s \big(\sigma_n(r,\X_r^{n,t,\eta}) - \sigma(r,\X_r^{t,\eta})\big) dW_r.
\]
Then, taking the $p$-th power, we get (recalling the standard inequality $(a+b)^p\leq2^{p-1}(a^p+b^p)$, for any $a,b\in\R$) that $|X_s^{n,t,\eta} - X_s^{t,\eta}|^p$ is less than or equal to
\[
2^{p-1}\bigg|\int_t^s \big(b_n(r,\X_r^{n,t,\eta}) - b(r,\X_r^{t,\eta})\big) dr\bigg|^p + 2^{p-1}\bigg|\int_t^s \big(\sigma_n(r,\X_r^{n,t,\eta}) - \sigma(r,\X_r^{t,\eta})\big) dW_r\bigg|^p.
\]
 In the sequel we shall denote $c_p$ a generic positive constant which may change from line to line, independent of $n$, depending only on $T$, $p$, and the Lipschitz constant of $b_n,\sigma_n$. Taking the supremum over $s\in[t,T]$, and 
applying H\"older's inequality to the drift term, we get
\begin{align}
\label{E:SupProof2}
\|\X_s^{n,t,\eta} - \X_s^{t,\eta}\|^p \ &\leq \ c_p\int_t^s \big|b_n(r,\X_r^{n,t,\eta}) - b(r,\X_r^{t,\eta})\big|^p dr \notag \\
&\quad \ + 2^{p-1}\sup_{t \leq u \leq s}\bigg|\int_t^u \big(\sigma_n(r,\X_r^{n,t,\eta}) - \sigma(r,\X_r^{t,\eta})\big) dW_r\bigg|^p.
\end{align}
Notice that
\begin{align}
\label{E:bn-b}
&\int_t^s \big|b_n(r,\X_r^{n,t,\eta}) - b(r,\X_r^{t,\eta})\big|^p dr \notag \\
&\leq \ 2^{p-1}\int_t^s \big|b_n(r,\X_r^{n,t,\eta}) - b_n(r,\X_r^{t,\eta})\big|^p dr + 2^{p-1}\int_t^s \big|b_n(r,\X_r^{t,\eta}) - b(r,\X_r^{t,\eta})\big|^p dr \notag \\
&\leq \ c_p\int_t^s \|\X_r^{n,t,\eta} - \X_r^{t,\eta}\|^p dr + 2^{p-1}\int_t^s \big|b_n(r,\X_r^{t,\eta}) - b(r,\X_r^{t,\eta})\big|^p dr.
\end{align}
In addition, from Burkholder-Davis-Gundy inequality we have
\begin{align}
\label{E:sigman-sigma}
&\E\bigg[\sup_{t \leq u \leq s}\bigg|\int_t^u \big(\sigma_n(r,\X_r^{n,t,\eta}) - \sigma(r,\X_r^{t,\eta})\big) dW_r\bigg|^p\bigg] \ \leq \ c_p \E\bigg[\int_t^s \big|\sigma_n(r,\X_r^{n,t,\eta})-\sigma(r,\X_r^{t,\eta})\big|^{p/2} dr\bigg] \notag \\
&\leq \ c_p \E\bigg[\int_t^s \big|\sigma_n(r,\X_r^{n,t,\eta})-\sigma_n(r,\X_r^{t,\eta})\big|^{p/2} dr\bigg] + c_p \E\bigg[\int_t^s \big|\sigma_n(r,\X_r^{t,\eta})-\sigma(r,\X_r^{t,\eta})\big|^{p/2} dr\bigg] \notag \\
&\leq \ c_p\int_t^s \E\big[\|\X_r^{n,t,\eta} - \X_r^{t,\eta}\|^p\big] dr + c_p \int_t^s \E\big[\big|\sigma_n(r,\X_r^{t,\eta})-\sigma(r,\X_r^{t,\eta})\big|^{p/2}\big] dr.
\end{align}
Taking the expectation in \eqref{E:SupProof2}, and using \eqref{E:bn-b} and \eqref{E:sigman-sigma}, we find
\begin{align*}
\E\big[\|\X_s^{n,t,\eta} - \X_s^{t,\eta}\|^p\big] \ &\leq \ c_p\int_t^s \E\big[\|\X_r^{n,t,\eta} - \X_r^{t,\eta}\|^p\big] dr + c_p \int_t^s \E\big[\big|b_n(r,\X_r^{t,\eta}) - b(r,\X_r^{t,\eta})\big|^p\big] dr \\
&\quad \ + c_p \int_t^T \E\big[\big|\sigma_n(r,\X_r^{t,\eta})-\sigma(r,\X_r^{t,\eta})\big|^{p/2}\big] dr.
\end{align*}
Then, applying Gronwall's lemma to the map $r\mapsto\E[\|\X_r^{n,t,\eta} - \X_r^{t,\eta}\|^p]$, we get
\begin{align*}
\E\Big[\sup_{t\leq s\leq T} |X_s^{n,t,\eta} - X_s^{t,\eta}|^p\Big] \ &\leq \ c_p \int_t^T \E\big[\big|b_n(r,\X_r^{t,\eta}) - b(r,\X_r^{t,\eta})\big|^p\big] dr \\
&\quad \ + c_p \int_t^T \E\big[\big|\sigma_n(r,\X_r^{t,\eta})-\sigma(r,\X_r^{t,\eta})\big|^{p/2}\big] dr.
\end{align*}
In conclusion, \eqref{limXn-X} follows from estimate \eqref{EstimateSupX} and Lebesgue's dominated convergence theorem.
\ep

\setcounter{Theorem}{0}
\setcounter{equation}{0}

\subsection{Estimates for backward stochastic differential equations}

We derive estimates for the norm of the $Z$ and $K$ components for supersolutions to backward stochastic differential equations, in terms of the norm of the $Y$ component. These results are standard, but seemingly not at disposal in the following form in the literature. Firstly, let us introduce a generator function $F\colon[0,T]\times\Omega\times\R\times\R^d\rightarrow\R$ satisfying the usual assumptions:
\begin{enumerate}
\item[(A.a)] $F(\cdot,y,z)$ is $\F$-predictable for every $(y,z)\in\R\times\R^d$.
\item[(A.b)] There exists a positive constant $C_F$ such that
\[
|F(s,y,z) - F(s,y',z')| \ \leq \ C_F\big(|y-y'| + |z-z'|\big),
\]
for all $y,y'\in\R$, $z,z'\in\R^d$, $ds\otimes d\P$-a.e.
\item[(A.c)] Integrability condition:
\[
\E\bigg[\int_t^T |F(s,0,0)|^2 ds\bigg] \ \leq \ M_F,
\]
for some positive constant $M_F$.
\end{enumerate}

\begin{Proposition}
\label{P:EstimateBSDEAppendix}
For any $t,T\in\R_+$, $t<T$, consider $(Y_s,Z_s,K_s)_{s\in[t,T]}$ satisfying the following.
\begin{enumerate}
\item[\textup{(i)}] $Y\in\S^2(t,T)$ and it is continuous.
\item[\textup{(ii)}] $Z$ is an $\R^d$-valued $\F$-predictable process such that $\P(\int_t^T |Z_s|^2 ds<\infty)=1$.
\item[\textup{(iii)}] $K$ is a real nondecreasing $($or nonincreasing$)$ continuous $\F$-predictable process such that $K_t = 0$.
\end{enumerate}
Suppose that $(Y_s,Z_s,K_s)_{s\in[t,T]}$ solves the BSDE, $\P$-a.s.,
\begin{equation}
\label{E:BSDE_Appendix}
Y_s \ = \ Y_T + \int_s^T F(r,Y_r,Z_r) dr + K_T - K_s - \int_s^T \langle Z_r, dW_r\rangle, \qquad t \leq s \leq T,
\end{equation}
for some generator function $F$ satisfying conditions \textup{(A.b)-(A.c)}. Then $(Z,K)\in\H^2(t,T)^d\times\A^{+,2}(t,T)$ and
\begin{equation}
\label{EstimateBSDE1}
\|Z\|_{\H^2(t,T)^d}^2 + \|K\|_{\S^2(t,T)}^2 \ \leq \ C\bigg(\|Y\|_{\S^2(t,T)}^2 + \E\int_t^T |F(s,0,0)|^2 ds\bigg),
\end{equation}
for some positive constant $C$ depending only on $T$ and $C_F$, the Lipschitz constant of $F$. If in addition $K\equiv0$, we have the standard estimate
\begin{equation}
\label{EstimateBSDE2}
\|Y\|_{\S^2(t,T)}^2 + \|Z\|_{\H^2(t,T)^d}^2 \ \leq \ C'\bigg(\E\big[|Y_T|^2\big] + \E\int_t^T |F(s,0,0)|^2 ds\bigg),
\end{equation}
for some positive constant $C'$ depending only on $T$ and $C_F$.
\end{Proposition}
\textbf{Proof.}
\emph{Proof of estimate \eqref{EstimateBSDE1}.} Let us consider the case where $K$ is nondecreasing. For every $k\in\N$, define the stopping time
\[
\tau_k \ = \ \inf\bigg\{s\geq t\colon \int_t^s |Z_r|^2 dr \geq k \bigg\} \wedge T.
\]
Then, the local martingale $(\int_t^s Y_r\langle 1_{[t,\tau_k]}(r)Z_r,dW_r\rangle)_{s\in[t,T]}$ satisfies, using Burkholder-Davis-Gundy inequality,
\[
\E\bigg[\sup_{t \leq s \leq T}\bigg|\int_t^s Y_r\langle 1_{[t,\tau_k]}(r)Z_r,dW_r\rangle\bigg|\bigg] \ < \ \infty,
\]
therefore it is a martingale. As a consequence, an application of It\^o's formula to $|Y_s|^2$ between $t$ and $\tau_k$ yields
\begin{align}
\label{E:Proof_ItoAppendix}
\E\big[|Y_t|^2\big] + \E\int_t^{\tau_k} |Z_r|^2 dr \ &= \ \E\big[|Y_{\tau_k}|^2\big] + 2\E\int_t^{\tau_k} Y_r F(r,Y_r,Z_r) dr + 2 \E\int_t^{\tau_k} Y_r dK_r.
\end{align}
In the sequel $c$ and $c'$ will be two strictly positive constants depending only on $C_F$, the Lipschitz constant of $F$.
 Using (A.b) and recalling the standard inequality $ab \leq a^2 + b^2/4$, for any $a,b\in\R$, we see that
\begin{equation}
\label{E:Proof_LipschitzAppendix}
2\E\int_t^{\tau_k} Y_r F(r,Y_r,Z_r) dr \ \leq \ cT\|Y\|_{\S^2(t,T)}^2 + \frac{1}{4}\E\int_t^{\tau_k} |Z_r|^2 dr + \E\int_t^T |F(r,0,0)|^2 dr.
\end{equation}
Regarding the last term on the right-hand side in \eqref{E:Proof_ItoAppendix}, for every $\eps>0$, recalling the standard inequality $2ab \leq \eps a^2 + b^2/\eps$, for any $a,b\in\R$, we have
\begin{equation}
\label{E:YdK_Appendix}
2\E\int_t^{\tau_k} Y_r dK_r \ \leq \ \frac{1}{\eps}\|Y\|_{\S^2(t,T)}^2 + \eps\E\big[|K_{\tau_k}|^2\big].
\end{equation}
Now, from \eqref{E:BSDE_Appendix} we get
\[
K_{\tau_k} \ = \ Y_t - Y_{\tau_k} - \int_t^{\tau_k} F(r,Y_r,Z_r) dr + \int_t^{\tau_k} \langle Z_r,dW_r\rangle.
\]
Therefore, recalling that $(x_1+\cdots+x_4)\leq4(x_1^2+\cdots+x_4^2)$, for any $x_1,\ldots,x_4\in\R$
\[
\E\big[|K_{\tau_k}|^2\big] \ \leq \ 8\|Y\|_{\S^2(t,T)}^2 + 4T\E\int_t^{\tau_k} |F(r,Y_r,Z_r)|^2 dr + 4\E\bigg|\int_t^{\tau_k} \langle Z_r,dW_r\rangle\bigg|^2.
\]
From It\^o's isometry and (A.b), we obtain
\begin{equation}
\label{E:K<C_Appendix}
\E\big[|K_{\tau_k}|^2\big] \ \leq \ c'(1+T^2)\|Y\|_{\S^2(t,T)}^2 + c'(1+T)\E\int_t^{\tau_k} |Z_r|^2 dr + c'T\E\int_t^T |F(r,0,0)|^2 dr.
\end{equation}
Then, taking $\eps=1/(4c'(1+T))$ in \eqref{E:YdK_Appendix} we get
\begin{align}
\label{E:YdK2_Appendix}
2\E\int_t^{\tau_k}\!\! Y_r dK_r \ &\leq \ \frac{16c'(1+T)^2+1+T^2}{4(1+T)}\|Y\|_{\S^2(t,T)}^2 + \frac{1}{4}\E\int_t^{\tau_k}\!\! |Z_r|^2 dr + \frac{T}{4(1+T)}\E\int_t^T\!\! |F(r,0,0)|^2 dr \notag \\
&\leq \ c(1+T^2)\|Y\|_{\S^2(t,T)}^2 + \frac{1}{4}\E\int_t^{\tau_k}\!\! |Z_r|^2 dr + cT\E\int_t^T\!\! |F(r,0,0)|^2 dr.
\end{align}
Plugging \eqref{E:Proof_LipschitzAppendix} and \eqref{E:YdK2_Appendix} into \eqref{E:Proof_ItoAppendix}, we end up with
\[
\E\big[|Y_{\tau_k}|^2\big] + \frac{1}{2}\E\int_t^{\tau_k} |Z_r|^2 dr \ \leq \ c(1+T^2)\|Y\|_{\S^2(t,T)}^2 + c(1+T)\E\int_t^T |F(r,0,0)|^2 dr.
\]
Then, from monotone convergence theorem,
\begin{equation}
\label{E:Z<C_Appendix}
\E\int_t^T |Z_r|^2 dr \ \leq \ c(1+T^2)\|Y\|_{\S^2(t,T)}^2 + c(1+T)\E\int_t^T |F(r,0,0)|^2 dr.
\end{equation}
Plugging \eqref{E:Z<C_Appendix} into \eqref{E:K<C_Appendix}, and using again monotone convergence theorem, we finally obtain
\[
\|K\|_{\S^2(t,T)}^2 \ = \ \E\big[|K_T|^2\big] \ \leq \ c(1+T^3)\|Y\|_{\S^2(t,T)}^2 + c(1+T^2)\E\int_t^T |F(r,0,0)|^2 dr.
\]

\vspace{1mm}

\noindent\emph{Proof of estimate \eqref{EstimateBSDE2}.} The proof of this estimate is standard, see, e.g., Remark (b) immediately after Proposition 2.1 in \cite{elkaroui_peng_quenez97}. We just recall that it can be done in the following \emph{two steps}: \emph{first}, we apply It\^o's formula to $|Y_s|^2$, afterwards we take the expectation, then we use the Lipschitz property of $F$ with respect to $(y,z)$, and finally we apply Gronwall's lemma to the map $v(s):=\E[|Y_s|^2]$, $s\in[t,T]$. Then, we end up with the estimate
\begin{equation}
\label{EstimateBSDE_Proof}
\sup_{s\in[t,T]}\E\big[|Y_s|^2\big] + \|Z\|_{\H^2(t,T)^d}^2 \ \leq \ \bar C\bigg(\E\big[|Y_T|^2\big] + \E\int_t^T |F(s,0,0)|^2 ds\bigg),
\end{equation}
for some positive constant $\bar C$ depending only on $T$ and $C_F$. In the \emph{second} step of the proof we estimate $\|Y\|_{\S^2(t,T)}^2=\E[\sup_{t\leq s\leq T}|Y_s|^2]$ proceeding as follows: we take the square in relation \eqref{E:BSDE_Appendix}, followed by the $\sup$ over $s$ and then the expectation. Finally, the claim follows exploiting the Lipschitz property of $F$ with respect to $(y,z)$, estimate \eqref{EstimateBSDE_Proof}, and Burkholder-Davis-Gundy inequality.
\ep

\setcounter{Theorem}{0}
\setcounter{equation}{0}

\subsection{Limit theorem for backward stochastic differential equations}

We prove a limit theorem for backward stochastic differential equations designed for our purposes, which is inspired by the monotonic limit theorem in \cite{peng00}, even if it is formulated under a different set of assumptions. In particular, the monotonicity of the sequence $(Y^n)_n$ is not assumed. On the other hand, we impose a uniform boundedness for the sequence $(Y^n)_n$ in $\S^p(t,T)$ for some $p>2$, instead of $p=2$ as in \cite{peng00}. Furthermore, unlike \cite{peng00}, the terminal condition and the generator function of the BSDE solved by $Y^n$ are allowed to vary with $n$.

\begin{Theorem}
\label{T:LimitThmBSDE}
Let $(F_n)_n$ be a sequence of generator functions satisfying assumption \textup{(Aa)-(Ac)}, with the same constants $C_F$ and $M_F$ for all $n$. For any $n$, let $(Y^n,Z^n,K^n)\in\S^2(t,T)\times\H^2(t,T)^d\times\A^{+,2}(t,T)$, with $Y^n$ and $K^n$ continuous, satisfying, $\P$-a.s.,
\[
Y_s^n \ = \ Y_T^n + \int_s^T F_n(r,Y_r^n,Z_r^n) dr + K_T^n - K_s^n - \int_s^T \langle Z_r^n, dW_r\rangle, \qquad t \leq s \leq T
\]
and
\[
\|Y^n\|_{\S^2(t,T)}^2 + \|Z^n\|_{\H^2(t,T)^d} + \|K^n\|_{\S^2(t,T)} \ \leq \ C, \qquad \forall\,n\in\N,
\]
for some positive constant $C$, independent of $n$. Suppose that there exist a generator function $F$ satisfying conditions \textup{(Aa)-(Ac)} and a continuous process $Y\in\S^2(t,T)$, in addition $\sup_n\|Y^n\|_{\S^p(t,T)}<\infty$ for some $p>2$, and, for some null measurable sets $N_F\subset[t,T]\times\Omega$ and $N_Y\subset\Omega$,
\begin{align*}
F_n(s,\omega,y,z) \ &\overset{n\rightarrow\infty}{\longrightarrow} \ F(s,\omega,y,z), \qquad \forall\,(s,\omega,y,z)\in(([t,T]\times\Omega)\backslash N_F)\times\R\times\R^d, \\
Y_s^n(\omega) \ &\overset{n\rightarrow\infty}{\longrightarrow} \ Y_s(\omega), \qquad\qquad\;\; \forall\,(s,\omega)\in[t,T]\times(\Omega\backslash N_Y).
\end{align*}
Then, there exists a unique pair $(Z,K)\in\H^2(t,T)^d\times\A^{+,2}(t,T)$ such that, $\P$-a.s.,
\begin{equation}
\label{E:BSDELimit}
Y_s \ = \ Y_T + \int_s^T F(r,Y_r,Z_r) dr + K_T - K_s - \int_s^T \langle Z_r, dW_r\rangle, \qquad t \leq s \leq T.
\end{equation}
In addition, $Z^n$ converges strongly $($resp. weakly$)$ to $Z$ in $\L^q(t,T;\R^d)$ $($resp. $\H^2(t,T)^d$$)$, for any $q\in[1,2[$, and $K_s^n$ converges weakly to $K_s$ in $L^2(\Omega,\Fc_s,\P)$, for any $s\in[t,T]$.
\end{Theorem}
\begin{Remark}
\label{R:Y_Sp<infty}
{\rm
Notice that, under the assumptions of Theorem \ref{T:LimitThmBSDE} (more precisely, given that $Y$ is continuous, $\sup_n\|Y^n\|_{\S^p(t,T)}<\infty$ for some $p>2$, $Y_s^n(\omega) \rightarrow Y_s(\omega)$ as $n$ tends to infinity for all $(s,\omega)\in[t,T]\times(\Omega\backslash N_Y)$), it follows that $\|Y\|_{\S^p(t,T)}<\infty$. Indeed, from Fatou's lemma we have
\begin{equation}
\label{E:Y_Sp<infty1}
\E\Big[\liminf_{n\rightarrow\infty}\sup_{t\leq s\leq T}|Y_s^n|^p\Big] \ \leq \ \liminf_{n\rightarrow\infty}\|Y^n\|_{\S^p(t,T)}^p \ < \ \infty.
\end{equation}
Moreover, since $Y$ is continuous, there exists a null measurable set $N_Y'\subset\Omega$ such that $s\mapsto Y_s(\omega)$ is continuous on $[t,T]$ for every $\omega\in\Omega\backslash N_Y'$. Then, for any $\omega\in\Omega\backslash(N_Y\cup N_Y')$, there exists $\tau(\omega)\in[t,T]$ such that
\begin{equation}
\label{E:Y_Sp<infty2}
\sup_{t \leq s \leq T}|Y_s(\omega)|^p \ = \ |Y_{\tau(\omega)}(\omega)|^p \ = \ \lim_{n\rightarrow\infty} |Y_{\tau(\omega)}^n(\omega)|^p \ \leq \ \liminf_{n\rightarrow\infty} \sup_{t\leq s\leq T}|Y_s^n(\omega)|^p.
\end{equation}
Therefore, combining \eqref{E:Y_Sp<infty1} with \eqref{E:Y_Sp<infty2}, we end up with $\|Y\|_{\S^p(t,T)}<\infty$.
\ep
}
\end{Remark}
\textbf{Proof.}
We begin proving the uniqueness of $(Z,K)$. Let $(Z,K),(Z',K')\in\H^2(t,T)^d\times\A^{+,2}(t,T)$ be two pairs satisfying \eqref{E:BSDELimit}. Taking the difference and rearranging the terms, we obtain
\[
\int_s^T \langle Z_r-Z_r',dW_r\rangle \ = \ \int_s^T \big(F(r,Y_r,Z_r)-F(r,Y_r,Z_r')\big) dr + K_T-K_s - (K_T'-K_s').
\]
Now, the right-hand side has finite variation, while the left-hand side has not finite variation, unless $Z=Z'$. This implies $Z=Z'$, from which we deduce $K=K'$.

The rest of the proof is devoted to the existence of $(Z,K)$ and it is divided in different steps.\\
\emph{Step 1. Limit BSDE.} From the assumptions, we see that there exists a positive constant $c$, independent of $n$, such that
\[
\E\int_t^T |F_n(r,Y_r^n,Z_r^n)|^2 dr \ \leq \ c, \qquad \forall\,n\in\N.
\]
It follows that the sequence $(Z_\cdot^n,F_n(\cdot,Y_\cdot^n,Z_\cdot^n))_n$ is bounded in the Hilbert space $\H^2(t,T)^d\times\L^2(t,T;\R)$. Therefore, there exists a subsequence $(Z_\cdot^{n_k},F_{n_k}(\cdot,Y_\cdot^{n_k},Z_\cdot^{n_k}))_k$ which converges weakly to some $(Z,G)\in\H^2(t,T)^d\times\L^2(t,T;\R)$. This implies that, for any $s\in[t,T]$, the following weak convergences hold in $L^2(\Omega,\Fc_s,\P)$ as $k\rightarrow\infty$:
\[
\int_t^s F_{n_k}(r,Y_r^{n_k},Z_r^{n_k}) dr \ \rightharpoonup \ \int_t^s G(r) dr, \qquad\qquad
\int_t^s \langle Z_r^{n_k},dW_r\rangle \ \rightharpoonup \ \int_t^s \langle Z_r,dW_r\rangle.
\]
Since
\[
K_s^n \ = \ Y_t^n - Y_s^n - \int_t^s F_n(r,Y_r^n,Z_r^n) dr + \int_t^s \langle Z_r^n,dW_r\rangle
\]
and, by assumption, $Y_s^n \rightarrow Y_s$ strongly in $L^2(\Omega,\Fc_s,\P)$, we also have the weak convergence, as $k\rightarrow\infty$,
\begin{equation}
\label{E:K_tau=Y}
K_s^{n_k} \ \rightharpoonup \ K_s,
\end{equation}
where
\[
K_s \ := \ Y_t - Y_s - \int_t^s G(r) dr + \int_t^s \langle Z_r,dW_r\rangle, \qquad t\leq s\leq T.
\]
Notice that $(K_s)_{t\leq s\leq T}$ is adapted and continuous, so that it is a predictable process. We have that $\E[|K_T|^2] < \infty$. Let us prove that $K$ is a nondecreasing process. For any pair $r,s$ with $t \leq r \leq s \leq T$, we have $K_r \leq K_s$, $\P$-a.s.. Indeed, let ${\xi}\in L^2(\Omega,\Fc_s,\P)$ be nonnegative, then, from the martingale representation theorem, we see that there exist a random variable $\zeta\in L^2(\Omega,\Fc_r,\P)$ and an $\F$-predictable square integrable process $\eta$ such that
\[
{\xi} \ = \ \zeta + \int_r^s \eta_u dW_u.
\]
Therefore
\begin{align*}
0 \ &\leq \ \E[{\xi}(K_s^n-K_r^n)] \ = \ \E[{\xi} K_s^n] - \E[\zeta K_r^n] - \E\bigg[\E\bigg[K_r^n\int_r^s \eta_u dW_u \bigg|\Fc_r\bigg]\bigg] \\
&= \ \E[{\xi} K_s^n] - \E[\zeta K_r^n] \ \overset{n\rightarrow\infty}{\longrightarrow} \ \E[{\xi} K_s] - \E[\zeta K_r] \ = \ \E[{\xi}(K_s-K_r)],
\end{align*}
which shows that $K_r \leq K_s$, $\P$-a.s.. As a consequence, there exists a null measurable set $N\subset\Omega$ such that $K_r(\omega) \leq K_s(\omega)$, for all $\omega\in\Omega\backslash N$, with $r,s\in\Q\cap[0,T]$, $r<s$. Then, from the continuity of $K$ it follows that it is a nondecreasing process, so that $K\in\A^{+,2}(t,T)$.

Finally, we notice that the process $Z$ in expression \eqref{E:K_tau=Y} is uniquely determined, as it can be seen identifying the Brownian parts and the finite variation parts in \eqref{E:K_tau=Y}. Thus, not only the subsequence $(Z^{n_k})_k$, but all the sequence $(Z^n)_n$ converges weakly to $Z$ in $\H^2(t,T)^d$. It remains to show that $G(r)$ in \eqref{E:K_tau=Y} is actually $F(r,Y_r,Z_r)$.\\
\emph{Step 2. Strong convergence of $(Z^n)_n$.} Let $\alpha\in(0,1)$ and consider the function $h_\alpha(y) = |\min(y-\alpha,0)|^2$, $y\in\R$. By applying Meyer-It\^o's formula combined with the occupation times formula (see, e.g., Theorem 70 and Corollary 1, Chapter IV, in \cite{protter05}) to $h_\alpha(Y_s^n-Y_s)$ between $t$ and $T$, observing that the second derivative of $h_\alpha$ in the sense of distributions is a $\sigma$-finite Borel measure on $\R$ absolutely continuous to the Lebesgue measure with density $2\cdot1_{]-\infty,\alpha[}(\cdot)$, we obtain
\begin{align*}
&\E\big[|\min(Y_t^n-Y_t-\alpha,0)|^2\big] + \E\int_t^T 1_{\{Y_s^n-Y_s<\alpha\}} |Z_s^n-Z_s|^2 ds \\
&= \ \E\big[|\min(Y_T^n-Y_T-\alpha,0)|^2\big] + 2\E\int_t^T \min(Y_s^n-Y_s-\alpha,0) \big(F_n(s,Y_s^n,Z_s^n) - G(s)\big)ds \\
&\quad \ + 2\E\int_t^T \min(Y_s^n - Y_s-\alpha,0) dK_s^n - 2\E\int_t^T \min(Y_s^n - Y_s-\alpha,0) dK_s.
\end{align*}
Since $\min(Y_s^n - Y_s-\alpha,0) dK_s^n \leq 0$, we get
\begin{align}
\label{E:min<K}
&\E\int_t^T 1_{\{Y_s^n-Y_s<\alpha\}} |Z_s^n-Z_s|^2 ds \ \leq \ \E\big[|\min(Y_T^n-Y_T-\alpha,0)|^2\big] \\
&+ 2\E\int_t^T \min(Y_s^n-Y_s-\alpha,0) \big(F_n(s,Y_s^n,Z_s^n) - G(s)\big)ds - 2\E\int_t^T \min(Y_s^n - Y_s-\alpha,0) dK_s. \notag
\end{align}
Let us study the behavior of the right-hand side of \eqref{E:min<K} as $n$ goes to infinity. We begin noting that
\begin{equation}
\label{E:ProofConvergence1}
\E\big[|\min(Y_T^n-Y_T-\alpha,0)|^2\big] \ \overset{n\rightarrow\infty}{\longrightarrow} \ \alpha^2.
\end{equation}
Regarding the second-term on the right-hand side of \eqref{E:min<K}, since the sequence $(F_n(\cdot,Y_\cdot^n,Z_\cdot^n) - G(\cdot))_n$ is bounded in $\L^2(t,T;\R)$, we have
\[
\sup_{n\in\N}\bigg(\E\bigg[\int_t^T|F_n(s,Y_s^n,Z_s^n) - G(s)|^2ds\bigg]\bigg)^{\frac{1}{2}} \ =: \ \bar c \ < \ \infty.
\]
Therefore, by Cauchy-Schwarz inequality we find
\begin{align}
\label{E:ProofConvergence2}
&\E\int_t^T |\min(Y_s^n-Y_s-\alpha,0)| |F_n(s,Y_s^n,Z_s^n) - G(s)|ds \notag \\
&\leq \bar c \bigg(\E\bigg[\int_t^T|\min(Y_s^n-Y_s-\alpha,0)|^2ds\bigg]\bigg)^{\frac{1}{2}} \  \overset{n\rightarrow\infty}{\longrightarrow} \ \bar c\sqrt{T-t}\,\alpha.
\end{align}
Concerning the last term on the right-hand side of \eqref{E:min<K}, we notice that, by assumption and Remark \ref{R:Y_Sp<infty}, there exists some $p>2$ such that, from Cauchy-Schwarz inequality,
\begin{align*}
&\sup_{n\in\N}\E\bigg[\int_t^T |\min(Y_s^n - Y_s-\alpha,0)|^{\frac{p}{2}} dK_s\bigg] \\
&\leq \ \sup_{n\in\N}\bigg(\E\Big[\sup_{t \leq s \leq T} |\min(Y_s^n - Y_s-\alpha,0)|^p \Big]\bigg)^{\frac{1}{2}} \big(\E\big[|K_T|^2\big]\big)^{\frac{1}{2}} \ < \infty.
\end{align*}
It follows that $(\min(Y_\cdot^n - Y_\cdot-\alpha,0))_n$ is a uniformly integrable sequence on $([t,T]\times\Omega,\Bc([t,T])\otimes\Fc,dK_s\otimes d\P)$. Moreover, by assumption, there exists a null measurable set $N_Y\subset\Omega$ such that $Y_s^n(\omega)$ converges to $Y_s(\omega)$, for any $(s,\omega)\notin[t,T]\times N_Y$. Notice that $dK_s\otimes d\P([t,T]\times N_Y)=0$, therefore $Y^n$ converges to $Y$ pointwise a.e. with respect to $dK_s\otimes d\P$. This implies that
\begin{equation}
\label{E:ProofConvergence3}
\E\bigg[\int_t^T |\min(Y_s^n - Y_s-\alpha,0)| dK_s\bigg] \ \overset{n\rightarrow\infty}{\longrightarrow} \ \alpha\E[K_T].
\end{equation}
By the convergence results \eqref{E:ProofConvergence1}, \eqref{E:ProofConvergence2}, and \eqref{E:ProofConvergence3}, \eqref{E:min<K} gives
\begin{equation}
\label{E:min<K2}
\limsup_{n\rightarrow\infty}\E\int_t^T 1_{\{Y_s^n-Y_s<\alpha\}} |Z_s^n-Z_s|^2 ds \ \leq \ \alpha^2 + 2\bar c\sqrt{T-t}\,\alpha + 2\alpha\E[K_T].
\end{equation}
From Egoroff's theorem, for any $\delta>0$ there exists a measurable set $A\subset[t,T]\times\Omega$, with $ds\otimes d\P(A)<\delta$, such that $(Y^n)_n$ converges uniformly to $Y$ on $([t,T]\times\Omega)\backslash A$. In particular, for any $\alpha\in]0,1[$ we have $|Y_s^n(\omega)-Y_s(\omega)| < \alpha$, for all $(s,\omega)\in([t,T]\times\Omega)\backslash A$, whenever $n$ is large enough. Therefore, from \eqref{E:min<K2} we get
\begin{align*}
&\limsup_{n\rightarrow\infty}\E\int_t^T 1_{([t,T]\times\Omega)\backslash A} |Z_s^n-Z_s|^2 ds \ = \ \limsup_{n\rightarrow\infty}\E\int_t^T 1_{([t,T]\times\Omega)\backslash A} 1_{\{Y_s^n-Y_s<\alpha\}} |Z_s^n-Z_s|^2 ds \notag \\
&\leq \ \limsup_{n\rightarrow\infty}\E\int_t^T 1_{\{Y_s^n-Y_s<\alpha\}} |Z_s^n-Z_s|^2 ds \ \leq \ \alpha^2 + 2\bar c\sqrt{T-t}\,\alpha + 2\alpha\E[K_T].
\end{align*}
Sending $\alpha\rightarrow0^+$, we obtain
\begin{equation}
\label{E:alpha-->0+}
\lim_{n\rightarrow\infty}\E\int_t^T 1_{([t,T]\times\Omega)\backslash A} |Z_s^n-Z_s|^2 ds \ = \ 0.
\end{equation}
Now, let $q\in[1,2[$; by H\"older's inequality,
\begin{align*}
&\E\int_t^T |Z_s^n-Z_s|^q ds \ = \ \E\int_t^T 1_{([t,T]\times\Omega)\backslash A} |Z_s^n-Z_s|^q ds + \E\int_t^T 1_A |Z_s^n-Z_s|^q ds \\
&\leq \ \bigg(\E\int_t^T 1_{([t,T]\times\Omega)\backslash A} |Z_s^n-Z_s|^2 ds\bigg)^{\frac{q}{2}}(T-t)^{\frac{2-q}{2}} + \bigg(\E\int_t^T |Z_s^n-Z_s|^2 ds\bigg)^{\frac{q}{2}}\delta^{\frac{2-q}{2}}.
\end{align*}
Since the sequence $(Z^n)_n$ is bounded in $\H^2(t,T)^d$, we have
\[
\sup_{n\in\N}\E\int_t^T |Z_s^n-Z_s|^2 ds \ =: \ \hat c < \infty.
\]
Therefore
\[
\E\int_t^T |Z_s^n-Z_s|^q ds \ \leq \ \bigg(\E\int_t^T 1_{([t,T]\times\Omega)\backslash A} |Z_s^n-Z_s|^2 ds\bigg)^{\frac{q}{2}}(T-t)^{\frac{2-q}{2}} + \hat c^{\frac{q}{2}}\delta^{\frac{2-q}{2}},
\]
which implies, by \eqref{E:alpha-->0+},
\[
\limsup_{n\rightarrow\infty}\E\int_t^T |Z_s^n-Z_s|^q ds \ \leq \ \hat c^{\frac{q}{2}}\delta^{\frac{2-q}{2}}.
\]
Sending $\delta\rightarrow0^+$ we deduce the strong convergence of $Z^n$ towards $Z$ in $\L^q(t,T;\R^d)$, for any $q\in[1,2[$.

Notice that, for any $q\in[1,2[$, we have (recalling the standard inequality $(x+y)^q \leq 2^{q-1}(x^q+y^q)$, for any $x,y\in\R_+$)
\begin{align*}
\E\bigg[\int_t^T |F_n(s,Y_s^n,Z_s^n) - F(s,Y_s,Z_s)|^q ds\bigg] \ &\leq \ 2^{q-1}\E\bigg[\int_t^T |F_n(s,Y_s^n,Z_s^n) - F_n(s,Y_s,Z_s)|^q ds\bigg] \\
&\;\;\; + 2^{q-1}\E\bigg[\int_t^T |F_n(s,Y_s,Z_s) - F(s,Y_s,Z_s)|^q ds\bigg].
\end{align*}
Therefore, by the uniform Lipschitz condition on $F_n$ with respect to $(y,z)$, and the convergence of $F_n$ towards $F$, we deduce the strong convergence of $(F_n(\cdot,Y_\cdot^n,Z_\cdot^n))_n$ to $F(\cdot,Y_\cdot,Z_\cdot)$ in $\L^q(t,T;\R)$, $q\in[1,2[$. Since $G(\cdot)$ is the weak limit of $(F_n(\cdot,Y_\cdot^n,Z_\cdot^n))_n$ in $\L^2(t,T;\R)$, we deduce that $G(\cdot)$ $=$ $F(\cdot,Y_\cdot,Z_\cdot)$. In conclusion, the triple $(Y,Z,K)$ solves the backward stochastic differential equation \eqref{E:BSDELimit}.
\ep

\setcounter{Theorem}{0}
\setcounter{equation}{0}

\subsection{An additional result in real analysis}

\begin{Lemma}
\label{L:StabilityApp}
Let $(G_{n,k})_{n,k\in\N}$, $(G_n)_{n\in\N}$, and $G$ be $\R^q$-valued continuous functions on $[0,T]\times X$, where $(X,d)$ is a separable metric space, and
\[
G_{n,k}(t,x) \ \overset{k\rightarrow\infty}{\longrightarrow} \ G_n(t,x), \qquad G_n(t,x) \ \overset{n\rightarrow\infty}{\longrightarrow} \ G(t,x), \qquad \forall\,(t,x)\in[0,T]\times X.
\]
Moreover, $G_{n,k}(t,x)\rightarrow G_n(t,x)$ as $k\rightarrow\infty$, for all $x\in X$, uniformly with respect to $t\in[0,T]$. Suppose also that the functions $G_{n,k}(t,\cdot)$, $n,k\in\N$, are equicontinuous on compact sets, uniformly with respect to $t\in[0,T]$. Then, there exists a subsequence $(G_{n,k_n})_{n\in\N}$ which converges pointwise to $G$ on $[0,T]\times X$.
\end{Lemma}
\textbf{Proof.}
We begin noting that, as a direct consequence of the assumptions of the lemma, the functions $G(t,\cdot)$, $G_n(t,\cdot)$, and $G_{n,k}(t,\cdot)$, for all $n,k\in\N$, are equicontinuous on compact sets, uniformly with respect to $t\in[0,T]$.

Let $D=\{x_1,x_2,\ldots,x_j,\ldots\}$ be a countable dense subset of $X$. Fix $n\in\N\backslash\{0\}$. Then, for any $j\in\N$ there exists $k_{n,j}\in\N$ such that
\[
|G_{n,k}(t,x_j) - G_n(t,x_j)| \ \leq \ \frac{1}{n}, \qquad \forall\,k\geq k_{n,j},\,\forall\,t\in[0,T].
\]
Set $k_n:=k_{n-1}\vee k_{n,1}\vee\cdots\vee k_{n,n}$, $\forall\,n\in\N$, where $k_{-1}:=0$. Then, we have
\[
|G_{n,k_n}(t,x_j) - G(t,x_j)| \ \overset{n\rightarrow\infty}{\longrightarrow} \ 0, \qquad \forall\,j\in\N,
\]
for all $t\in[0,T]$. It remains to prove that the convergence holds for all $(t,x)\in[0,T]\times X$. To this end, fix $x\in X$ and consider a subsequence $(x_{j_m})_{m\in\N}\subset D$ which converges to $x$. Then, the set $K$ defined by
\[
K \ := \ (x_{j_m})_{m\in\N}\cup\{x\}
\]
is a compact subset of $X$. Recall that the functions $G(t,\cdot)$ and $G_{n,k_n}(t,\cdot)$, for all $n\in\N$, are equicontinuous on $K$, uniformly with respect to $t\in[0,T]$. Then, for every $\eps>0$, there exists $\delta>0$ such that, for all $n\in\N$,
\[
|G_{n,k_n}(t,x_1) - G_{n,k_n}(t,x_2)| \ \leq \ \frac{\eps}{3}, \qquad |G(t,x_1) - G(t,x_2)| \ \leq \ \frac{\eps}{3},
\]
whenever $\|x_1-x_2\|\leq\delta$, $x_1,x_2\in K$, for all $t\in[0,T]$. Fix $t\in[0,T]$ and $x_{j_{m_0}}\in(x_{j_m})_{m\in\N}$ such that $\|x-x_{j_{m_0}}\|\leq\delta$. Then, we can find $n_0\in\N$ (possibly depending on $t$) for which $|G_{n,k_n}(t,x_{j_{m_0}})-G(t,x_{j_{m_0}})|\leq\eps/3$ for any $n\geq n_0$. Therefore, given $n\geq n_0$ we obtain
\begin{align*}
|G_{n,k_n}(t,x) - G(t,x)| \ &\leq \ |G_{n,k_n}(t,x) - G_{n,k_n}(t,x_{j_{m_0}})| + |G_{n,k_n}(t,x_{j_{m_0}}) - G(t,x_{j_{m_0}})| \\
&\quad \ + |G(t,x_{j_{m_0}}) - G(t,x)| \ \leq \ \eps.
\end{align*}
This implies that $G_{n,k_n}$ converges to $G$ at $(t,x)$, and the claim follows from the arbitrariness of $(t,x)$.
\ep

\vspace{5mm}

\noindent\textbf{Acknowledgments.}
The authors are very grateful to the Referee for her/his 
very stimulating remarks which helped in 
drastically improve the first version of the paper.\\ 
 Part of the paper was done during the visit
of the second named author at the Centre for Advanced Study (CAS) at the Norwegian Academy of Science and Letters in Oslo.

\small
\bibliographystyle{plain}
\bibliography{biblio}

\end{document}